\input amstex


\def\b1{\text{\bf 1}}

\def\BP{{\Bbb P}}
\def\BQ{{\Bbb Q}}
\def\bR{\bar{R}}

\def\bx{\bar{x}}
\def\BZ{{\Bbb Z}}
\def\CA{{\Cal A}}
\def\CB{{\Cal B}}
\def\CC{{\Cal C}}
\def\CD{{\Cal D}}

\def\CK{{\Cal K}}
\def\CL{{\Cal L}}
\def\CM{{\Cal M}}

\def\CO{{\Cal O}}
\def\CP{{\Cal P}}

\def\CT{{\Cal T}}
\def\CU{{\Cal U}}
\def\CV{{\Cal V}}

\def\CX{{\Cal X}}

\def\dpar{\partial}
\def\dplus{\buildrel\cdot\over{+}}

\def\fg{{\frak g}}

\def\Id{\text{Id}}

\def\tAlg{\widetilde{{\Cal A}lg}}

\def\tc{\tilde c}

\def\tg{\tilde g}

\def\tJ{\tilde J}

\def\tOmega{\tilde{\Omega}}

\def\tx{\tilde x}


\def\btu{\bigtriangleup}
\def\hra{\hookrightarrow}
\def\iso{\buildrel\sim\over\longrightarrow} 

\def\lra{\longrightarrow}

\parskip=6pt

\documentstyle{amsppt}
\document
\NoBlackBoxes


\centerline{\bf Gerbes of chiral differential operators. II}

\bigskip
\centerline{Vassily Gorbounov, Fyodor Malikov, Vadim Schechtman}
\bigskip


\bigskip\bigskip 

\centerline{\bf Contents}

\bigskip\bigskip

Introduction

\S 0. Recollections and notation

\S 1. Vertex algebroids

\S 2. From vertex algebras to vertex algebroids 

\S 3. Category of vertex algebroids 

\S 4. Cofibered structure 

\S 5. Chern-Simons 

\S 6. Atiyah 

\S 7. Gerbes of vertex algebroids 

\S 8. Vertex envelope of a conformal algebra 

\S 9. Enveloping algebra of a vertex algebroid 

References

\bigskip\bigskip


\centerline{\bf Introdution}

\bigskip\bigskip

The aim of this paper is a study of certain class 
of chiral (vertex) algebras which we call 
{\it chiral algebras of differential operators}. 

{\bf I.0.} To explain this notion, let us take up a "naive" 
approach to vertex algebras. Let us fix (in this Introduction) 
a ground field $k$ of characteristic $0$. Recall 
(cf. [K]) that a {\it $\BZ_{\geq 0}$-graded vertex algebra} 
is a $\BZ_{\geq 0}$-graded $k$-vector space $V=\oplus_{i\geq 0}\ V_i$ 
equipped with a distinguished {\it vacuum vector} $\b1\in V_0$ and 
a family of bilinear operations
$$
_{(n)}:\ V\otimes_k V\lra V,\ \ x\otimes y\mapsto x_{(n)}y, 
$$
$(n\in\BZ)$, the operation $_{(n)}$ having the degree $-n-1$, that is, 
$V_{i(n)}V_j\subset V_{i+j-n-1}$ where by definition $V_i=0$ for $i<0$.   
These operations must satisfy a family of quadratic 
{\it Borcherds identities}. Throughout this paper we shall deal only with 
$\BZ_{\geq 0}$-graded vertex algebras, so we shall omit the words 
"$\BZ_{\geq 0}$-graded" when speaking about them. 

Note that in particular the operation $_{(-1)}$ has degree $0$. 
It is neither commutative nor associative in general. 
 
We are interested in a question: what structure on the subspace 
$V_{\leq 1}:=V_0\oplus V_1$ is induced by the structure of a vertex 
algebra on $V$? 

The claims in italics below are easily proved using the Borcherds 
identities, see the main body of the paper.  

First of all, {\it $V_0$ is a commutative associative $k$-algebra with 
unit $\b1$ with respect to $_{(-1)}$}. Let us denote this algebra by 
$A$, and the operation $a_{(-1)}b$ on it by $ab$. Note that for each 
$i$ the operation $_{(-1)}$ gives a mapping 
$$
A\otimes V_i\lra V_i,
\eqno{(I1)} 
$$ 
to be denoted $ax$,  
but this mapping does {\it not} make $V_i$ an $A$-module since 
the associativity $abx=(ab)x$ is not satisfied in general.  

We have a mapping $\dpar:\ V_i\lra V_{i+1}$ defined by $\dpar x=x_{(-2)}\b1$. 
Let us denote by $\Omega\subset V_1$ a $k$-subspace generated by all 
elements $a\dpar b,\ a,b\in A$. {\it The operation} (I1) {\it makes $\Omega$ an 
$A$-module (i.e. the associativity holds true on this subspace) and $\dpar:\ 
A\lra\Omega$ is a derivation}. 

Let $T$ be the quotient space $V_1/\Omega$. {\it The operation} 
(I1) {\it induces the structure of an $A$-module on $T$ and operation 
$_{(0)}$ induces a structure of a Lie algebra on $T$ and an action of 
$T$ on $A$ by derivations. These structures are compatible: 
they make $T$ a Lie algebroid} (see 0.2) {\it over $A$.} 

{\it The operation $_{(0)}$ induces an action of $T$ on $\Omega$ which 
makes $\Omega$ a module over the Lie algebra $T$, and $\dpar$ is a 
morphism of $T$-modules.} 

{\it The operation $_{(1)}$ induces an $A$-bilinear pairing 
$T\otimes\Omega\lra A$ which is a morphism of $T$-modules.} 

Let us choose a splitting 
$$
s:\ T\lra V_1
\eqno{(I2)}
$$ 
of the canonical projection 
$\pi:\ V_1\lra T$. Define a mapping $\gamma:\ A\times T\lra\Omega$ by 
$\gamma(a,\tau)=s(a\tau)-as(\tau)$. Define a symmetric mapping 
$\langle,\rangle:\ T\times T\lra A$ by $\langle\tau,\tau'\rangle=
s(\tau)_{(1)}s(\tau')$. Finally, define a skew symmetric operation 
$c:\ T\times T\lra \Omega$ by $c(\tau,\tau')=s([\tau,\tau'])-
[s(\tau),s(\tau')]$ where we set $[x,y]:=\frac{1}{2}(x_{(0)}y-y_{(0)}x)$. 

These three operations, $\gamma,\ \langle,\rangle$ and $c$, satisfy 
certain compatibilities, listed below in 1.3, (A1) --- (A5). 
 
Thus, we have assigned to a  vertex algebra $V$ 
(with a splitting (I2)) a collection of data $\CA=(A,T,\Omega,\dpar,
\gamma,\langle,\rangle,c)$ as above satisfying properties (A1) --- (A5). 
We call such collection of data a {\bf vertex algebroid} which is a central 
hero of our story. 

The above procedure gives a functor form the category of  
vertex algebras to the (appropriately defined) 
category of vertex algebroids. This functor admits a left adjoint $U$, 
called {\it vertex envelope}. A vertex algebra of the form $U\CA$ 
is called a {\it chiral algebra of differential operators}. 

{\bf I.2.} Let $X$ be a smooth algebraic variety over $k$, $U\subset X$ 
an affine Zariski open subset, $A=\CO(U)$. Let us consider the class of vertex 
algebroids of the form $\CA=(A,T,\Omega,\dpar,\ldots)$ with $T=Der_k(A),\ 
\Omega=\Omega^1_k(A),\ \dpar=d_{DR}$ --- the usual de Rham differential, 
and with the standard action of $T$ on $\Omega$ by the Lie derivative. 

Such vertex algebroids (with suitably defined morphisms) form a 
{\it groupoid} (i.e. a category where all morphisms are invertible) 
$\CA lg(A)$. This category is nonempty if there exists an $A$-base of $T$ 
consisting of commuting vector fields.  

This groupoid has a remarkable additional structure. 
To describe this structure, 
let us define another groupoid $\Omega^{[2,3\rangle}(A)$ as follows. 
By definition, its objects will be all closed $3$-forms on $A$. 
If $\omega,\omega'$ are two such $3$-forms, the set of morphisms 
between them is defined as 
$$
Hom_{\Omega^{[2,3\rangle}(A)}(\omega,\omega')=
\{\eta\in\Omega^2_{A/k}|\ d_{DR}\eta=\omega'-\omega\}
$$
The compostion is defined in the obvious way. The addition of $3$-forms 
induces a functor
$$
\Omega^{[2,3\rangle}(A)\times\Omega^{[2,3\rangle}(A)\lra 
\Omega^{[2,3\rangle}(A)
$$
which makes $\Omega^{[2,3\rangle}(A)$ an {\it Abelian group in categories}. 

The first main result of this paper is  

{\bf Theorem I1.} {\it The groupoid $\CA lg(A)$ admits a canonical structure 
of a $\Omega^{[2,3\rangle}(A)$-Torseur.} 

This is Theorem 7.2 of the paper. It means that one has a functor 
$$
\dplus:\ \CA lg(A)\times\Omega^{[2,3\rangle}(A)\lra\CA lg(A),\ 
(\CA,\omega)\mapsto \CA\dplus\omega
$$
which defines an Action of the Abelian group $\Omega^{[2,3\rangle}(A)$ 
on $\CA lg(A)$, such that for each $\CA\in\CO b\CA lg(A)$ the induced 
functor 
$$
\CA\dplus ?:\ \Omega^{[2,3\rangle}(A)\lra\CA lg(A)
$$
is an equivalence of categories. 

{\bf I.3.} Sheafifying the previous construction one gets a sheaf 
of groupoids, i.e. {\it gerbe} $\CA lg_X$ over $X$, with the {\it lien} 
$\Omega^{[2,3\rangle}_X:=\Omega^2_X\lra\Omega^{3,closed}$ (a short 
complex of sheaves over $X$, with $\Omega^2_X$ sitting in degree $0$). 

By a simple general homological construction (see 7.3), this gerbe 
gives rise to a characteristic class
$$
c(\CA lg_X)\in H^2(X;\Omega^{[2,3\rangle}_X)
$$
with the following property. 

{\it The groupoid of global sections $\CA lg_X(X)$ is nonempty iff 
$c(\CA lg_X)=0$.
If this is the case then its set of connected components is a nonempty  
$H^1(X;\Omega^{[2,3\rangle}_X)$-torseur and the group 
of automorphisms of its object is isomorphic to 
$H^0(X;\Omega^{[2,3\rangle}_X)$.} 

Our next aim is to calculate the class $c(\CA lg_X)$. One shows 
(see Theorem 7.10) that it essentially coincides with the second component 
of the Chern character of the tangent bundle $\CT_X$. Let us 
describe an "explicit formula" for it. 

Let $E$ be an arbitrary finite dimensional vector bundle over $X$, 
given by a Cech cocycle $(g_{ij})\in Z^1(\CU;GL_r(\CO_X))$ on 
a Zariski open covering $\CU=\{U_i\}$ of $X$. Define a Cech two-cocycle 
$ch_2(E)\in Z^2(\CU;\Omega^{[2,3\rangle}_X)$ by 
$$
ch_2(E)=\bigl(\frac{1}{2}tr(g_{ij}^{-1}g_{jk}^{-1}dg_{jk}dg_{ij}), 
\frac{1}{6}tr((g_{ij}^{-1}dg_{ij})^3)\bigr)
\eqno{(ACS)}
$$
This class may be called the {\bf "Atiyah-Chern-Simons class"} of $E$. 
Its first component, $\alpha(E)$, is an element of 
$Z^2(\CU;\Omega^2_X)$ 
which is the usual "Atiyah" representative of the degree $2$ part 
of the Chern character of $E$ "style Hodge" living in $H^2(X;\Omega^2_X)$,  
while the second component, $\beta(E)$ ("Chern-Simons"), is a Cech 
$1$-cochain with coefficients in $\Omega^{3,closed}_X$. 
The de Rham differential of $\alpha(E)$ is equal to the Cech coboundary 
of $\beta(E)$. One checks that the cohomology class of (ACS) 
does not depend on the choice of trivialization of $E$. 

Our second main result (Theorem 7.5) is 

{\bf Theorem I.2.} {\it The class $c(\CA lg)$ is equal to (the 
cohomology class of) $2ch_2(\CT_X)$.}

I.4. Our third main topic is an explicit construction of the enveloping 
algebra $U\CA$ of a vertex algebroid $\CA$ and "Poincar\'e-Birkhoff-Witt" 
type theorem for it. Let us formulate the last theorem (see Theorem 9.18). 

{\bf Theorem I.3.} {\it Each sheaf of algebras of chiral differential 
operators $\CD=U\CA,\ \CA\in\CA lg_X(X)$ admits a canonical filtration $F$, 
compatible with the conformal grading and finite on each component 
of fixed conformal weight, such that the corresponding graded algebra 
(which is a sheaf of $\BZ_{\geq 0}$-graded commutative vertex algebras) is 
$$
gr_F(\CD)=Sym_{\CO_X}\bigl\{(\oplus_{i\geq 1}\CT^{(i)})\oplus
(\oplus_{i\geq 1}\Omega^{(i)})\bigr\}
$$
where $\CT^{(i)}$ (resp. $\Omega^{(i)}$) is a copy of $\CT_X$ (resp. 
of $\Omega^1_X$) living in conformal weight $i$.} 

As a preparation to this theorem, we study in Section 8 the {\it vertex 
envelope} $UC$ of an arbitrary conformal algebra $C$. Although its 
construction  
is more or less contained in Kac's book [K], we present here in a sense 
more direct construction of $C$, maybe of independent interest. 
Note in particular an amusing Lie bracket (8.21.1) which is defined 
on an arbitrary conformal algebra. 

In the hope of possible arithmetical applications, we adopted in this 
paper a "characteristic free" approach to vertex algebras. Contrary to this 
Introduction (and to the tradition), in the main 
body of this work the ground ring will be an arbitrary commutative 
ring containing $1/2$ (with one exception: the PBW theorem 9.18). 
This generalization is achieved without great effort. Note that 
the original Borcherds' definition of vertex algebras was characteristic-free.        

I.4. This paper may be regarded as an "algebraic" version of our last 
note [GMS] where we worked in the analytical category. However, 
the approach adopted here is quite different from {\it op. cit.} 
(cf. also [MSV] and [MS]). 

During the preparation of this work we have enourmously benefited 
from the discussions and correspondence with Sasha Beilinson. He was 
the first to suggest (by analogy with the classical picture, [BB]) that 
in the algebraic category the situation is more subtle than in the 
analytical one. 
He has sent us a note [BD1] where similar questions 
are treated from a different point of view. In a sense, a great part of 
this work 
is a result of the attempts to understand this note (which still 
remains mysterious for us). For the details on Beilinson-Drinfeld approach, 
see [BD2], 3.8. 

Another special gratitude goes to H\'el\`ene Esnault. She has found a mistake 
in the previous formulation of our main result whose correction 
lead to the discovery of the class (ACS), and greatly 
helped in some computations. Theorem 7.10 was obtained in collaboraion 
with her.  

This work was mostly done while V.S. visited IHES. He is grateful 
to this institute for the support and excellent working conditions.

\bigskip\bigskip

\newpage


\centerline{\bf \S\ 0. Recollections and Notation}

\bigskip 

{\bf 0.0.} Throughout this paper, a {\it commutative ring (algebra)} 
will mean a commutative associative ring (algebra) with unit. 

$k$ will denote a fixed ground commutative ring. We will often assume 
that $1/2\in k$; we will indicate this assumption when necessary. 
 
{\it Algebra} will mean a $k$-algebra; 
$\otimes$ will mean the tensor product over $k$. 

In a nonassociatve algebra, $abc\ldots$ will mean $a(b(c(\ldots))\ldots)$. 

$\BZ_{\geq 0}$ will denote the set of nonnegative integers. The binomial 
coefficients are defined by 
$$
\binom{a}{b}=\frac{a(a-1)\cdot\ldots (a-b+1)}{b!},\ \ a\in\BZ,\ b\in\BZ_{\geq 
0}, 
\eqno{(0.0.1)}
$$
and $0$ if $b<0$. We have 
$$
\binom{-a-1}{b}=(-1)^b\binom{a+b}{b}
\eqno{(0.0.2)}
$$
for all $a\geq 0$.    

{\bf 0.1.} Let $A$ be a commutative algebra. The Lie algebra of 
$k$-derivations $T_{A}:=Der_k(A,A)$ acts on the $A$-module of K\"ahler 
$1$-differentials 
$\Omega^1_A:=\Omega^1_{A/k}$ according to the usual rule 
$$
\tau(a db)=\tau(a) db + a d\tau(b)
\eqno{(0.1.1)}
$$
The de Rham differential $d:\ A\lra\Omega^1_A$ commutes with the $T_A$-action. 

We have the canonical $A$-bilinear pairing 
$$
\langle\ ,\ \rangle:\ T_A\times\Omega^1_A\lra A,\ \ \langle\tau, a db\rangle= 
a\tau(b)
\eqno{(0.1.2)}
$$
We have 
$$
\tau(a\omega)=\tau(a)\omega+a\tau(\omega)
\eqno{(0.1.3)}
$$
$$
(a\tau)(\omega)=a\tau(\omega)+\langle\tau,\omega\rangle da
\eqno{(0.1.4)}
$$
$$
\tau(\langle\nu,\omega\rangle)=\langle [\tau,\nu],\omega\rangle+ 
\langle\nu,\tau(\omega)\rangle
\eqno{(0.1.5)}
$$
$$
\langle\tau,da\rangle=\tau(a)
\eqno{(0.1.6)}
$$
$(a\in A, \tau,\nu\in T_A, \omega\in\Omega^1_A)$. (Of course the last 
formula is a particular case of (0.1.2).) 

{\bf 0.2.} A {\it Lie $A$-algebroid} is a Lie algebra $T$ acting 
by derivations 
on $A$ and equipped with a structure of an $A$-module, such that 
$$
[\tau, a\nu]=a[\tau,\nu]+\tau(a)\nu
\eqno{(0.2.1)}
$$
and 
$$
(a\tau)(b)=a\tau(b)
\eqno{(0.2.2)}
$$
for all $\tau,\nu\in T;\ a,b\in A$. 

{\bf 0.2.1.} {\it Direct image (pushout).} Let $B$ be a commutative 
$A$-algebra, $i: A\lra B$ be the structure morphism and $T$ be an $A$-algebroid 
Lie. Assume that $T$, as a Lie algebra, acts on $B$ by derivations in such 
a way that $i$ is a morphism of $T$-modules and 
$$
(a\tau)(b)=a\tau(b)
$$
$(a\in A,\ b\in B, \tau\in T)$. 

Then the $B$-module $T_B:=B\otimes_AT$ admits a canonical structure 
of a $B$-algebroid Lie. Namely, the Lie bracket on $T_B$ is given by 
$$
[b_1\otimes \tau_1,b_2\otimes\tau_2]=b_1b_2\otimes [\tau_1,\tau_2] + 
\tau_1(b_2)b_1\otimes\tau_2 - 
\tau_2(b_1)b_2\otimes\tau_1
\eqno{(0.2.1.1)}
$$
and the action of $T_B$ on $B$ is defined by 
$$
(b_1\otimes\tau)(b_2)=b_1\tau(b_2)
\eqno{(0.2.1.2)}
$$

\bigskip 

{\it Vertex algebras} 

\bigskip

{\bf 0.3.} In his original paper, [B], Borcherds defined 
a vertex algebra over an arbitrary commutative ring. However, 
later most people preferred to work over the complex numbers. 
It is fairly obvious that all the general theorems of Kac's 
book [K] are true over an arbitrary field of characteristic $0$. 
A little less obvious, but true, is that with a minor modification 
of the definitions, they remain true over an arbitrary  
commutative ring. 
Below we recall the definitions and results to be used in the sequel, 
and explain these modifications. 

Throughout this work, we will deal only with  
{\it $\BZ_{\geq 0}$-graded} vertex algebras. 

{\bf 0.4. Definition.} A {\it $\BZ_{\geq 0}$-graded conformal algebra} 
is a $\BZ_{\geq 0}$-graded $k$-module $C=\oplus\ C_i$, together with a 
family of endomorphisms
$$
\dpar^{(j)}:\ C\lra C,\ \text{of\ degree}\ j,\ j\in\BZ_{\geq 0},
$$
such that 
$$
\dpar^{(i)}\cdot\dpar^{(j)}=\binom{i+j}{i}\dpar^{(i+j)};\ 
\dpar^{(0)}=Id, 
\eqno{(0.4.1)}
$$ 
and a family of bilinear operations
$$
_{(n)}:\ C\times C\lra C,\ \ (a,b)\mapsto a_{(n)}b,\ 
\text{of\ degree\ }-n-1,\ \ n\in\BZ_{\geq 0},
$$
such that 
$$
(\dpar^{(i)}a)_{(n)}b=(-1)^i\binom{n}{i}a_{(n-i)}b
\eqno{(0.4.2)}
$$
$$
a_{(n)}b=(-1)^{n+1}\sum_{j=0}^\infty\ 
(-1)^j\dpar^{(j)} (b_{(n+j)}a)
\eqno{(0.4.3)}
$$
$$
a_{(m)}b_{(n)}c=b_{(n)}a_{(m)}c+
\sum_{j=0}^m\ \binom{m}{j}(a_{(j)}b)_{(m+n-j)}c
\eqno{(0.4.4)}
$$
for all $a,b,c\in C;\ m,n,i\in\BZ_{\geq 0}$. 

Cf. [K], Definition 2.7b. In each conformal algebra we have the following 
identities: 
$$
(a_{(m)}b)_{(n)}c=\sum_{j=0}^m\ (-1)^j\binom{m}{j}\bigl\{a_{(m-j)}b_{(n+j)}c
-(-1)^mb_{(m+n-j)}a_{(j)}c\bigl\}
\eqno{(0.4.5)}
$$
and
$$
\dpar^{(j)}(a_{(n)}b)=\sum_{p=0}^j\ \dpar^{(p)}a_{(n)}\dpar^{(j-p)}b
\eqno{(0.4.6)}
$$
We leave their direct proof to the reader. The proof of (0.4.5) uses 
only the axiom (0.4.4) (and (0.4.1)), and the proof of (0.4.6) uses the 
axioms (0.4.2) and (0.4.3). 

We also have 
$$
a_{(n)}\dpar^{(j)}b=\sum_{p=0}^j\ \binom{n}{p}\dpar^{(j-p)}(a_{(n-p)}b)
\eqno{(0.4.7)}
$$
This is proven by induction on $j$.   

{\bf 0.5.} {\bf First Definition of a Vertex Algebra.} A {\it $\BZ_{\geq 
0}$-graded 
vertex algebra} is a 
$\BZ_{\geq 0}$-graded $k$-module $V=\oplus\ V_i$, equipped with 
a distinguished vector $\b1\in V_0$ ({\it vacuum vector}) 
and a family of bilinear operations 
$$
_{(n)}:\ V\times V\lra V,\ (a, b)\mapsto a_{(n)}b,\ 
\text{\ of degree\ }-n-1,\ 
n\in \BZ,
\eqno{(0.5.1)}
$$
such that 
$$
\b1_{(n)}a=\delta_{n,-1}a;\ a_{(-1)}\b1=a;\ a_{(n)}\b1=0\ 
\text{\ if\ }n\geq 0,
\eqno{(0.5.2)}
$$
and 
$$
\sum_{j=0}^\infty\ \binom{m}{j}(a_{(n+j)}b)_{(m+l-j)}c=
$$
$$
=
\sum_{j=0}^\infty\ (-1)^j\binom{n}{j} 
\bigl\{a_{(m+n-j)}b_{(l+j)}c-(-1)^nb_{(n+l-j)}a_{(m+j)}c\bigr\}
\eqno{(0.5.3)}
$$
for all $a,b,c\in V,\ m,n,l\in \BZ$.   

Cf. [B], Section 4, [K], Prop. 4.8. (b). The important particular case of 
(0.5.3) corresponds to $m=0$: 
$$
(a_{(n)}b)_{(l)}c=\sum_{j=0}^\infty\ (-1)^j\binom{n}{j}
\bigl\{ a_{(n-j)}b_{(l+j)}c-(-1)^nb_{(n+l-j)}a_{(j)}c\bigr\}
\eqno{(0.5.4)}
$$
cf. (0.4.5). 
 
Set 
$$
\dpar^{(j)}a:=a_{(-1-j)}\b1,\ \ j\in\BZ_{\geq 0}
\eqno{(0.5.5)}
$$
This way we get endomorphisms $\dpar^{(j)}$ of $V$ of degree $j$. It follows 
from (0.5.2) that 
$$
\dpar^{(j)}\b1=\delta_{j,0}\b1
\eqno{(0.5.6)}
$$
and 
$$
\dpar^{(0)}=Id
\eqno{(0.5.7)}
$$
and (0.5.4) applied to $b=c=\b1$ gives 
$$
\dpar^{(i)}\cdot\dpar^{(j)}=\binom{i+j}{i}\dpar^{(i+j)}
\eqno{(0.5.8)}
$$
We have the {\it commutativity formula}   
$$
a_{(n)}b=(-1)^{n+1}\sum_{j=0}^\infty\ (-1)^j \dpar^{(j)}(b_{(n+j)}a)
\eqno{(0.5.9)}
$$
for all $a,b\in V,\ n\in\BZ$. The proof will be given below, see the 
last paragraph of the next subsection. 

One deduces from (0.5.4) that 
$$
(\dpar^{(j)}a)_{(n)}b=(-1)^j\binom{n}{j}a_{(n-j)}b
\eqno{(0.5.10)}
$$
and 
$$
\dpar^{(j)}(a_{(n)}b)=\sum_{p=0}^j\ (\dpar^{(p)}a)_{(n)}\dpar^{(j-p)}b
\eqno{(0.5.11)}
$$
for all $n\in\BZ$. 
We have the {\it Operator Product Expansion (OPE)} formula, 
$$
[x_{(m)},y_{(n)}]=\sum_{j\geq 0}\ \binom{m}{j}(x_{(j)}y)_{(m+n-j)}\ \ 
(m,n\in\BZ)
\eqno{(0.5.12)}
$$
cf. [K], (4.6.7) and the end of the next Subsection.      

{\bf 0.6.} Following [K], let us give an equivalent definition of a 
vertex algebra. If $M$ is a $k$-module, $M[[z,z^{-1}]]$ will denote 
as usually the module of formal power series $\sum_{n=-\infty}^\infty 
a_nz^n,\ a_n\in M$. Let us define endomorphisms 
$\dpar^{(j)},\ j\in\BZ_{\geq 0}$, of $M[[z,z^{-1}]]$ by 
$$
\dpar^{(j)}(\sum_n a_nz^n)=\sum_n\ \binom{n}{j}a_nz^{n-j}
\eqno{(0.6.1)}
$$

{\bf Second Definition of a Vertex Algebra.} A $\BZ_{\geq 0}$-graded vertex 
algebra is a $\BZ_{\geq 0}$-graded $k$-module $V=\oplus_i\ V_i$ equipped 
with 

--- a distingushed vector $\b1\in V_0$ ({\it vacuum vector}); 

--- a family of endomorphisms $\dpar^{(j)}:\ V\lra V$ of degree $j$, 
$j\in \BZ_{\geq 0}$, such that $\dpar^{(0)}=Id$ and 
$\dpar^{(i)}\cdot\dpar^{(j)}=
\binom{i+j}{i}\dpar^{(i+j)}$; 

--- a linear mapping $V\lra End(V)[[z,z^{-1}]],\ a\mapsto a(z)=\sum_n\ 
a_{(n)}z^{-n-1}$ such that $deg\ a_{(n)}=deg\ a -n-1$. 

These data should satisfy the following axioms. 

{\it Translation invariance.} For all $j\in\BZ_{\geq 0}$, $[\dpar^{(j)}, a(z)]=
\dpar^{(j)}a(z)$.

{\it Vacuum.} $\dpar^{(j)}\b1=\delta_{j,0}\b1;\ \b1(z)=\b1;\ 
a_{(n)}\b1=0$ for $n\geq 0$ and $a_{(-1)}\b1=a$. 

{\it Locality.} $(z-w)^N[a(z),b(w)]=0$ for $N>>0$. 

This is a modification of the definition given in [K], 4.1, which 
works over an arbitrary base commutative ring. Proposition 4.1 of 
{\it op. cit.} remains true if we understand by $e^{z\dpar}$ the expression
$$
e^{z\dpar}:=\sum_{j=0}^\infty\ z^j\dpar^{(j)}
\eqno{(0.6.2)}
$$
(we cannot use Lemma 4.1 of {\it op.cit.} anymore!). The argument 
of {\it op. cit.} 4.2 shows that the commutativity formula (0.5.9) is true 
for the vertex algebras in the second definition. Finally, the proof of 
{\it op. cit.}, Proposition 4.8 (b) shows that the first and the second 
definitions are equivalent. The proof from [K], 4.6 works to give the 
proof of the OPE formula (0.5.12).  

{\bf 0.7. Theorem.} {\it Let $V=\oplus V_i$ be a $\BZ_{\geq 0}$-graded 
$k$-module 
equipped with a distinguished vector $\b1\in V_0$ and a family of 
endomorphisms $\dpar^{(j)}$ of degree $j$, $j\in \BZ_{\geq 0}$ such that 
$\dpar^{(0)}=Id,\ \dpar^{(i)}\cdot\dpar^{(j)}=\binom{i+j}{i}\dpar^{(i+j)}$ 
and $\dpar^{(j)}\b1=\delta_{j,0}\b1$. 

Assume that we are given a family of homogeneous vectors $\{a^\alpha\}
\subset V$ and a family of formal power series ("distributions")  
$$
\{ a^\alpha(z)=\sum_n\ a^\alpha_{(n)}z^{-n-1}\}\subset End(V)[[z,z^{-1}]],\ 
deg\ a^\alpha_{(n)}=deg\ a^\alpha -n-1,
$$
satisfying the following conditions} 

{\bf (t)} {\it $[\dpar^{(j)},a^\alpha(z)]=\dpar^{(j)}a^\alpha(z)$};

{\bf (v)} {\it $a^\alpha_{(n)}\b1=0$ for $n\geq 0$; $a^\alpha_{-1}
\b1=a^\alpha$};  

{\bf (l)} {\it The distributions $a^\alpha(z)$ are mutually local, i.e. 
for all $\alpha, \beta$, $(z-w)^N[a^\alpha(z),a^\beta(w)]=0$ for sufficiently 
large $N$.}

{\it Let $S$ denote the set of all vectors of the form  
$$
a^{\alpha_1\ldots\alpha_N}_{j_1\ldots j_N}:=
a^{\alpha_1}_{(-1-j_1)}\ldots a^{\alpha_N}_{(-1-j_N)}\b1,\ 
N, j_i\in\BZ_{\geq 0},
\eqno{(0.7.1)}
$$
Assume that we are given a map 
$$
V\lra\ End(V)[[z,z^{-1}]],\ a\mapsto a(z)=\sum\ a_{(n)}z^{-n-1}
\eqno{(*)}
$$
having the following property:}

{\bf (P)} {\it There exists a subset $S'\subset S$ which generates $V$ 
as a $k$-module such that for each 
$a^{\alpha_1\ldots\alpha_N}_{j_1\ldots j_N}\in S'$, we have 
$$
a^{\alpha_1\ldots\alpha_N}_{j_1\ldots j_N}(z)=
:\dpar^{j_1}\alpha^{\alpha_1}(z)\ldots\dpar^{j_N}a^{\alpha_N}(z):
\eqno{(0.7.2)}
$$

Then the mapping} (*) {\it defines the structure of a $\BZ_{\geq 0}$-graded 
vertex algebra on $V$, and} (0.7.2) {\it holds true for all 
$a^{\alpha_1\ldots\alpha_N}_{j_1\ldots j_N}\in S$.} 

This is a version of {\it "Existence Theorem"}, [K], Theorem 4.5 and 
Corollary 4.5, whose proof goes through. 

{\bf 0.8.} Let $\CV ert,\ \CC onf$ denote the categories of $\BZ_{\geq 
0}$-graded 
vertex and conformal algebras respectively. We have a functor 
$$
c:\ \CV ert\lra \CC onf
\eqno{(0.8.1)}
$$
which assigns to a vertex algebra $V$ the same space $V$, with the 
operations $_{(n)}$, $n<0,$ forgotten. The axioms of a conformal algebra 
are satisfied due to (0.5.10), (0.5.9) and (0.5.12).  
This functor admits a left adjoint, to be constructed in Section 8 below.   

In the sequel we will call $\BZ_{\geq 0}$-graded vertex (conformal) 
algebras simply vertex (conformal) algebras. 

{\bf 0.9.} A vertex algebra $V$ is called {\it commutative} 
if $a_{(n)}b=0$ for all $a, b\in V,\ n\geq 0$ (Kac uses the term 
"holomorphic"). 

Let $V$ be a commutative vertex algebra. Then, with respect to 
the operation $ab:=a_{(-1)}b$, $V$ becomes a commutative associative 
algebra with the unity $\b1$. The operations with negative indices 
are recovered from the formula 
$$
a_{(-1-j)}b=(\dpar^{(j)}a)b
\eqno{(0.9.1)}
$$
This way we get an equivalence of categories 
$$
\CV ert\CA b\iso \dpar-\CA lg
\eqno{(0.9.2)}
$$
Here $\CV ert\CA b$ denotes the category of commutative vertex algebras 
and $\dpar-\CA lg$ denotes the category whose objects are $\BZ_{\geq 0}$-graded 
vector spaces $V=\oplus V_i$ equipped with a structure of a commutative 
algebra such that $V_i\cdot V_j\subset V_{i+j}$ and $\b1\in V_0$, and 
with a family of endomorphisms $\dpar^{(j)}$ of degree $j$, $j\in \BZ_{\geq 0}$, 
such that 
$$
\dpar^{(i)}\cdot\dpar^{(j)}=\binom{i+j}{i}\dpar^{(i+j)},\ 
\dpar^{(0)}=Id
\eqno{(0.9.3)}
$$
and
$$
\dpar^{(j)}(ab)=\sum_{p=0}^j\ \dpar^{(p)}a\cdot\dpar^{(j-p)}b
\eqno{(0.9.4)}
$$
Cf. [B]. Objects of $\dpar-\CA lg$ will be called {\it $\dpar$-algebras}.  

\bigskip\bigskip


\newpage 

\centerline{\bf \S\ 1. Vertex Algebroids}

\bigskip\bigskip

{\bf 1.1.} Let us define an {\it extended 
Lie algebroid} to be a quintuple $\CT=(A,T,\Omega,\dpar,\langle\ ,\ \rangle)$ 
where $A$ is a commutative $k$-algebra,  
$T$ is a Lie $A$-algebroid, $\Omega$ is an $A$-module equipped 
with a structure of a module over the Lie algebra $T$, $\dpar:\ A\lra\Omega$ 
is an $A$-derivation and a morphism of $T$-modules, 
$\langle\ ,\ \rangle :\ T\times\Omega\lra A$ is an $A$-bilinear pairing. 

These data must satisfy the following properties $(a\in A,\ 
\tau, \nu\in T,\ \omega\in\Omega)$:
$$
\langle\tau,\dpar a\rangle=\tau(a)
\eqno{(1.1.1)}
$$ 
$$
\tau(a\omega)=\tau(a)\omega+a\tau(\omega)
\eqno{(1.1.2)}
$$
$$
(a\tau)(\omega)=a\tau(\omega)+\langle\tau,\omega\rangle\dpar a
\eqno{(1.1.3)}
$$
$$
\tau(\langle\nu,\omega\rangle)=\langle [\tau,\nu],\omega\rangle+ 
\langle\nu,\tau(\omega)\rangle
\eqno{(1.1.4)}
$$
Cf. 0.1. We will also say that $\CT=(A,T,\ldots)$ is {\it an extended Lie 
$A$-algebroid}.  

Let define a {\it morphism} between two extended Lie algebroids 
$\CT=(A,T,\ldots)$ and $\CT'=(A',T',\ldots)$ to be a triple 
$g=(g_A,g_T,g_\Omega)$ 
where $g_A:\ A\lra A'$ is a morphism of $k$-algebras, $g_T:\ T
\lra T'$ is a morphism of Lie algebras and $A$-modules, $g_\Omega:\ 
\Omega\lra \Omega'$ is a morphism of $A$-modules. We require that the 
following properties should hold:
$$
g_A(\tau(a))=g_T(\tau)(g_A(a))
\eqno{(1.1.5)}
$$ 
$$
g_\Omega(\dpar a)=\dpar g_A(a)
\eqno{(1.1.6)}
$$
$$
g_A(\langle\tau,\omega\rangle=\langle g_T(\tau),g_\Omega(\omega)\rangle
\eqno{(1.1.7)}
$$
Composition of morphisms is defined in the obvious way. 
This way we get a category $\CL ieAlg$ of extended Lie algebroids.

{\bf 1.2.} Let $T$ be a Lie $A$-algebroid. Set $\Omega := 
Hom_A(T,A)$. Define $\dpar:\ A\lra\Omega$ by (1.1.1); let $\langle\ ,\ \rangle$ 
be the evident pairing. Action of $T$ on $\Omega$ is defined 
by (1.1.4). This way we get an extended Lie $A$-algebroid $\CT_T$.

Let us call an extended Lie algebroid $\CT=
(A,T,\Omega,\langle,\rangle)$ {\it perfect} if the pairing 
$\langle,\rangle$ induces an isomorphism $\Omega\iso Hom_A(T,A)$. 
 
The correspondence $T\mapsto\CT_T$ provides an equivalence of the 
category of Lie algebroids with the full subcategory $\CL ieAlg^{perf}
\subset\CL ieAlg$ of perfect extended Lie algebroids.  

{\bf 1.3. De Rham - Chevalley complex.} Let $\CT=(A,T,\Omega,\ldots)$ be 
an extended Lie $A$-algebroid. Let us define $A$-modules $\Omega^i=
\Omega^i(\CT),\ 
i\in\BZ_{\geq 0}$, as follows. Set $\Omega^0=A,\ \Omega^1=\Omega$. For 
$i\geq 2$, $\Omega^i$ is the submodule of the module of $A$-polylinear 
homomorphisms $h$ from $T^{i-1}$ to $\Omega$ such that the function 
$\langle\tau_1,h(\tau_2,\ldots,\tau_i)\rangle$ is skew symmetric with 
respect to all permutations of $(\tau_1,\ldots,\tau_i)$. 

For example, if $\CT$ is as in the previous example, then  
$\Omega^i=Hom_A(\Lambda^i_AT,A)$. 

Let us define the maps $d_{DR}=d_{DR}^i:\ \Omega^i\lra\Omega^{i+1}$ as follows. 
For $i=0$ we set $d_{DR}a=-\dpar a$. For $i\geq 1$ we set  
$$
d_{DR}h(\tau_1,\ldots,\tau_i)=d_{Lie}h(\tau_1,\ldots,\tau_i)-
\dpar\langle\tau_1,h(\tau_2,\ldots,\tau_i)\rangle
\eqno{(1.3.1)}
$$
where
$$
d_{Lie}h(\tau_1,\ldots,\tau_i)=\sum_{p=1}^i\ (-1)^{p+1}
\tau_p(h(\tau_1,\ldots,\widehat{\tau}_p,\ldots,\tau_i))+
$$
$$
+\sum_{1\leq p<q\leq i}\ (-1)^{p+q} 
h([\tau_p,\tau_q],\tau_1,\ldots,\widehat{\tau}_p,\ldots,\widehat{\tau}_q,
\ldots,\tau_i)
\eqno{(1.3.2)}
$$
For example, 
$$
d_{DR}\omega(\tau)=\tau(\omega)-\dpar\langle\tau,\omega\rangle,\ 
\eqno{(1.3.3)}
$$
for $\omega\in\Omega^1=\Omega$; and 
$$
d_{DR}h(\tau_1,\tau_2)=-h([\tau_1,\tau_2])+\tau_1(h(\tau_2))- 
\tau_2(h(\tau_1))-\dpar\langle\tau_1,h(\tau_2)\rangle,\ 
\eqno{(1.3.4)}
$$
for $h\in\Omega^2$. 

Let us introduce the action of the Lie algebra $T$ on the modules 
$\Omega^i$ by 
$$
\tau(h)(\tau_1,\ldots,\tau_{i-1})=\tau(h(\tau_1,\ldots,\tau_{i-1}))-
\sum_{p=1}^{i-1}\ h(\tau_1,\ldots,[\tau,\tau_p],\ldots,\tau_i)
\eqno{(1.3.5)}
$$
Let us define the convolution operators $\langle\tau,\cdot\rangle:\ 
\Omega^i\lra\Omega^{i-1}$ by
$$
\langle\tau,h\rangle(\tau_1,\ldots,\tau_{i-2}\rangle=
h(\tau,\tau_1,\ldots,\tau_{i-2})
\eqno{(1.3.6)}
$$
The maps $\{d_{DR}^i\}$ may be characterized as a unique collection 
of maps such that $d_{DR}^0=-\dpar$ and the {\it Cartan formula} 
$$
\tau(h)=\langle\tau,d_{DR}h\rangle+d_{DR}\langle\tau,h\rangle
\eqno{(1.3.7)}
$$
holds true. 

The maps $d_{DR}$ commute with the action of $T$. One checks that 
$d_{DR}^2=0$, so we get a complex $(\Omega^\bullet(\CT),d_{DR})$ 
called the {\it de Rham-Chevalley complex} of $\CT$. 

{\bf 1.4.} In the definition below we assume that $1/2\in k$. 
 
A {\it vertex algebroid} is a septuple $\CA=(A,T,\Omega,\dpar,
\gamma,\langle\ ,\ \rangle,c)$ where $A$ is a commutative $k$-algebra, 
$T$ is a Lie $A$-algebroid, $\Omega$ is an $A$-module equipped 
with an action of the Lie algebra $T$, $\dpar:\ A\lra\Omega$ is a derivation 
commuting with the $T$-action, 
$$
\langle\ ,\ \rangle:\ (T\oplus\Omega)\times (T\oplus\Omega)\lra A
$$
is a symmetric $k$-bilinear pairing equal to zero on $\Omega\times\Omega$ 
and such 
that $\CT_{\CA}=(A,T,\Omega,\dpar,\langle\ ,\ \rangle|_{T\times\Omega})$ 
is an extended Lie $A$-algebroid; $c:\ T\times T\lra\Omega$ is a 
skew symmetric $k$-bilinear pairing and $\gamma:\ A\times T\lra\Omega$ 
is a $k$-bilinear map. 

The following axioms must hold $(a,b\in A;\ \tau,\tau_i\in T)$: 
$$
\gamma(a,b\tau)=\gamma(ab,\tau)-a\gamma(b,\tau)-\tau(a)\dpar b-\tau(b)\dpar a
\eqno{(A1)}
$$
$$
\langle a\tau_1,\tau_2\rangle=a\langle\tau_1,\tau_2\rangle+
\langle\gamma(a,\tau_1),\tau_2\rangle-\tau_1\tau_2(a)
\eqno{(A2)}
$$
$$
c(a\tau_1,\tau_2)=ac(\tau_1,\tau_2)+\gamma(a,[\tau_1,\tau_2])-
\gamma(\tau_2(a),\tau_1)+\tau_2(\gamma(a,\tau_1))-
$$
$$
-\frac{1}{2}\langle\tau_1,\tau_2\rangle\dpar a+
\frac{1}{2}\dpar\tau_1\tau_2(a)-
\frac{1}{2}\dpar\langle\tau_2,\gamma(a,\tau_1)\rangle
\eqno{(A3)}
$$
$$
\langle [\tau_1,\tau_2],\tau_3\rangle+\langle\tau_2,[\tau_1,\tau_3]\rangle=
\tau_1(\langle\tau_2,\tau_3\rangle)
-\frac{1}{2}\tau_2(\langle\tau_1,\tau_3\rangle)-\frac{1}{2}\tau_3(
\langle\tau_1,\tau_2\rangle)+
$$
$$
+\langle\tau_2,c(\tau_1,\tau_3)\rangle+
\langle\tau_3,c(\tau_1,\tau_2)\rangle
\eqno{(A4)}
$$
$$
d_{Lie}c(\tau_1,\tau_2,\tau_3)=-\frac{1}{2}\dpar\bigl\{
\langle [\tau_1,\tau_2],\tau_3\rangle+\langle [\tau_1,\tau_3],\tau_2\rangle
-\langle [\tau_2,\tau_3],\tau_1\rangle-
$$
$$
-\tau_1(\langle\tau_2,\tau_3\rangle)+
\tau_2(\langle\tau_1,\tau_3\rangle)-2\langle\tau_3,c(\tau_1,\tau_2)\rangle
\bigr\}
\eqno{(A5)}
$$
where $d_{Lie}$ is defined by (1.3.2). 

To stress the dependence on $A$, we shall sometimes say that $\CA$ is a 
{\it vertex $A$-algebroid}.  

{\bf 1.5.} Multiplying the identity (A4) by $l$ and subtracting the $j$-th 
multiple of  
(A4) corresponding to the triple $(\tau_2,\tau_1,\tau_3)$ we get 
an equivalent form of (A4): 

for every $l,j\in k$, 
$$
-j\langle\tau_1,[\tau_2,\tau_3]\rangle+l\langle\tau_2,[\tau_1,\tau_3]\rangle+
(l+j)\langle\tau_3,[\tau_1,\tau_2]\rangle
$$
$$
-(l+\frac{1}{2}j)\tau_1(\langle\tau_2,\tau_3\rangle)+(\frac{1}{2}l+j)\tau_2
(\langle\tau_1,\tau_3\rangle)+\frac{1}{2}(l-j)\tau_3(\langle\tau_1,\tau_2
\rangle)
\eqno{(A4)^{bis}_{l,j}}
$$
$$
+j\langle\tau_1,c(\tau_2,\tau_3)\rangle-l\langle\tau_2,c(\tau_1,\tau_3)
\rangle-(l+j)\langle\tau_3,c(\tau_1,\tau_2)\rangle\ =\ 0 
$$

{\bf 1.6.} The left hand side of (A5) is skew symmetric with respect to all 
permutations of $(\tau_1,\tau_2,\tau_3)$. The right hand side 
is manifestly symmetric only with respect to the transposition of 
$(\tau_1,\tau_2)$. However, if we replace $\tau_1(\langle\tau_2,\tau_3
\rangle)$ in the right hand side by its expression from (A4) we get the 
following equivalent form of (A5): 
$$
d_{Lie}c(\tau_1,\tau_2,\tau_3)=-\frac{1}{2}\dpar\bigl\{
-\langle\tau_1,[\tau_2,\tau_3]\rangle+\frac{1}{2}\tau_2(\langle\tau_1,
\tau_3\rangle)-\frac{1}{2}\tau_3(\langle\tau_1,\tau_2\rangle)+
$$
$$
+\langle\tau_2,c(\tau_1,\tau_3)\rangle-\langle\tau_3,c(\tau_1,\tau_2)\rangle
\bigr\}
\eqno{(A5)^{bis}}
$$
The rhs of $(A5)^{bis}$ is skew symmetric with respect to the transposition 
of $(\tau_2,\tau_3)$. Consequently, the rhs of (A5) is completely 
skew symmetric since the symmetric group $S_3$ is generated by two 
transpositions $(12)$ and $(23)$. 

Let us replace in $(A5)^{bis}$ the triple $(\tau_1,\tau_2,\tau_3)$ 
by the triples $(\tau_2,\tau_3,\tau_1)$ and $(\tau_3,\tau_1,\tau_2)$ and 
sum up the three identities. We will get another equivalent for of 
(A5): 
$$
3d_{Lie}c(\tau_1,\tau_2,\tau_3)=\dpar\bigl\{\langle\tau_1,\frac{1}{2}
[\tau_2,\tau_3]+c(\tau_2,\tau_3)\rangle+
$$
$$
+\langle\tau_2,\frac{1}{2}[\tau_3,\tau_1]+c(\tau_3,\tau_1)\rangle
+\langle\tau_3,\frac{1}{2}[\tau_1,\tau_2]+c(\tau_1,\tau_2)\rangle
\bigr\}
\eqno{(A5)^{ter}}
$$

{\bf 1.7. Example.} Let $\fg$ be a Lie algebra over $k$ equipped with a 
symmetric invariant form $\langle\ ,\ \rangle:\ \fg\times\fg\lra k$. 
Then $\CA_{\fg,\langle,\rangle}=(k,\fg,0,0,0,\langle\ ,\ \rangle,0)$ 
is a vertex $k$-algebroid. 

{\bf 1.8.} Let $A$ be a commutative $k$-algebra and $\CT=(A,T,\Omega,\dpar,
\langle\ ,\ \rangle)$ be an extended Lie $A$-algebroid. Let $T_0\subset T$ 
be a $k$-submodule which generates $T$ as an $A$-module, and assume that we are 
given $k$-bilinear mappings 
$$
\gamma:\ A\times T_0\lra\Omega;\ 
\langle\ ,\ \rangle:\ T_0\times T_0\lra A,\ \text{symmetric};\ 
$$
$$
c:\ T_0\times T_0\lra\Omega,\ \text{skew symmetric}
\eqno{(1.8.1)}
$$
It is clear that there exists not more than one extension of the mappings 
to mappings
$$
\gamma:\ A\times T\lra\Omega;\ 
\langle\ ,\ \rangle:\ T\times T\lra A,\  
c:\ T\times T\lra\Omega,\ 
\eqno{(1.8.2)}
$$
satisfying (A1), (A2), (A3). 
This extension, if it exists, must be given by the formulas $(a,b\in A,\ \tau, 
\tau_i\in T_0)$:  
$$
\gamma(a,b\tau)=\gamma(ab,\tau)-a\gamma(b,\tau)-\tau(a)\dpar b-
\tau(b)\dpar a
\eqno{(1.8.3)_{\gamma}}
$$
$$
\langle a\tau_1,b\tau_2\rangle=ab\langle\tau_1,\tau_2\rangle+a\langle
\tau_1,\gamma(b,\tau_2)\rangle+b\langle\tau_2,\gamma(a,\tau_1)\rangle-
$$
$$
-a\tau_2\tau_1(b)-b\tau_1\tau_2(a)-\tau_1(b)\tau_2(a)
\eqno{(1.8.3)_{\langle,\rangle}}
$$
$$
c(a\tau_1,b\tau_2)=ab\cdot c(\tau_1,\tau_2)+\gamma(ab,[\tau_1,\tau_2])+
$$
$$
+\gamma(a\tau_1(b),\tau_2)-\gamma(b\tau_2(a),\tau_1)-a\tau_1(\gamma(b,\tau_2))+
b\tau_2(\gamma(a,\tau_1))+
$$
$$
+\frac{1}{2}\langle\tau_1,\tau_2\rangle\bigl\{a\dpar b-b\dpar a\bigr\}
-[\tau_1,\tau_2](a)\dpar b-[\tau_1,\tau_2](b)\dpar a+
\frac{1}{2}\bigl\{\tau_1(b)\dpar\tau_2(a)-\tau_2(a)\dpar\tau_1(b)\bigr\}+
$$
$$
+\frac{1}{2}\bigl\{-\langle\tau_1,\gamma(b,\tau_2)\rangle\dpar a+
\langle\tau_2,\gamma(a,\tau_1)\rangle\dpar b+
a\dpar\langle\tau_1,\gamma(b,\tau_2)\rangle-b\dpar\langle\tau_2,\gamma(a,
\tau_1)\rangle\bigr\}+
$$
$$
+\frac{1}{2}\dpar\bigl\{b\tau_1\tau_2(a)-a\tau_2\tau_1(b)\bigr\}
\eqno{(1.8.3)_c}
$$
When deducing the last formula, one should take into account (1.1.3).

{\bf 1.9. Theorem.} {(\it Extension of Identities)} 
{\it Assume that (A4) and (A5) hold true 
for all $\tau_i\in T_0$ and that the formulas} (1.8.3) {\it provide well 
defined mappings} (1.8.2). 

{\it Then the axioms (A1) --- (A5) hold true for all $\tau, \tau_i\in T$, 
i.e. $\CA=(A,T,\ldots)$ is a vertex $A$-algebroid.} 

{\it Proof.} Let us check the axiom (A1). It is enough to show that if 
(A1) holds true for some $\tau$ and all $a, b$ then it holds true 
for $c\tau\ (c\in A)$ and all $a,b$. Thus, we have to check the identity 
$$
\gamma(ab,c\tau)=\gamma(a,bc\tau)+a\gamma(b,c\tau)+c\tau(a)\dpar b+
c\tau(b)\dpar a
\eqno{(?)}
$$
The left hand side of it equal to 
$$
\gamma(abc,\tau)-ab\gamma(c,\tau)-\tau(ab)\dpar c-\tau(c)\dpar(ab),
$$
cf. $(1.8.3)_\gamma$. On the other hand, the first two terms in the right 
hand side are equal to
$$
\gamma(a,bc\tau)=\gamma(abc,\tau)-a\gamma(bc,\tau)-\tau(a)\dpar(bc)-
\tau(bc)\dpar a
$$
and
$$
a\gamma(b,c\tau)=a\gamma(bc,\tau)-ab\gamma(c,\tau)-a\tau(b)\dpar c-
a\tau(c)\dpar b,
$$
again by $(1.8.3)_\gamma$. Comparing, we get the desired identity (?). 

The other axioms are checked in a similar way. This is a tiresome but 
straightforward calculation. When checking (A5) one should use the 
identity $(A4)_{1,1}^{bis}$. $\btu$ 

{\bf 1.10. Pushout.} Let $\CA=(A,T,\Omega,\ldots)$ be a vertex 
$A$-algebroid. Let $B$ be a commutative $A$-algebra, $i:\ A\lra B$ 
the structure morphism. Set $\Omega_B:=B\otimes_A\Omega,\ 
T_B:=B\otimes T$. The $A$-derivation $\dpar:\ A\lra \Omega$ induces a 
$B$-derivation $\dpar_B:\ B\lra \Omega_B$. The $A$-bilinear pairing 
$\langle , \rangle:\ T\times \Omega\lra A$ uniquely extends to a 
$B$-bilinear pairing $\langle,\rangle_B:\ T_B\times\Omega_B\lra B$.  

Assume that the Lie algebra $T$ acts on $B$ by derivations in such a way that 
$\tau(i(a))=i(\tau(a))$ and $(a\tau)(b)=a\tau(b)\ (a\in A,\ b\in B,\ 
\tau\in T)$. Then $T_B$ acquires a canonical structure of a Lie $B$-algebroid, 
cf. 0.2.1, and $(T_B,\Omega_B,\dpar_B,\langle,\rangle_B)$ becomes an 
extended Lie $B$-algebroid.  

{\bf 1.10.1. Theorem.} {\it 
Assume that we are given a $k$-bilinear mapping $\gamma:\ 
B\times T\lra\Omega_B$ such that $\gamma(i(a),\tau)=1\otimes\gamma(a,\tau)$ and 
that (A1) holds true for all $\tau\in T,\ a\in B,\ b\in A$. 

Then there exists a unique extension of $\gamma$ to a $k$-bilinear 
mapping $\gamma_B:\ B\times T_B\lra\Omega_B$ satisfying (A1) for all 
$a,b\in B,\ \tau\in T_B$; 

there exists a unique extension of the pairing $\langle, \rangle:\ 
T\times T\lra A$ to a pairing $\langle,\rangle_B:\ T_B\times T_B\lra B$ 
satisfying (A2) for all $a\in B,\ \tau_i\in T_B$; 

there exists a unique extension of the pairing $c:\ T\times T\lra\Omega$ to 
a pairing $c_B:\ T_B\times T_B\lra\Omega_B$ satisfying (A3) for all 
$a\in B,\ \tau_i\in T_B$. 

The septuple $\CA_B=(B,T_B,\Omega_B,\dpar_B,\gamma_B,\langle,\rangle_B,c_B)$ 
is a vertex $B$-algebroid.}

{\it Proof.} Apply 1.9 to $T_0:=Im(T\lra T_B)$. $\btu$ 

{\bf 1.11. Example.} Let $A$ be a commutative $k$-algebra, $\Omega=\Omega^1_k$ 
--- 
the $A$-module of K\"ahlerian $1$-differentials, $\dpar:\ A\lra\Omega$ 
the canonical $A$-derivation. Let $\fg$ be a 
Lie $k$-algebra acting on $A$ by derivations and equipped with an 
invariant bilinear form $\langle,\rangle:\ \fg\times\fg\lra k$. 

Due to the morphism $\fg\lra Der_k(A)$,  
$\Omega$ is equipped with a canonical structure of a $\fg$-module  
such that $\dpar$ is a morphism of $\fg$-modules and there is a canonical 
pairing $\fg\times\Omega\lra A$, cf. 0.1. 

Set $T=A\otimes\fg$. Then $T$ is canonically a Lie $A$-algebroid, there 
exists a unique extension of the zero map $A\times\fg\lra A$ (resp., of  
$\langle,\rangle$ and  of the zero map $\fg\times\fg\lra\Omega$) to the map 
$\gamma:\ A\times T\lra\Omega$ (resp., to the pairing $\langle,\rangle:\ 
T\times T\lra A$ and to the map $c:\ T\times T\lra\Omega$) which satisfies 
(A1) (resp., (A2) and (A3)). 

This way we get a vertex $A$-algebroid $\CA=(A,T,\Omega,\dpar,\gamma,
\langle,\rangle,c)$. 

For the proof, apply Theorem 1.9 to $T_0=\fg\subset T$. 

{\bf 1.12. Example.} Let $A$ be a commutative $k$ algebra, $\fg\subset T=
Der_k(A)$ a Lie $k$-subalgebra such that $T=A\otimes\fg$, $\Omega=\Omega^1_k$.  
 
For example, $A$ is the coordinate ring of an algebraic group, 
$\fg$ is the Lie algebra of left invariant vector fields. 

Assume that $\fg$ is equipped with an invariant bilinear form $\langle,\rangle$.
Then we are in the situation of 1.11 and get a vertex $A$-algebroid 
$\CA=\CA_{A,\fg,\langle,\rangle}$.  

{\bf 1.13. Example.} Let $A$ be a commutative $k$-algebra equipped with 
an \'etale morphism $\text{{\bf{f}}}=(f_1,\ldots,f_n):\ A_0=k[x_1,\ldots,x_n]
\lra A$. 
Such an {\bf{f}} is called an {\it \'etale coordinate system} on $A$; 
it exists Zariski locally for every $A$ smooth over $k$.  

The commuting vector fields $\dpar/\dpar x_i\in Der_k(A_0)$ admit 
unique liftings to vector fields 
$\tau_i\in T=Der_k(A)$. The $k$-submodule $\fg\subset T$ spanned by $\tau_i$ 
is an abelian Lie subalgebra; let us equip it with the trivial bilinear form. 
Then we are in the situation of 1.11 and we get a vertex $A$-algebroid 
$\CA=\CA_{A;\text{{\bf{f}}}}$. 

\bigskip\bigskip



\centerline{\bf \S 2. From Vertex Algebras to Vertex Algebroids}

\bigskip\bigskip

{\bf 2.1.} Let $V$ be a vertex algebra over $k$, cf. 0.5. In this section 
we assume that $1/2\in k$; recall that we always deal with 
$\BZ_{\geq 0}$-graded vertex algebras. We denote $\dpar:=\dpar^{(1)}$. 

The identity (0.5.4) at $n=0$ gives
$$
a_{(0)}b_{(n)}c=(a_{(0)}b)_{(n)}c+b_{(n)}a_{(0)}c
\eqno{(2.1.1)}
$$
for all $n$, i.e. the operation $_{(0)}$ is a derivation with respect 
to all operations $_{(n)}$.  

We shall denote the operation $a_{(-1)}b$ by $ab$ or $a\cdot b$.  

Set $A=V_0$. In the sequel 
$a,b,c$ will denote elements of $A$.  Since 
$a_{(n)}b=0$ for all $n\geq 0$, (0.5.4) and (0.5.9) imply that $A$ is a 
commutative associative $k$-algebra with the unity $\b1$ with respect to 
the operation $ab$. 

Elements of $V_1$ will be denoted $x,y,z$. It follows from (0.5.4), 
taking into account that $\binom{-1}{j}=(-1)^j$, that
$$
(ab)x=abx+\dpar a\ b_{(0)}x+\dpar b\ a_{(0)}x
\eqno{(2.1.2)}
$$
The identity (0.5.9) implies
$$
a_{(0)}x=-x_{(0)}a
\eqno{(2.1.3)}
$$
and
$$
xa=ax+\dpar(x_{(0)}a)
\eqno{(2.1.4)}
$$
Let $\Omega\subset V_1$ be the $k$-submodule generated by all elements 
$a\dpar b$; $\omega$ will denote an element of $\Omega$.  

From (2.1.2) we have 
$$
(ab)\dpar c=ab\dpar c+\dpar a\ b_{(0)}\dpar c+\dpar b\ a_{(0)}\dpar c
$$
Note that 
$$
\dpar u_{(0)}v=0
\eqno{(2.1.5)}
$$
for all $u,v\in V$, due to (0.5.10). Hence $b_{(0)}\dpar c=-\dpar c_{(0)}b=0$; 
therefore $(ab)\dpar c=ab\dpar c$. On the other hand, $\b1 \dpar c=
\dpar c \b1-\dpar c_{(0)}\b1=\dpar c \b1=\dpar c$. It follows that 
the operation $a\omega$ gives a structure of an $A$-module on $\Omega$. 
(It is not true for $V_1$ since $ax$ is not associative.) 

Note that by (0.5.4) $(a\dpar b)_{(0)}c=a\dpar b_{(0)}c=0$, in other 
words 
$$
\omega_{(0)}a=0
\eqno{(2.1.6)}
$$
It follows from (0.5.9) that 
$$
a\omega=\omega a
\eqno{(2.1.7)}
$$
We have the map $\dpar:\ A\lra\Omega$. By (0.5.11) and (2.1.7) 
$\dpar$ is an $A$-derivation. 

Let $T$ denote the quotient module $V_1/\Omega$; let $\pi:\ V_1\lra T$ 
be the canonical projection. Elements of $T$ will be denoted $\tau$. 
The operation $ax$ induces an operation $A\times T\lra T,\ (a,\tau)\mapsto 
a\tau$. By (2.1.2) $(ab)\tau=ab\tau$. On the other hand, 
$$
\b1 x=x \b1-\dpar(x_{(0)}\b1)=x\b1=x
\eqno{(2.1.8)}
$$
i.e. the operation $a\tau$ provides $T$ with a structure of an $A$-module. 

By (0.5.4)
$$  
(a\dpar b)_{(0)}x=a\dpar b_{(0)}x-a_{(-2)}\dpar b_{(1)}x+
\dpar b\ a_{(0)}x
$$
Note that by (0.5.10) and (2.1.7) 
$$
a_{(-2)}b=\dpar a\ b=b\dpar a
\eqno{(2.1.9)}
$$
Therefore 
$$
(a\dpar b)_{(0)}x=(b_{(0)}x)\dpar a+(a_{(0}x)\dpar b
$$
whence
$$
\Omega_{(0)}V_1\subset\Omega
\eqno{(2.1.10)}
$$ 
On the other hand, 
$$
x_{(0)}a\dpar b=(x_{(0)}a)\dpar b+ax_{(0)}\dpar b
$$
by (2.1.1), and 
$$
x_{(0)}\dpar b=-\dpar b_{(0)}x+\dpar(\dpar b_{(1)}x)=-\dpar(b_{(0)}x)=
\dpar(x_{(0)}b),
$$
i.e. 
$$
x_{(0)}(a\dpar b)=(x_{(0)}a)\dpar b+a\dpar(x_{(0)}b)
\eqno{(2.1.11)}
$$
It follows that 
$$
V_{1(0)}\Omega\subset\Omega
\eqno{(2.1.12)}
$$
Therefore the operation $_{(0)}$ induces an operation $T\times T\lra T$, 
to be denoted $[,]$. By (0.5.9)
$$
x_{(0)}y=-y_{(0)}x+\dpar(y_{(1)}x)
\eqno{(2.1.13)}
$$
hence $[,]$ is skew symmetric. (2.1.1) implies that $[,]$ satisfies the 
Jacobi identity, hence it provides a structure of a Lie algebra (over $k$) 
on $T$. 

By (0.5.4) and (2.1.5) $(a\dpar b)_{(0)}c=a\dpar b_{(0)}c=0$, i.e. 
$$
\Omega_{(0)}A=0
\eqno{(2.1.14)}
$$
It follows that the operation $x_{(0)}a$ induces a pairing $T\times A\lra A$, 
to be denoted $\tau(a)$. By (2.1.1), $\tau_1\tau_2(a)=[\tau_1,\tau_2](a)+
\tau_2\tau_1(a)$ and $\tau(ab)=\tau(a)b+a\tau(b)$. By (0.5.4), $(a\tau)(b)=
a\tau(b)$. Thus, $T$ is equipped with a structure of a Lie $A$-algebroid. 

By (2.1.11) and (2.1.14) 
$$
\Omega_{(0)}\Omega=0
\eqno{(2.1.15)}
$$
It follows that the operation $x_{(0)}\omega$ induces an operation 
$T\times\Omega\lra\Omega$, to be denoted $\tau(\omega)$, which makes $\Omega$ 
a module over the Lie algebra $T$, by (2.1.1). Again by (2.1.1), 
$\tau(a\omega)=\tau(a)\omega+a\tau(\omega)$. 

It follows from (2.1.11) that 
$$
\tau(a\dpar b)=\tau(a)\dpar b+a\dpar\tau(b)
\eqno{(2.1.16)}
$$
In particular $\tau(\dpar b)=\dpar\tau(b)$, i.e. $\dpar$ is a morphism 
of $T$-modules. 

By (0.5.9) 
$$
x_{(1)}y=y_{(1)}x
\eqno{(2.1.17)}
$$
We have 
$$
(a\dpar b)_{(1)}x=a\dpar b_{(1)}x=-a\dpar b_{(0)}x=ax_{(0)}b
\eqno{(2.1.18)}
$$
by (0.5.4), (2.1.5), (0.5.10) and (2.1.3). It follows that 
$$
\Omega_{(1)}\Omega=0
\eqno{(2.1.19)}
$$
by (2.1.14). Therefore the operation $x_{(1)}y$ induces an operation 
$T\times\Omega\lra A$, to be denoted $\langle\tau,\omega\rangle$. 

By (2.1.17) and (2.1.18) 
$$ 
\langle\tau,a\dpar b\rangle=a\tau(b) 
\eqno{(2.1.20)}
$$ 
From (0.5.4), $\langle a\tau,\omega\rangle=a\langle\tau,\omega\rangle$.  From
(2.1.1), 
$$
\tau(\langle\tau',\omega\rangle)=\langle [\tau,\tau'],
\omega\rangle+\langle\tau',\tau(\omega)\rangle
$$  
Finally, again from (0.5.4),
$(ax)_{(0)}\omega =ax_{(0)}\omega +a_{(-2)}x_{(1)}\omega$ which implies
$(a\tau)(\omega)=a\tau(\omega)+\langle\tau,\omega\rangle\dpar a$, by (2.1.9).

Therefore we have canonically associated with our vertex algebra an 
extended Lie $A$-algebroid $\CT(V)=(T,\Omega,\dpar,\langle,\rangle)$. 

{\bf 2.2.} Let us assume that the projection $\pi:\ V_1\lra T$ admits 
a splitting, i.e. there exists a morphism of $k$-modules 
$s:\ T\lra V_1$ such that $\pi\circ s=Id_T$.  
Let us fix such a splitting. A vertex algebra with a chosen $s$ will 
be called {\it split}. 

Define a skew symmetric operation $[,]:\ V_1\times V_1\lra V_1$ by 
$$
[x,y]:=\frac{1}{2}(x_{(0)}y-y_{(0)}x)
\eqno{(2.2.1)}
$$
We have $\pi([x,y])=[\pi(x),\pi(y)]$. 

Let us use the notation 
$\langle x,y\rangle$ for $x_{(1)}y$. It follows form (0.5.9) that 
$$
[x,y]=x_{(0)}y-\frac{1}{2}\dpar\langle x,y\rangle
\eqno{(2.2.2)}
$$
Set
$$
\gamma(a,\tau)=s(a\tau)-as(\tau)
\eqno{(2.2.3)_{\gamma}}
$$
$$
\langle\tau_1,\tau_2\rangle=\langle s(\tau_1),s(\tau_2)\rangle
\eqno{(2.2.3)_{\langle,\rangle}}
$$
$$
c(\tau_1,\tau_2)=s([\tau_1,\tau_2])-[s(\tau_1),s(\tau_2)]
\eqno{(2.2.3)_c}
$$

{\bf 2.3. Theorem.} {\it The septuple $\CA=(A,T,\Omega,\dpar,\gamma,
\langle,\rangle,c)$ is a vertex algebroid.} 

{\it Proof.} We have to check the axioms (A1) --- (A5) from 1.4. 

{\it Check of} (A1). We have $\gamma(ab,\tau)=s(ab\tau)-(ab)s(\tau)
=s(ab\tau)-as(b\tau)+as(b\tau)-(ab)s(\tau)$. By (2.1.2) and (2.1.3), 
$$
(ab)s(\tau)=ab s(\tau)+\dpar a\ b_{(0)}s(\tau)+\dpar b\ a_{(0)}s(\tau)=
ab s(\tau)-\tau(b)\dpar a-\tau(a)\dpar b
$$
which implies (A1). 

{\it Check of} (A2). We have $\langle a\tau_1,\tau_2\rangle=\langle s(a\tau_1),
s(\tau_2)\rangle=\langle as(\tau_1)+\gamma(a,\tau_1),s(\tau_2)\rangle$. 
By (0.5.4), 
$$
\langle ax,y\rangle=a\langle x,y\rangle-x_{(0)}y_{(0)}a
\eqno{(2.3.1)}
$$
which implies (A2). 

{\it Check of} (A3). We have $c(a\tau_1,\tau_2)=s([a\tau_1,\tau_2])-
[s(a\tau_1),s(\tau_2)]$; 
$$
s([a\tau_1,\tau_2])=s(a[\tau_1,\tau_2]-\tau_2(a)\tau_1)=as([\tau_1,\tau_2])+
\gamma(a,[\tau_1,\tau_2])-\tau_2(a)s(\tau_1)-\gamma(\tau_2(a),\tau_1); 
$$
$[s(a\tau_1),s(\tau_2)]=[as(\tau_1)+\gamma(a,\tau_1),s(\tau_2)]$. 

We have $[ax,y]=(ax)_{(0)}y-\frac{1}{2}\dpar\langle ax,y\rangle$ 
(see (2.2.2)). By (0.5.4), 
$$
(ax)_{(0)}y=ax_{(0)}y+a_{(-2)}x_{(1)}y+xa_{(0)}y=
$$
$$
=ax_{(0)}y+\langle x,y
\rangle\dpar a-\pi(y)(a)x-\dpar\pi(x)\pi(y)(a)
\eqno{(2.3.2)}
$$
By (2.3.1) 
$$
-\frac{1}{2}\dpar\langle ax,y\rangle=-\frac{1}{2}\dpar\bigl\{a\langle x,y
\rangle-x_{(0)}y_{(0)}a\bigr\}=
-\frac{1}{2}\langle x,y\rangle\dpar a-\frac{1}{2}a\dpar\langle x,y\rangle+
\frac{1}{2}\dpar\pi(x)\pi(y)(a)
$$
Therefore
$$
[ax,y]=a[x,y]+\frac{1}{2}\langle x,y\rangle\dpar a-\frac{1}{2}\dpar\pi(x)
\pi(y)(a)-\pi(y)(a)x
$$
On the other hand, 
$$
[\omega,x]=-[x,\omega]=-\pi(x)\omega+\frac{1}{2}\dpar\langle x,\omega\rangle
$$
It follows that 
$$
[s(a\tau_1),s(\tau_2)]=a[s(\tau_1),s(\tau_2)]+\frac{1}{2}\langle\tau_1,
\tau_2\rangle\dpar a-\frac{1}{2}\dpar\tau_1\tau_2(a)-\tau_2(a)s(\tau_1)-
\tau_2(\gamma(a,\tau_1))+
$$
$+\frac{1}{2}\dpar\langle\gamma(a,\tau_1),\tau_2\rangle$. This implies 
(A3). 

{\it Check of} (A4). We shall make use of the formula 
$$
s(\tau_1)_{(0)}\langle s(\tau_2),s(\tau_3)\rangle=\langle s(\tau_1)
_{(0)}s(\tau_2),s(\tau_3)\rangle+\langle s(\tau_2),s(\tau_1)_{(0)}s(\tau_3)
\rangle
\eqno{(2.3.3)}
$$
which is a particular case of (2.1.1). We have $s(\tau_1)_{(0)}
\langle s(\tau_2),s(\tau_3)\rangle=\tau_1(\langle\tau_2,\tau_3\rangle)$. 
On the other hand, 
$$
s(\tau_1)_{(0)}s(\tau_2)=[s(\tau_1),s(\tau_2)]+\frac{1}{2}\dpar\langle
\tau_1,\tau_2\rangle=s([\tau_1,\tau_2])-c(\tau_1,\tau_2)+\frac{1}{2}
\dpar\langle\tau_1,\tau_2\rangle
\eqno{(2.3.4)}
$$
Therefore
$$
\langle s(\tau_1)_{(0)}s(\tau_2),s(\tau_3)\rangle=\langle[\tau_1,\tau_2],
\tau_3\rangle-\langle c(\tau_1,\tau_2),\tau_3\rangle+\frac{1}{2}
\tau_3(\langle\tau_1,\tau_2\rangle)
$$
Interchanging $\tau_2$ with $\tau_3$ and plugging these expressions into 
(2.3.3) we get (A4). 

{\it Check of} (A5). We shall use the formulas 
$$
s(\tau_1)_{(0)}s(\tau_2)_{(0)}s(\tau_3)=(s(\tau_1)_{(0)}s(\tau_2))_{(0)}
s(\tau_3)+s(\tau_2)_{(0)}s(\tau_1)_{(0)}s(\tau_3) 
$$
and (2.3.4). We get
$$
s(\tau_1)_{(0)}s(\tau_2)_{(0)}s(\tau_3)=s(\tau_1)_{(0)}\bigl\{s([\tau_2,
\tau_3])-c(\tau_2,\tau_3)+\frac{1}{2}\dpar\langle\tau_2,\tau_3\rangle
\bigr\}=
$$
$$
=s([\tau_1,[\tau_2,\tau_3])-c(\tau_1,[\tau_2,\tau_3])+\frac{1}{2}\dpar
\langle\tau_1,[\tau_2,\tau_3]\rangle-\tau_1(c(\tau_2,\tau_3))+\frac{1}{2}
\dpar\tau_1(\langle\tau_2,\tau_3\rangle)
$$
Interchanging $\tau_1$ with $\tau_2$ we get
$$
s(\tau_2)_{(0)}s(\tau_1)_{(0)}s(\tau_3)=s([\tau_2,[\tau_1,\tau_3])-
$$
$$
-c(\tau_2,[\tau_1,\tau_3])+\frac{1}{2}\dpar
\langle\tau_2,[\tau_1,\tau_3]\rangle-\tau_2(c(\tau_1,\tau_3))+\frac{1}{2}
\dpar\tau_2(\langle\tau_1,\tau_3\rangle)
$$
Similarly, 
$$
(s(\tau_1)_{(0)}s(\tau_2))_{(0)}s(\tau_3)=\bigl\{s([\tau_1,\tau_2])-
c(\tau_1,\tau_2)+\frac{1}{2}\dpar\langle\tau_1,\tau_2\rangle\bigr\}_{(0)}
s(\tau_3)=
$$
$$
=s([\tau_1,\tau_2],\tau_3])-c([\tau_1,\tau_2],\tau_3)+\frac{1}{2}
\dpar\langle [\tau_1,\tau_2],\tau_3\rangle+\tau_3(c(\tau_1,\tau_2))-
\dpar\langle\tau_3,c(\tau_1,\tau_2)\rangle
$$
where we have used the formula 
$$
\omega_{(0)}x=-\pi(x)(\omega)+\dpar\langle\pi(x),\omega\rangle
\eqno{(2.3.5)}
$$ 
which is a consequence of (0.5.9). 
The axiom (A5) follows. 
This completes the proof of the theorem. $\btu$ 

{\bf 2.4.} Thus, to a pair $(V,s)$, where $V$ is a vertex algebra and  
$s:\ T\lra V_1$ is a splitting on $V$,  
we have assigned 
a vertex algebroid $\CA(V,s)=(A,T,\Omega,\dpar,\gamma,\langle,\rangle,c)$. 
Note that by definition $\CA(V,s)$ has the property 

(Sur) The module $\Omega$ is generated as an $A$-module by the subspace  
$\dpar A$. Equivalently, the  
morphism $\Omega^1_{A/k}\lra\Omega$ induced by $\dpar$ is epimorphic.

\bigskip\bigskip



\centerline{\bf \S 3. Category of Vertex Algebroids}

\bigskip\bigskip

{\bf 3.1.} Let us define a {\it $1$-truncated vertex algebra} to be a 
septuple $v=(V_0,V_1,\b1,\dpar,_{(-1)},_{(0)},_{(1)})$ where $V_0, V_1$ are 
two $k$-modules, $\b1$ an element of $V_0$ ({\it vacuum vector}), 
$\dpar:\ V_0\lra V_1$ a 
morphism of $k$-modules, 

$_{(i)}:\ (V_0\oplus V_1)\times (V_0\oplus V_1)\lra V_0\oplus V_1\ 
(i=-1,0,1)$

are $k$-bilinear operations of degree $-i-1$. Elements of $V_0$ 
(resp., $V_1$) will be denoted $a,b,c$ (resp., $x,y,z$). So, we have 
$7$ operations: $a_{(-1)}b, a_{(-1)}x, x_{(-1)}a, a_{(0)}x$, $
x_{(0)}a, x_{(0)}y$ and $x_{(1)}y$. The following axioms must be 
satisfied: 

{\it (Vacuum)} 
$$
a_{(-1)}\b1=a;\ x_{(-1)}\b1=x;\ x_{(0)}\b1=0
\eqno{(Vac)}
$$ 

{\it (Derivation)} 

$$
(\dpar a)_{(0)}b=0;\ (\dpar a)_{(0)}x=0;\ 
(\dpar a)_{(1)}x=-a_{(0)}x
\eqno{(Der)_1}
$$
$$
\dpar(a_{(-1)}b)=(\dpar a)_{(-1)}b+a_{(-1)}\dpar b;\ 
\dpar(x_{(0)}a)=x_{(0)}\dpar a
\eqno{(Der)_2}
$$

{\it (Commutativity)} 

$$
a_{(-1)}b=b_{(-1)}a;\ a_{(-1)}x=x_{(-1)}a-\dpar(x_{(0)}a)
\eqno{(Com)_{-1}}
$$
$$
x_{(0)}a=-a_{(0)}x;\ x_{(0)}y=-y_{(0)}x+\dpar(y_{(1)}x)
\eqno{(Com)_0}
$$ 
$$
x_{(1)}y=y_{(1)}x
\eqno{(Com)_1}
$$ 

{\it (Associativity)}

$$
(a_{(-1)}b)_{(-1)}c=a_{(-1)}b_{(-1)}c
\eqno{(Ass)_{-1}}
$$
{\it Operation $_{(0)}$ is a derivation with respect to all 
operations $_{(i)}$, i.e. 
$$
\alpha_{(0)}\beta_{(i)}\gamma=
(\alpha_{(0)}\beta)_{(i)}\gamma+\beta_{(i)}\alpha_{(0)}\gamma,\ (\alpha, 
\beta, \gamma\in V_0\oplus V_1)
\eqno{(Ass)_0}
$$
whenever the both sides are defined.} 

$$
(a_{(-1)}x)_{(0)}b=a_{(-1)}x_{(0)}b
\eqno{(Ass)_1}
$$ 
$$
(a_{(-1)}b)_{(-1)}x=a_{(-1)}b_{(-1)}x+(\dpar a)_{(-1)}b_{(0)}x+ 
(\dpar b)_{(-1)}a_{(0)}x
\eqno{(Ass)_2}
$$ 
$$
(a_{(-1)}x)_{(1)}y=a_{(-1)}x_{(1)}y-x_{(0)}y_{(0)}a
\eqno{(Ass)_3}
$$ 

A {\it morphism} between two $1$-truncated vertex algebras 
$f:\ v=(V_0,V_1,\ldots)
\lra v'=(V'_0,V'_1,\ldots)$ is a pair of maps of $k$-modules 
$f=(f_0,f_1),\ f_i:\ V_i\lra V'_i$ such that $f_0(\b1)=\b1',\ f_1(\dpar a)=
\dpar f_0(a)$ and $f(\alpha_{(i)}\beta)=f(\alpha)_{(i)}f(\beta)$, 
whenever both sides are defined. 

This way we get a category $\CV ert_{\leq 1}$ of $1$-truncated vertex algebras. 
We have an obvious truncation functor 
$$
t:\ \CV ert\lra \CV ert_{\leq 1}
\eqno{(3.1.1)}
$$
which assignes to a vertex algebra $V=\oplus V_i$ the truncated algebra 
$tV:=(V_0,V_1,\ldots)$. In Section 9 below we shall construct a 
left adjoint to $t$.  

{\bf 3.2.} Let $v=(V_0,V_1,\ldots)$ be a $1$-truncated vertex algebra. 
Let $\Omega_v\subset V_1$ be the $k$-submodule generated by all elements 
$a_{(-1)}\dpar b$. Set $T_v=V_1/\Omega_v$; let $\pi:\ V_1\lra T_v$ be 
the canonical projection. 

Let us call $v$ {\it splittable} if $\pi$ admits a $k$-linear 
splitting $s:\ T\lra V_1$, cf. 2.2. This is of course a weak condition; 
it holds true for example if $T$ is a projective $k$-module. Let 
$\CV ert'_{\leq 1}\subset \CV ert_{\leq 1}$ denote the full subcategory 
of splittable $1$-truncated vertex algebras. 

A $1$-truncated vertex algebra with a chosen splitting $s$ will 
be called {\it split}. The argument of Section 2 assigns a vertex algebroid 
$\CA(v;s)$ to every split $1$-truncated vertex 
algebra $(v;s)$. 

We have by definition 
$$
\CA(V;s)=\CA(tV;s)
\eqno{(3.2.2)}
$$
for every split vertex algebra $(V;s)$. 

{\bf 3.3.} Conversely, let $\CA=(A,T,\Omega,\ldots)$ be a vertex 
algebroid. We want to assign to it a split $1$-truncated vertex algebra. 
To this end one needs simply to invert the construction of the previous 
Section.  

Namely $V_0=A,\ V_1=T\oplus\Omega$; let $\dpar:\ V_0\lra V_1$ be 
the composition of $\dpar:\ A\lra\Omega$ with the obvious embedding 
$\Omega\subset V_1$; let $s:\ T\lra V_1$ be the obvious embedding. 

Let us define the operations $_{(i)}$ as follows: 
$$
a_{(-1)}b=ab;\ a_{(-1)}\omega=a\omega;\ a_{(-1)}\tau=a\tau-\gamma(a,\tau)
\eqno{(3.3.1)}
$$
$$
a_{(0)}b=a_{(0)}\omega=\omega_{(0)}\omega'=0;\ 
\tau_{(0)}a=\tau(a);\ \tau_{(0)}\omega=\tau(\omega)
\eqno{(3.3.2)}
$$
$$
\tau_{(0)}\tau'=[\tau,\tau']-c(\tau,\tau')+\frac{1}{2}\dpar\langle\tau,
\tau'\rangle
\eqno{(3.3.3)}
$$
$$
x_{(1)}y=\langle x,y\rangle
\eqno{(3.3.4)}
$$
By inverting the argument of the previous Section, one sees easily 
that we get a split $1$-truncated vertex algebra $(V_0,V_1,\ldots)$, to 
be denoted by $u\CA$. 

Let $v\in\CV ert'_{\leq 1}$. For each splitting $s$ of $v$ we have 
by construction a canonical isomorphism 
$$
v\iso u\CA(v;s)
\eqno{(3.3.5)}
$$  
For two vertex algebroids $\CA,\CA'$ define the set of morphisms between 
$\CA$ and $\CA'$ to be the subset of $Hom_{\CV ert_{\leq 1}}(u\CA,u\CA')$ 
which consists of all morphisms $g:\ u\CA\lra u\CA'$ such that 
$g(\Omega)\subset\Omega'$. 
This way we get a category of vertex algebroids, to be denoted $\CA lg$. The 
mapping $u$ induces a functor  
$$
u:\ \CA lg\lra\CV ert'_{\leq 1}
\eqno{(3.3.6)}
$$
which is in fact an equivalence of categories, due to (3.3.5). 

{\bf 3.4.}  Let $\CA=(A,T,\Omega,\ldots)$ and $\CA'=(A',T',\Omega',\ldots)$ 
be two vertex algebroids. We want to describe explicitly the set of 
morphisms $Hom_{\CA lg}(\CA,\CA')\subset Hom_{\CV ert_{\leq 1}}(u\CA,u\CA')$. 

Let $g:\ u\CA\lra u\CA'$ be a morphism belonging to $Hom_{\CA lg}(\CA,\CA')$. 
Let $g_A:\ A\lra A'$ be its component of weight $0$. Since $g$ preserves 
the operation $_{(-1)}$ and the vacuum, $g_A$ is a morphism of commutative 
$k$-algebras. Let $g_1:\ T\oplus\Omega\lra T'\oplus\Omega'$ be 
its component of weight $1$. We have $\dpar\circ g=g\circ\dpar$; therefore 
$g_1(\dpar A)\subset \dpar A'$. By definition, $g_1(\Omega)\subset\Omega'$. 
Let us denote by $g_{\Omega}:\ \Omega\lra\Omega'$ and $g_T:\ T=U\CA/\Omega
\lra U\CA'/\Omega'$ the morphisms induced by $g_1$. Thus, in components 
$g_1$ has the form 
$$
g_1(\tau,\omega)=(g_T(\tau),g_\Omega(\omega)+h(\tau))
\eqno{(3.4.1)}
$$
where $h:\ T\lra\Omega'$ is some $k$-linear mapping. 

{\bf 3.5. Theorem.} {\it The correspondence described in} 3.4 
{\it provides a canonical 
isomorphism of the set $Hom_{\CA lg}(\CA,\CA')$  
with the set of all 
quadruples $(g_A,g_T,g_\Omega,h)$ where 

$(3.5A)$ $g_A:\ A\lra A'$ is a morphism of $k$-algebras; 

$(3.5\Omega)$ $g_\Omega:\ \Omega\lra\Omega'$ is a morphism of $k$-modules 
such that $g_\Omega(\dpar a)=\dpar g_A(a)$ and $g_\Omega(a\omega)=g_A(a)
g_\Omega(\omega)$; 

$(3.5T)$ $g_T:\ T\lra T'$ is a morphism of Lie $k$-algebras such that 
$g_T(a\tau)=g_A(a)g_T(\tau)$,\ $g_A(\tau(a))=g_T(\tau)(g_A(a))$,\ 
$g_\Omega(\tau(\omega))=g_T(\tau)(g_\Omega(\omega))$ and $g_A(\langle
\tau,\omega\rangle)=\langle g_T(\tau),g_\Omega(\omega)\rangle$; 

$(3.5\gamma)$ $h:\ T\lra\Omega'$ is a morphism of $k$-modules such that 
$$
 h(a\tau)=g_A(a)h(\tau)-\gamma'(g_A(a),g_T(\tau))+g_\Omega(\gamma(a,
\tau)); 
$$

$(3.5\langle,\rangle)$ $g_A(\langle\tau_1,\tau_2\rangle)=\langle g_T(\tau_1),
g_T(\tau_2)\rangle'+\langle g_T(\tau_1),h(\tau_2)\rangle+
\langle g_T(\tau_2),h(\tau_1)\rangle$; 

$(3.5c)$\ $g_\Omega(c(\tau_1,\tau_2))=c'(g_T(\tau_1),g_T(\tau_2))+
\frac{1}{2}\dpar\langle g_T(\tau_1),h(\tau_2)\rangle-\frac{1}{2}\dpar
\langle g_T(\tau_2),h(\tau_1)\rangle-$

$-g_T(\tau_1)(h(\tau_2))+g_T(\tau_2)
(h(\tau_1))+h([\tau_1,\tau_2])$.}

{\bf 3.6.} {\it Proof.} By definition, a quadruple $g=(g_A,g_T,g_\Omega,h)$ 
defines a morphism $u\CA\lra u\CA'$ iff it satisfies the identities 
(3.6.1) through (3.6.6) below $(a,b\in A,\ x,y\in u\CA_1)$:  
$$
g_0(\b1)=\b1;\ g_0(a_{(-1)}b)=g_0(a)_{(-1)}g_0(b)
\eqno{(3.6.1)}
$$
$$
g_1(\dpar a)=\dpar g_0(a)
\eqno{(3.6.2)}
$$
$$
g_1(a_{(-1)}x)=g_0(a)_{(-1)}g_0(x)
\eqno{(3.6.3)}
$$
$$
g_1(x_{(0)}a)=g_1(x)_{(0)}g_0(a)
\eqno{(3.6.4)}
$$
$$
g_1(x_{(0)}y)=g_1(x)_{(0)}g_1(y)
\eqno{(3.6.5)}
$$
$$
g_0(x_{(1)}y)=g_1(x)_{(1)}g_1(y)
\eqno{(3.6.6)}
$$
The condition (3.6.1) is equivalent to (3.5A); (3.6.2) and (3.6.3) for 
$x=\omega$ is equivalent to $(3.5\Omega)$. Let us write down (3.6.3) 
for $x=\tau$. We have
$$
a_{(-1)}\tau=(a\tau,-\gamma(a,\tau))
\eqno{(3.6.7)}
$$
Therefore $g_1(a_{(-1)}\tau)=(g(a\tau),h(a\tau)-g(\gamma(a,\tau))$ and 
$$
g_1(a)_{(-1)}g_1(\tau)=g_1(a)_{(-1)}(g(\tau),h(\tau))=(g(a)g(\tau),
-\gamma(g(a),g(\tau))+g(a)h(\tau))
$$
It follows that (3.6.3) for $x=\tau$ is equivalent to $g(a\tau)=g(a)g(\tau)$ 
and $(3.5\gamma)$. (3.6.4) for $x=\omega$ is vacuous. (3.6.4) for 
$x=\tau$ is quivalent to $g(a\tau)=g(a)g(\tau)$. 

(3.6.6) for $x,y\in\Omega$ 
is vacuous. (3.6.6) for $x=\omega, y=\tau$ (or vice versa) is equivalent 
to $g(\langle\tau,\omega\rangle)=\langle g(\tau),g(\omega)\rangle$. We have 
$$
g_1(\tau_1)_{(1)}g_1(\tau_2)=\langle g(\tau)_1,g(\tau_2)\rangle'+
\langle g(\tau_1),h(\tau_2)\rangle+\langle g(\tau_2),h(\tau_1)\rangle
$$
Therefore (3.6.6) for $x,y\in T$ is equivalent to $(3.5\langle,\rangle)$. 

(3.6.5) for $x=\tau, y=\omega$ (or vice versa) is equivalent to 
$g(\tau(\omega))=g(\tau)(g(\omega))$. Finally, let us write down (3.6.5) 
for $x,y\in T$. We have 
$$
\tau_{1(0)}\tau_2=([\tau_1,\tau_2],-c(\tau_1,\tau_2)+\frac{1}{2}
\dpar\langle\tau_1,\tau_2\rangle)
\eqno{(3.6.8)}
$$
cf. (2.2.2) and $(2.2.3)_c$. Therefore
$$
g_1(\tau_{1(0)}\tau_2)=(g([\tau_1,\tau_2],h([\tau_1,\tau_2])-g(c(\tau_1,
\tau_2)+\frac{1}{2}\dpar\langle g(\tau_1),g(\tau_2)\rangle+
$$
$$
+\frac{1}{2}
\dpar\langle g(\tau_1),h(\tau_2)\rangle+\frac{1}{2}\dpar\langle g(\tau_2),
h(\tau_1)\rangle)
$$
where we have used $(3.5\langle,\rangle)$. On the other hand, 
$$
g_1(\tau_1)_{(0)}g_1(\tau_2)=([g(\tau_1),g(\tau_2)],-c'(g(\tau_1),g(\tau_2))+ 
\frac{1}{2}\dpar\langle g(\tau_1),g(\tau_2)\rangle+
$$
$$
+g(\tau_1)(h(\tau_2))-
g(\tau_2)(h(\tau_1))+\dpar\langle h(\tau_1),g(\tau_2)\rangle) 
$$
where we have used the formula 
$$
\omega_{(0)}\tau=-\tau(\omega)+\dpar\langle\tau,\omega\rangle
\eqno{(3.6.9)}
$$
following from (0.5.9) (cf. (2.3.5)). Therefore (3.6.5) for $x,y\in T$ 
is equivalent to the requirement that $g_T$ is a morphism of Lie 
algebras, and to $(3.5c)$. The theorem is proved. $\btu$ 

{\bf 3.7.} Let $(A,T,\ldots)\buildrel{g}\over\lra (A',T',\ldots)
\buildrel{g'}\over\lra (A'',T'',\ldots)$ be two morphisms of 
vertex algebroids. Then their composition is obviously equal to 
$$
g'\circ g=(g'_A g_A, g'_T g_T, g'_\Omega g_\Omega,g'_\Omega h+h'g_T)
\eqno{(3.7.1)}
$$
cf. (3.4.1). The identity morphisms are 
$$
Id_\CA=(Id_A,Id_T,Id_\Omega,0)
\eqno{(3.7.2)}
$$
  
{\bf 3.8. Theorem.} {\it (Extension of Mappings.)}\ {\it In the conditions of} 
(3.5), {\it let $T_0\subset T$ 
be a $k$-submodule which generates $T$ as an $A$-module. Let $g=(g_A,g_T,
g_\Omega,h)$ be such a quadruple that $(3.5A), (3.5\Omega)$ and $(3.5T)$ 
are fulfilled and such that $(3.5\gamma), (3.5\langle,\rangle)$ and  
$(3.5c)$ are fulfilled for all $\tau, \tau_i\in T_0$. 

Then $(3.5\gamma), (3.5\langle,\rangle)$ and $(3.5c)$ hold true for all 
$\tau, \tau_i\in T$, i.e. $g$ defines a morphism $\CA\lra \CA'$.} 

Cf. Theorem 1.9. 

{\it Proof.} Let us prove that if $(3.5\gamma)$ is true $(a,\tau)$ 
with a fixed $\tau$ and all $a$, then it is true for all couples 
$(a,b\tau)$. We have to prove that
$$
h(ab\tau)=g(a)h(b\tau)-\gamma(g(a),g(b\tau))+g(\gamma(a,b\tau))
\eqno{(?)}
$$
The left hand side is equal to 
$$
h(ab\tau)=g(ab)h(\tau)-\gamma(g(ab),g(\tau))+g(\gamma(ab,\tau))=
$$
$$  
=g(ab)h(\tau)-\gamma(g(a),g(b)g(\tau))-g(a)\gamma(g(b),g(\tau))-
g(\tau)(g(a))\dpar g(b)-g(\tau)(g(b))\dpar g(a)+
$$
$$
+g(\gamma(a,b\tau))+g(a\gamma(b,\tau))+g(\tau(a)\dpar b)+
g(\tau(b))\dpar a)
$$
where we have used (A1). On the other hand, the first summand in the rhs 
is equal to 
$$
g(a)h(b\tau)=g(a)g(b)h(\tau)-g(a)\gamma(g(b),g(\tau))+g(a)g(\gamma(b,\tau)), 
$$
and we easily see that we have the required identity indeed. 

The similar claims connected with the equations $(3.5\langle,\rangle)$ 
and $(3.5c)$ are proved analogously, and we leave them to 
the reader. $\btu$ 

\bigskip\bigskip


\newpage 

\centerline{\bf \S 4. Cofibered Structure} 

\bigskip\bigskip

{\bf 4.1.} Let us introduce objects which lie "in between" extended 
Lie algebroids and vertex algebroid. Let us define a 
{\it vertex prealgebroid} to be a sextuple $\CB=(A,T,\Omega,\dpar,\gamma,
\langle,\rangle)$ where $A, T,\Omega,\dpar,\gamma,\langle,\rangle$ are 
as in 1.4. 

The data $(A,T,\Omega,\langle,\rangle|_{T\times\Omega})$ should 
form an extended Lie algebroids and the axioms 1.4 (A1) and (A2) must 
be satisfied. 

Let us define a {\it morphism} between two vertex prealgebroids 
$\CB=(A,T,\ldots)$ and $\CB'=(A',T',\ldots)$ to be a quadruple $g=(g_A,g_T,
g_\Omega,h)$ as in Theorem 3.5, satisfying the properties 
$(3.5A)$ --- $(3.5\langle,\rangle)$. The composition of morphisms 
and the identity morphisms are defined by the rules (3.7.1), (3.7.2).  
This way we get a category 
$\CP reAlg$ of vertex prealgebroids.    

We have obvious forgetful functors
$$
\CA lg\buildrel{P}\over\lra\CP reAlg\buildrel{Q}\over\lra\CL ieAlg
\eqno{(4.1.1)}
$$
Note that $P$ is injective on $Hom$'s. 
The composition $Q\circ P:\ \CA lg\lra \CL ieAlg$ will be denoted by 
$\Theta$. 

We will use the standard notation for fibers. For example, if 
$\CB\in \CP reAlg$ then $\CA lg_\CB$ will denote the category 
whose objects are $\CA\in \CA lg$ such that $P(\CA)=\CB$ and 
morphisms are morphisms in $\CA$ which go to $Id_\CB$ after 
applying $P$, and similarly with the other fibers. 

{\bf 4.2.} Let  
$\CB=(A,T,\Omega,\dpar,\gamma,\langle,\rangle)\in\CP reAlg$; 
let $\CT=Q(\CB)=(A,T,\Omega,\langle,\rangle|_{T\times\Omega})$. 
Let $\tAlg_\CB$ be a full subcategory of $\CA lg_\CT$ with 
$Ob\tAlg_\CB=Ob\CA lg_\CB$. In other words, $Hom_{\tAlg_\CB}(\CA,\CA')$ 
consists of quadruples of the form $g=(Id_A,Id_T,Id_\Omega,h)$.
    
Note that all such $g$ are invertible, 
namely $g^{-1}=(Id_A,Id_T,\Id_\Omega,-h)$, cf. (3.7.1), (3.7.2). 
In other words, $\tAlg_\CB$ is a {\it groupoid}.  
    
Consider the de Rham-Chevalley complex of $\CT$: 
$$
\Omega^\cdot(\CT):\ 0\lra A\lra\Omega\buildrel{d_{DR}}\over{\lra}
\Omega^2\lra\ldots
$$
cf. 1.3. Let $\CA_i=(...,c_i)\in Ob \tAlg_\CB,\ i=1,2$. Set $\alpha=c_1-c_2$. 
It follows from (A4) that 
$$
\langle\tau_2,\alpha(\tau_1,\tau_3)\rangle+\langle\tau_3,\alpha(\tau_1,
\tau_2)\rangle=0,
$$
i.e. $\alpha\in\Omega^3$. It follows from (A5) that $\alpha\in\Omega^{3,cl}
:=Ker(\Omega^3\lra\Omega^4)$, cf. (1.3.2). 

Conversely, if $\CA=(\ldots,c)\in Ob \tAlg_\CB$ and $\alpha\in\Omega^{3,cl}$ 
then 
$\CA\dplus\alpha:=(\ldots,c+\alpha)\in Ob \CA lg_\CB$. We have proven that 

{\bf 4.2.1.} {\it The set $Ob \tAlg_\CB$ is canonically an 
$\Omega^{3,cl}$-torseur.} 

Let $\CA\in Ob \tAlg_\CB$,  
$\CA'=\CA\dplus\alpha,\ \alpha\in\Omega^{3,cl}$, $g=(Id_A,Id_T,
Id_\Omega,h)\in 
Hom_{\tAlg_\CB}(\CA,\CA')$. Then, due to $(3.5\gamma)$, $h$ must be 
$A$-linear; by $(3.5\langle,\rangle)$, $h\in\Omega^2$ and by 
$(3.5c)$, $d_{DR}h=-\alpha$. Therefore, 

{\bf 4.2.2.} {\it the set $Hom_{\tAlg_\CB}(\CA,\CA\dplus\alpha)$ may be 
canonically identified with the set of $h\in\Omega^2$ such that 
$d_{DR}h=-\alpha$. 

Consequently, the set of isomorphism classes of objects 
of the groupoid $\tAlg_\CB$ is equal to the third de Rham cohomology 
$H^3(\Omega^\cdot(\CT))$.} 

{\bf 4.3.} Let $\CT=(A,T,\ldots), \CT'=(A',T',\ldots)\in\CL ieAlg$, 
$g=(g_A,g_T,g_\Omega)\in Hom_{\CL ie Alg}(\CT',\CT)$. Let  
$\CB\in \CP re Alg_\CT, \CB'\in \CP reAlg_{\CT'}$. 

Let us consider the set $Hom_g(\CB',\CB)\subset Hom_{\CP reAlg}(\CB',\CB)$ 
consisting of all morphisms $\tg$ such that $Q(\tg)=g$. Let 
$\tg=(g_A,g_T,g_\Omega,h), \tg'=(g_A,g_T,g_\Omega,h')\in Hom_g(\CB',\CB)$. 
Set $\beta:=h-h':\ T\lra\Omega'$.  
Due to $(3.5\gamma)$ and $(3.5\langle,\rangle)$, $\beta$ will satisfy 
the  properties
$$
\beta(a\tau)=g_A(a)\beta(\tau)
\eqno{(4.3.1)}
$$
and 
$$
\langle g_T(\tau_1),\beta(\tau_2)\rangle+
\langle g_T(\tau_2),\beta(\tau_1)\rangle=0
\eqno{(4.3.2)}
$$            
Let us denote by $\Omega^2_g$ the set of all maps 
$\beta:\ T'\lra \Omega$ satisfying (4.3.1) and (4.3.2). 
It is an $A$-module in the obvious way, and consequently an $A'$-module,  
by restriction of scalars.  
In particular, $\Omega^2_{Id_\CT}=\Omega^2(\CT)$. 

Vice versa, for each $\beta\in \Omega^2_g$, the map 
$\tg'\dplus\beta:=(g_A,g_T,g_\Omega,h+\beta)$ belongs to $Hom_g(\CB',\CB)$. 
Thus we have proven that 

{\bf 4.3.1.} {\it the set $Hom_g(\CB',\CB)$ is canonically an 
$\Omega^2_g$-torseur. 
In particular, if $\CT=\CT'$ then $Hom_{\CP reAlg_\CT}(\CB',\CB)$ is 
an $\Omega^2(\CT)$-torseur.} 

Similarly, let $\CA\in\CA lg_\CT,\ \CA'\in\CA lg_{\CT'}$. Let $Hom_g(\CA',
\CA)$ be the subset of $Hom_{\CA lg}(\CA',\CA)$ consisting of all morphisms 
$\tg$ such that $\Theta(\tg)=g$. Let $\Omega^{2,cl}_g$ be  
the $k$-module of all maps $\beta:\ T'\lra\Omega$ satisfying 
(4.3.1), (4.3.2) and 
$$
-\beta([\tau_1,\tau_2])+g_T(\tau_1)(h(\tau_2))-g_T(\tau_2)(h(\tau_1))-
\dpar\langle g_T(\tau_1),h(\tau_2)\rangle=0
\eqno{(4.3.3)}
$$
In particular, $\Omega^{2,cl}_{Id_\CT}=\Omega^{2,cl}(\CT)$, cf. (1.3.4). 
As above, due to (3.5c), we get 

{\bf 4.3.2.} {\it the set $Hom_g(\CA',\CA)$ is canonically an 
$\Omega^{2,cl}_g$-torseur. In particular, if $\CT=\CT'$ then 
$Hom_{\CA lg_\CT}(\CA',\CA)$ is an $\Omega^{2,cl}(\CT)$-torseur.} 

If $g':\ \CT''\lra\CT'$ is another morphism of extended Lie algebroids, 
we have a composition map
$$
\nu_{g,g'}:\ \Omega^2_g\times\Omega^2_{g'}\lra\Omega^2_{g g'},\ \ 
\nu_{g,g'}(\beta,\beta'):=g_\Omega\beta'+\beta g'_T
\eqno{(4.3.4)}
$$
It is a morphism of $A'$-modules. We have $\nu_{g,g'}(\Omega^{2,cl}_g\times
\Omega^{2,cl}_{g'})\subset\Omega^{2,cl}_{gg'}$ (cf. Remark 4.4 below). 
The maps (4.3.4) are associative in the obvious sense (with 
respect to triples of composable morphisms in $\CL ieAlg$). 

On the other hand, given $\CB''\in\CP reAlg_{\CT''}$, we have 
the composition 
$$
Hom_g(\CB',\CB)\times Hom_{g'}(\CB'',\CB')\lra 
Hom_{gg'}(\CB'',\CB)
\eqno{(4.3.5)}
$$

{\bf 4.3.3.} {\it The map} (4.3.5) {\it is compatible with} (4.3.4). 
{\it Therefore, we have a canonical isomorphism of $\Omega^2_{gg'}$-torseurs
$$
\rho_{g,g'}:\ Hom_{gg'}(\CB'',\CB)\iso\nu_{g,g'*}(Hom_g(\CB',\CB)\times
Hom_{g'}(\CB'',\CB))
\eqno{(4.3.6)}
$$   
\it The isomorphisms $\rho_{g,g'}$  
satisfy obvious $2$-cocycle equations connected with  
triples of composable morphisms $(g,g',g'')$.}

This follows from (3.7.1).

{\bf 4.4. Remark.} Generalizing 1.3, one can define the de Rham-Chevalley 
complexes $\Omega^\cdot_g$ such that the modules $\Omega^2_g$ and 
$\Omega^{2,cl}_g$ from the previous no. become really the module 
of two-forms and closed ones respectively. These complexes come 
equipped with the composition maps $\Omega^\cdot_g\times\Omega^\cdot_{g'}\lra
\Omega^\cdot_{gg'}$ which are compatible with the de Rham differentials 
and satisfy associativity. 

We leave the necessary definitions as an exercise to the reader.

{\bf 4.5. Theorem.} {\it Let $g=(g_A,g_T,g_\Omega,h):\ \CB\lra\CB'$ 
be a morphism between vertex prealgebroids such that $g_T$ is an 
isomorphism. Let $\CA\in \CA lg_\CB$. There exists a unique pair 
$(g_*\CA,\tg)$ where $g_*\CA\in\CA lg_{\CB'}$ and $\tg:\ \CA\lra 
g_*\CA$ is a morphism of vertex algebroids such that $P(\tg)=\tg$.} 

{\it Proof.} {\it Uniqueness.} Let $\CA=(A,T,\Omega,\dpar,\gamma,
\langle,\rangle,c)$;
\newline let $\CA'=(A',T',\Omega',\dpar',\langle,\rangle',
\gamma',c')\in\CA lg_{\CB'}$. A morphism $\tg:\ \CA\lra\CA'$ such 
that $P(\tg)=g$ must of course be  
represented by the same quadruple as $g$, i.e. 
$\tg=(g_A,g_T,g_\Omega,h)$. It is a morphism of vertex algebroids iff 
the condition $(3.5c)$ is fulfilled, i.e. iff 
$$
c'(g_T(\tau_1),g_T(\tau_2))=g_\Omega(c(\tau_1,\tau_2))-\frac{1}{2}
\dpar'\langle g_T(\tau_1),h(\tau_2)\rangle'+\frac{1}{2}
\dpar'\langle g_T(\tau_2),h(\tau_1)\rangle'+
$$
$$
+g_T(\tau_1)(h(\tau_2))-
g_T(\tau_2)(h(\tau_1))-h([\tau_1,\tau_2])
\eqno{(4.5.1)}
$$  
This equation defines $c'$ uniquely since $g_T$ is an isomorphism 
by assumption. 

{\it Existence.} We have to check that the map $c'$ defined by (4.5.1) 
satisfies axioms (A3), (A4), (A5). Let us check (A4) for example. 

To unburden the notation, let us assume that $Q(\CB)=Q(\CB')$ and $g=(Id_A,
Id_T,Id_\Omega,h)$ (the general case being treated by an identical 
computation). We have to prove that
$$
\langle\tau_2,c'(\tau_1,\tau_3)\rangle+\langle\tau_3,c'(\tau_1,\tau_2)
\rangle=
$$ 
$$
=\langle [\tau_1,\tau_2],\tau_3\rangle'+\langle [\tau_1,\tau_3],\tau_2
\rangle'-\tau_1(\langle\tau_2,\tau_3\rangle')+\frac{1}{2}\tau_2(
\langle\tau_1,\tau_3)\rangle')+\frac{1}{2}\tau_3(\langle\tau_1,\tau_2
\rangle')
\eqno{(?)}
$$
We have 
$$
\langle\tau_2,c'(\tau_1,\tau_3)\rangle=\langle\tau_2,c(\tau_1,\tau_3)-
\frac{1}{2}\dpar\langle\tau_1,h(\tau_3)\tau_3)\rangle+\frac{1}{2}
\dpar\langle\tau_3,h(\tau_1)\rangle+\tau_1(h(\tau_3))-
$$
$$
-\tau_3(h(\tau_1))-h([\tau_1,\tau_3])\rangle=\langle\tau_2,c(\tau_1,\tau_3)
\rangle-\frac{1}{2}\tau_2(\langle\tau_1,h(\tau_3)\rangle)+\frac{1}{2}
\tau_2(\langle\tau_3,h(\tau_1)\rangle)+
$$
$$
+\tau_1(\langle\tau_2,h(\tau_3)\rangle)-
\langle [\tau_1,\tau_2],h(\tau_3)\rangle
-\tau_3(\langle\tau_2,h(\tau_1)\rangle)-
\langle [\tau_3,\tau_2],h(\tau_1)\rangle-\langle\tau_2,h([\tau_1,\tau_3])
\rangle  
$$
Similarly, interchanging $\tau_2$ with $\tau_3$, 
$$
\langle\tau_3,c'(\tau_1,\tau_2)\rangle=
\langle\tau_3,c(\tau_1,\tau_2)
\rangle-\frac{1}{2}\tau_3(\langle\tau_1,h(\tau_2)\rangle)+\frac{1}{2}
\tau_3(\langle\tau_2,h(\tau_1)\rangle)+
$$
$$
+\tau_1(\langle\tau_3,h(\tau_2)\rangle)-
\langle [\tau_1,\tau_3],h(\tau_2)\rangle
-\tau_2(\langle\tau_3,h(\tau_1)\rangle)-
\langle [\tau_2,\tau_3],h(\tau_1)\rangle-\langle\tau_3,h([\tau_1,\tau_2])
\rangle
$$
By assumption,
$$
\langle\tau_2,c(\tau_1,\tau_3)\rangle+\langle\tau_3,c(\tau_1,\tau_2)
\rangle=
$$ 
$$
=\langle [\tau_1,\tau_2],\tau_3\rangle+\langle [\tau_1,\tau_3],\tau_2
\rangle-\tau_1(\langle\tau_2,\tau_3\rangle)+\frac{1}{2}\tau_2(
\langle\tau_1,\tau_3)\rangle)+\frac{1}{2}\tau_3(\langle\tau_1,\tau_2
\rangle)
$$
Using this and the axiom $(3.5\langle,\rangle)$ which takes the form 
$$
\langle\tau_1,\tau_2\rangle'=\langle\tau_1,\tau_2\rangle
-\langle\tau_1,h(\tau_2)\rangle-\langle\tau_2,h(\tau_1)\rangle
$$
we get the required identity (?). 

The other two axioms, (A3) and (A5), are checked in a similar manner, 
and we leave them to the reader. It is convenient to 
check (A5) in its equivalent form $(A5)^{ter}$, see 1.6. $\btu$ 

\bigskip\bigskip


\centerline{\bf \S 5. Chern-Simons} 

\bigskip\bigskip

{\bf 5.1.} Let us define a {\it frame} of a Lie $A$-algebroid $T$ to be  
a $k$-submodule 
$\fg\subset T$ such that $A\otimes_k\fg=T$. For example, if $T$ is a free 
$A$-module and $\{\tau_i\}$ is some $A$-base of $T$ then $\fg_{\{\tau_i\}}:=
\sum k\tau_i$ is a frame in $T$.
 
A frame of a vertex $A$-algebroid (resp. prealgebroid, extended Lie 
algebroid) is by definition a frame of the underlying  
Lie algebroid $T$.  
A vertex algebroid (resp. prealgebroid, ...) equipped 
with a frame will be called {\it framed}. 

We will call an extended Lie algebroid $\CT$ {\it quasiregular} if it 
is perfect (see 1.2) and admits a frame.     
 
{\bf 5.2.} Let us consider the situation of 4.3, 
and assume that $\CT'$ is equipped with a frame $\fg'$ and $\CT$ is perfect. 
Let us assume that $g_T$ is an isomorphism. Then $\fg:=g_T(\fg')$ is 
a frame in $T$. Let us define a $k$-linear map $h_{\fg'}:\ \fg'\lra\Omega$ 
by the condition 
$$
\langle g_T(\tau_1),h_{\fg'}(\tau_2)\rangle=
\frac{1}{2}\bigl\{ g_A(\langle\tau_1,\tau_2\rangle')-
\langle g_T(\tau_1),g_T(\tau_2)\rangle\}
\ \ (\tau_i\in\fg')
\eqno{(5.2.1)}
$$
Since $g_T$ is an isomorphism, (5.2.1) defines $h_{\fg'}$ uniquely. 
Then the condition $(3.5\langle,\rangle)$ is obviously satisfied 
for all $\tau_i\in\fg'$. There is a unique extension of $h_{\fg'}$ to a 
map $h_{\fg'}:\ T'\lra\Omega$ satisfying $(3.5\gamma)$ for all $\tau\in\fg'$. 

By Theorem 3.8, the conditions $(3.5\gamma)$ and $(3.5\langle,\rangle)$ 
are then fulfilled for all $\tau, \tau_i\in T$, that is, 
$g_{\fg'}:=(g,h_{\fg'})$ is a morphism of vertex prealgebroids $\CB'\lra\CB$. 

In other words, 

{\bf 5.2.1.} {\it the $\Omega^2(\CT)$-torseur $Hom_g(\CB',\CB)$ is 
equipped with a trivialization $g_{\fg'}$.}

{\bf 5.3.} Let $(\CT,\fg)$ be a framed extended Lie algebroid. 
We can define a vertex algebroid 
$\CA_{\CT,\fg}$ by setting $\gamma(a,\tau)=0,\ \langle\tau,\tau'
\rangle=c(\tau,\tau')=0$ for all $a\in A,\ \tau,\tau'\in\fg$ and then extending 
the operations $\gamma,\ \langle,\rangle,\ c$ to the whole of $T$ using 
the axioms (A1) --- (A3). By Theorem 1.9 we get a vertex algebroid, 
to be denoted $\CA_{\CT;\fg}$. 

If $\fg$ is a Lie subalgebra of $T$ then $\CA_{\CT;\fg}$ may be 
defined as follows. 
Take the vertex $k$-algebroid $\CA_\fg=(k,\fg,\fg^*,0,0,\langle,\rangle,0)$ 
where $\fg$ acts on $\fg^*$ in the coadjoint way, 
$\langle,\rangle|_{\fg\times\fg^*}$ is the obvious pairing, 
the other components of $\langle,\rangle$ being zero. Then apply 
to $\CA_\fg$ the pushout with respect to the structure morphism $k\lra A$, 
where we 
set $\gamma(a,\tau)=0$ 
for $\tau\in\fg$ and the Lie algebra acts on $A$ due to the embedding 
$\fg\subset T$, cf. 1.10. We get a vertex $A$-algebroid which is equal to 
$\CA_{\CT;\fg}$. 

We set $\CB_{\CT;\fg}:=P(\CA_{\CT;\fg})$. 

Let us 
call the frame $\fg$ {\it abelian} if $\fg$ is an abelian Lie 
subalgebra of $T$, i.e. $[\tau,\tau']=0$ for all $\tau, \tau'\in\fg$. 
We call an extended Lie algebroid $\CT$ {\it regular} if it is 
perfect and admits an abelian frame.   

{\bf 5.4.} Let $\CT=(A,T,\Omega,\ldots)$ be a regular extended Lie algebroid.  
Let $\fg, \fg'$ be two abelian frames of $\CT$. Let us assume that 
$\fg$ and $\fg'$ are free $k$-modules of finite rank $n$.   
In the sequel we will 
need some formulas pertaining to this situation. 

Choose some $k$-bases $\{\tau_i\},\ \{\tau'_i\},\ i=1,\ldots,n,$ of $\fg$ and 
$\fg'$ respectively; let $\{\omega_i\}
\subset\fg^*,\ \{\omega'_i\}\subset\Omega'$ be the dual bases. Note that 
$$
\tau_i(\omega_j)=0\ \text{for all\ }i,j, 
\eqno{(5.4.1)}
$$
since for all $p$, $\langle\tau_p,\tau_i(\omega_j)\rangle=\tau_i(\langle
\tau_p,\omega_j\rangle)-\langle [\tau_i,\tau_p],\omega_j\rangle=0$ 
because $\fg$ is abelian and $\langle\tau_p,\omega_j\rangle=\delta_{pj}$. 

Recall that 
$$
\langle\tau,a\dpar b\rangle=a\tau(b)
\eqno{(5.4.2)}
$$
and 
$$
(a\tau)(\omega)=a\tau(\omega)+\langle\tau,\omega\rangle\dpar a
\eqno{(5.4.3)}
$$
Since $\langle\tau_i,\dpar a\rangle=\tau_i(a)$, 
$$
\dpar a=\tau_i(a)\omega_i
\eqno{(5.4.4)}
$$
where we always imply the summation over repeating indices.  

Define the matrices $\phi=(\phi^{ij}), \rho=(\rho^{ij})\in GL_n(A)$ 
by $\tau_i'=\phi^{ij}\tau_j,\ \omega_i'=\rho^{ij}\omega_j$.  
Since $\langle \tau_i,\omega_j\rangle=
\langle \tau_i',\omega_j'\rangle=\delta_{ij}$, we have 
$$
\phi\cdot\rho^t=1
\eqno{(5.4.5)}
$$
where $(\rho^t)^{ij}=\rho^{ij}$. Since $[\tau_i,\tau_j]=
[\tau'_i,\tau'_j]=0$, 
$$
\phi^{ip}\tau_p(\phi^{jq})=\phi^{jp}\tau_p(\phi^{iq})
\eqno{(5.4.6)}
$$
for all $i,j,q$. Applying $\tau_r$ to (5.4.6) we get 
$$
\phi^{ip}\tau_r\tau_p(\phi^{jq})-\phi^{jp}\tau_r\tau_p(\phi^{iq})=
\tau_r(\phi^{jp})\tau_p(\phi^{iq})-\tau_r(\phi^{ip})\tau_p(\phi^{jq})
\eqno{(5.4.7)}
$$
for all $i,j,q,r$. Setting here $r=q$ and summing up by $q$ we 
get 
$$
\phi^{ip}\tau_q\tau_p(\phi^{jq})=\phi^{jq}\tau_p\tau_q(\phi^{ip})
\eqno{(5.4.8)}
$$
for all $i,j$. Applying $\tau_r$ we get 
$$
\phi^{ip}\tau_r\tau_q\tau_p(\phi^{jq})-
\phi^{jq}\tau_r\tau_p\tau_q(\phi^{ip})=
\tau_r(\phi^{jq})\tau_p\tau_q(\phi^{ip})-
\tau_r(\phi^{ip})\tau_q\tau_p(\phi^{jq})
\eqno{(5.4.9)}
$$
It follows from (5.4.5) that 
$$
\tau_r(\rho)=-\rho\tau_r(\phi^t)\rho;\ 
\tau_r(\phi)=-\phi\tau_r(\rho^t)\phi
\eqno{(5.4.10)}
$$
for all $r$. Multiplying (5.4.6) by $\rho^{ia}$ and summing up by $i$, 
we get
$$
\tau_a(\phi^{jq})=\phi^{jp}\rho^{ia}\tau_p(\phi^{iq})
\eqno{(5.4.11)}
$$
whence
$$
\rho^{ib}\tau_a(\phi^{iq})=\rho^{ia}\tau_b(\phi^{iq})
\eqno{(5.4.12)}
$$    
By (5.4.10) 
$$
\tau_c(\rho^{ub})=-\rho^{uk}\tau_c(\phi^{lk})\rho^{lb}
\buildrel{(5.4.11)}\over{=}
$$
$$
=-\rho^{uk}\phi^{lp}\rho^{ic}\tau_p(\phi^{ik})\rho^{lb}=
-\rho^{uk}\tau_b(\phi^{ik})\rho^{ic}=\tau_b(\rho^{uc})
$$
Thus, 
$$
\tau_a(\rho^{bc})=\tau_c(\rho^{ba})
\eqno{(5.4.13)}
$$

{\bf 5.5.} In the situation 5.4, consider the vertex prealgebroids 
$\CB_{\CT;\fg}=(A,T,\Omega,\dpar,\gamma,\langle,\rangle)$ and 
$\CB_{\CT;\fg'}=(A,T,\Omega,\dpar,\gamma',\langle,\rangle')$. 
According to 5.2, we have an isomorphism 
$$
g_{\fg,\fg'}=(Id_A,Id_T,Id_\Omega,h_{\fg,\fg'}):\ \CB_{\CT;\fg'}\iso
\CB_{\CT;\fg}
\eqno{(5.5.1)}
$$
where $h=h_{\fg;\fg'}$ is defined by
$$
\langle\tau'_i,h(\tau'_j)\rangle=-\frac{1}{2}\langle\tau'_i,\tau'_j\rangle
\eqno{(5.5.2)}
$$
which is (5.2.1) in our situation. 

Now consider the vertex algebroids $\CA_{\CT;\fg}=(\ldots,c),\ 
\CA_{\CT;\fg'}=(\ldots,c')$. We have $g_*\CA_{\CT;\fg'},\ 
\CA_{\CT;\fg}\in\CA lg_{\CB_{\CT;\fg}}$. Therefore, by 
4.2.1 $\CA_{\CT;\fg}=g_*\CA_{\CT;\fg'}\dplus \beta$ for some closed 
$3$-form $\beta$.    

Let us define a form $\beta=\beta_{\fg,\fg'}\in\Omega^3(\CT)$ by 
$$
\beta(\tau'_i,\tau'_j)=\beta^{ijr}\omega'_r=
\frac{1}{2}\bigl\{\tau_u(\phi^{jp})\tau_p(\phi^{iq})\tau_q(\phi^{ru})-
\tau_u(\phi^{ip})\tau_p(\phi^{jq})\tau_q(\phi^{ru})\bigr\}\omega'_r
\eqno{(5.5.3)}
$$

{\bf 5.6. Magic Lemma.} {\it The form $\beta_{\fg,\fg'}$ is closed and  
we have $\CA_{\CT;\fg}=g_*\CA_{\CT;\fg'}\dplus 
\beta_{\fg,\fg'}$.}

{\bf 5.7.} {\it Proof.}  Let us write down the things explicitly. 
The following formulas hold true in the algebroid $\CA_{\CT;\fg}$ 
(they are (1.8.3) in our situation):
$$
\gamma(a,b\tau_i)=-\tau_i(a)\dpar b-\tau_i(b)\dpar a
\eqno{(5.7.1)_\gamma}
$$
$$
\langle a\tau_i,b\tau_j\rangle=-b\tau_i\tau_j(a)-a\tau_j\tau_i(b)-
\tau_i(b)\tau_j(a)
\eqno{(5.7.1)_{\langle,\rangle}}
$$
$$
c(a\tau_i,b\tau_j)=\frac{1}{2}\{\tau_i(b)\dpar\tau_j(a)-\tau_j(a)\dpar
\tau_i(b)\}+\frac{1}{2}\dpar\{b\tau_i\tau_j(a)-a\tau_j\tau_i(b)\}
\eqno{(5.7.1)_c}
$$
Define the matrix $(h^{ij})\in Mat_n(A)$ by $h(\tau'_i)=h^{ij}\omega_j$. 
The left hand side of (5.5.2) is equal to 
$$
\langle\phi^{ip}\tau_p,h^{jq}\omega_q\rangle=\phi^{ip}h^{jp}
$$
By $(5.7.1)_{\langle,\rangle}$, the right hand side is equal to 
$$
\frac{1}{2}\bigl\{\phi^{jq}\tau_p\tau_q(\phi^{ip})+
\phi^{ip}\tau_q\tau_p(\phi^{jq})+\tau_p(\phi^{jq})\tau_q(\phi^{ip})\}
\buildrel{(5.4.8)}\over{=}
$$
$$
=\phi^{ip}\tau_q\tau_p(\phi^{jq})+\frac{1}{2}\tau_p(\phi^{jq})
\tau_q(\phi^{ip})
$$
Thus, the equation (5.5.2) takes the form
$$
\phi^{ip}h^{jp}=
\phi^{ip}\tau_q\tau_p(\phi^{jq})+\frac{1}{2}\tau_p(\phi^{jq})
\tau_q(\phi^{ip})
$$
wherefrom, applying (5.4.5), 
$$
h^{ij}=\tau_p\tau_j(\phi^{ip})+\frac{1}{2}\tau_q(\phi^{ip})\tau_p(\phi^{rq})
\rho^{rj}
\eqno{(5.7.2)}
$$
We have to prove that 
$$
c(\tau'_i,\tau'_j)=g_*c'(\tau'_i,\tau'_j)+\beta(\tau'_i,\tau'_j)
$$
where $g_*c'$ is defined by 
$$
g_*c'(\tau'_i,\tau'_j)=c'(\tau'_i,\tau'_j)+\tau'_i(h(\tau'_j))-
\tau'_j(h(\tau'_i)),
$$
by (4.5.1). Thus, we have to prove that
$$
c(\tau'_i,\tau'_j)-c'(\tau'_i,\tau'_j)
-\tau'_i(h(\tau'_j))+
\tau'_j(h(\tau'_i))=\beta(\tau'_i,\tau'_j)
\eqno{(5.7.3)}
$$ 
By definition, $c'(\tau'_i,\tau'_j)=0$. 
By $(5.7.1)_c$ and (5.4.8) we have
$$
c(\tau'_i,\tau'_j)=c(\phi^{ip}\tau_p,\phi^{jq}\tau_q)=
\frac{1}{2}\bigl\{\tau_p(\phi^{jq})\dpar\tau_q(\phi^{ip})-
\tau_q(\phi^{ip})\dpar\tau_p(\phi^{jq})\bigr\}
\eqno{(5.7.4)}
$$
On the other hand, by (5.4.1) and (5.4.3), 
$$
\tau'_i(h(\tau'_j))=(\phi^{ip}\tau_p)(h^{jq}\omega_q)=
\phi^{ip}\tau_p(h^{jq}\omega_q)+\langle\tau_p,h^{jq}\omega_q\rangle
\dpar\phi^{ip}=\phi^{ip}\tau_p(h^{jq})\omega_q+h^{jp}\dpar\phi^{ip}
$$
where we have used (5.4.1). Thus, (5.7.3) takes the form 
$$
\frac{1}{2}\bigl\{\tau_p(\phi^{jq})\dpar\tau_q(\phi^{ip})-
\tau_q(\phi^{ip})\dpar\tau_p(\phi^{jq})\bigr\}-
$$
$$
-\phi^{ip}\tau_p(h^{jq})\omega_q-h^{jp}\dpar\phi^{ip}+
\phi^{jp}\tau_p(h^{iq})\omega_q+h^{ip}\dpar\phi^{jp}=\beta^{ijr}\omega'_r
\eqno{(5.7.5)}
$$
We have to prove that the matrix $(h^{ij})$ defined by (5.7.2) satisfies 
the differential equation (5.7.5). Using (5.4.4), rewrite (5.7.5) 
in the form 
$$
\frac{1}{2}\bigl\{\tau_p(\phi^{jq})\tau_l\tau_q(\phi^{ip})-
\tau_q(\phi^{ip})\tau_l\tau_p(\phi^{jq})\bigr\}-
$$
$$
-\phi^{ip}\tau_p(h^{jl})-h^{jp}\tau_l(\phi^{ip})
+\phi^{jp}\tau_p(h^{il})+h^{ip}\tau_l(\phi^{jp})=\beta^{ijr}\rho^{rl}\omega_l
\eqno{(5.7.6)}
$$
Denote
$$
A=\frac{1}{2}\bigl\{\tau_p(\phi^{jq})\tau_l\tau_q(\phi^{ip})-
\tau_q(\phi^{ip})\tau_l\tau_p(\phi^{jq})\bigr\},
$$
$$
B=B^{ij}=-\phi^{ip}\tau_p(h^{jl})-h^{jp}\tau_l(\phi^{ip}),
$$
and
$$
C=-B^{ji}=\phi^{jp}\tau_p(h^{il})+h^{ip}\tau_l(\phi^{jp})
$$
We have to prove that $A+B+C=\beta^{ijr}\rho^{rl}\omega_l$ where $h^{ij}$ 
is given by (5.7.2) and $\beta^{ijr}$ is given by (5.5.3). 
Thus, we have 
$$
B=-\phi^{ip}\tau_p\bigl\{\tau_q\tau_l(\phi^{jq})+
\frac{1}{2}\tau_u(\phi^{jq})\tau_q(\phi^{ru})\rho^{rl}\bigr\}
-\tau_l(\phi^{ip})\bigl\{\tau_q\tau_p(\phi^{jq})+
\frac{1}{2}\tau_u(\phi^{jq})\tau_q(\phi^{ru})\rho^{rp}\bigr\}
$$
and
$$
C=\phi^{jp}\tau_p\bigl\{\tau_q\tau_l(\phi^{iq})+
\frac{1}{2}\tau_u(\phi^{iq})\tau_q(\phi^{ru})\rho^{rl}\bigr\}
+\tau_l(\phi^{jp})\bigl\{\tau_q\tau_p(\phi^{iq})+
\frac{1}{2}\tau_u(\phi^{iq})\tau_q(\phi^{ru})\rho^{rp}\bigr\}
$$
Let us denote the $n$-th summand in an expression $X$ by $Xn$ (where we open 
the brackets).  
We have 
$$
B1+C1=-\phi^{ip}\tau_p\tau_q\tau_l(\phi^{jq})
+\phi^{jp}\tau_p\tau_q\tau_l(\phi^{iq})\buildrel{(5.4.9)}\over{=}
$$
$$
=
-\tau_l(\phi^{jq})\tau_p\tau_q(\phi^{ip})+\tau_l(\phi^{iq})
\tau_p\tau_q(\phi^{jp})=-C3-B3
$$
Next, 
$$
B2=-\frac{1}{2}\phi^{ip}\bigl\{
\tau_p\tau_u(\phi^{jq})\tau_q(\phi^{ru})\rho^{rl}+
\tau_u(\phi^{jq})\tau_p\tau_q(\phi^{ru})\rho^{rl}+
\tau_u(\phi^{jq})\tau_q(\phi^{ru})\tau_p(\rho^{rl})\bigr\}
$$
We have
$$
B21=-\frac{1}{2}\phi^{ip}\tau_p\tau_u(\phi^{jq})\tau_q(\phi^{ru})\rho^{rl}
\buildrel{(5.4.7)}\over{=}
$$
$$
=-\frac{1}{2}\bigl\{\phi^{jp}\tau_p\tau_u(\phi^{iq})+
\tau_u(\phi^{jp})\tau_p(\phi^{iq})-
\tau_u(\phi^{ip})\tau_p(\phi^{jq})\bigr\}
\tau_q(\phi^{ru})\rho^{rl}
$$
Next, 
$$
B22=-\frac{1}{2}\phi^{ip}\tau_p\tau_q(\phi^{ru})\tau_u(\phi^{jq})\rho^{rl}
\buildrel{(5.4.7)}\over{=}
$$
$$
=-\frac{1}{2}\tau_u(\phi^{jq})\rho^{rl}\bigl\{
\phi^{rp}\tau_p\tau_q(\phi^{iu})+\tau_q(\phi^{rp})\tau_p(\phi^{iu})-
\tau_q(\phi^{ip})\tau_p(\phi^{ru})\bigr\}=
$$
$$
=-\frac{1}{2}\tau_u(\phi^{jq})\tau_l\tau_q(\phi^{iu})
-\frac{1}{2}\tau_u(\phi^{jq})\tau_q(\phi^{rp})\tau_p(\phi^{iu})\rho^{rl}
+\frac{1}{2}\tau_u(\phi^{jq})\tau_q(\phi^{ip})\tau_p(\phi^{ru})\rho^{rl}
$$
We see that $B221=-A1,\ B222=-B213$ and $B223=-B212$. Similarly, 
$A2=-C221$. We compute $B23$ using (5.4.10): 
$$
B23=\frac{1}{2}\phi^{ip}\tau_u(\phi^{jq})\tau_q(\phi^{ru})\rho^{ra}
\tau_p(\phi^{ba})\rho^{bl}\buildrel{(5.4.6)}\over{=}
$$
$$
=\frac{1}{2}\phi^{bp}\tau_p(\phi^{ia})\tau_u(\phi^{jq})\tau_q(\phi^{ru})
\rho^{ra}\rho^{bl}=
\frac{1}{2}\tau_l(\phi^{ia})\tau_u(\phi^{jq})\tau_q(\phi^{ru})\rho^{ra}=
-B4
$$
Finally, $B211+C211=\beta^{ijr}\rho^{rl}\omega_l$. Everything except these 
last terms cancels out, and this proves the Lemma. $\btu$ 

\bigskip\bigskip


\centerline{\bf \S 6. Atiyah}

\bigskip\bigskip

{\bf 6.1.} Let $\CT=(A,T,\Omega,\dpar)$ be a perfect extended Lie 
algebroid. Let $\fg, \fg', \fg''$ be three frames in $T$. According 
to 5.2, we have the morphisms of the corresponding vertex prealgebroids 
$$
\CB_{\CT;\fg''}\buildrel{g_{\fg',\fg''}}\over{\lra}
\CB_{\CT;\fg'}\buildrel{g_{\fg,\fg'}}\over{\lra}\CB_{\CT;\fg}
\eqno{(6.1.1)}
$$
as well as the morphism $g_{\fg;\fg''}:\ \CB_{\CT;\fg''}\lra
\CB_{\CT;\fg}$, all of them over $Id_\CT$. 
Recall that $Hom_{\CP re \CA lg_\CT}(\CB_{\CT;\fg''},\CB_{\CT;\fg'})$ is an  
$\Omega^{2}(\CT)$-torseur, cf. 4.3.1. We are aiming to compute the 
discrepancy
$$
\alpha_{\fg,\fg',\fg''}:=g_{\fg;\fg'}\circ g_{\fg',\fg''}-g_{\fg,\fg''}
\in\Omega^{2}(\CT)
\eqno{(6.1.2)}
$$
We have the functions $h_{\fg,\fg'}$, etc., acting from $T$ to $\Omega$ 
(not $A$-linear!),  
as in the previous Section, which define our morphisms.  
By (3.7.1) 
the composition $g_{\fg;\fg'}\circ g_{\fg',\fg''}$ is defined 
by the function $h_{\fg;\fg'}+h_{\fg';\fg''}$; therefore the 
discrepancy (6.1.2) is defined by the $A$-linear function 
$$
\alpha_{\fg,\fg',\fg''}=h_{\fg,\fg'}+h_{\fg',\fg''}-h_{\fg,\fg''}\in
\Omega^{2}(\CT)\subset Hom_A(T,\Omega)
\eqno{(6.1.3)}
$$
which by definition coincides with (6.1.2). 

Note that the functions (6.1.3) obviously satisfy the $2$-cocycle condition 
$$
\alpha_{\fg',\fg'',\fg'''}-\alpha_{\fg,\fg'',\fg'''}+
\alpha_{\fg,\fg',\fg'''}-\alpha_{\fg,\fg',\fg''}=0
\eqno{(6.1.4)}
$$

{\bf 6.2.} Choose some bases $\{\tau_i\}, \{\tau'_i\}$ and $\{\tau''_i\}$ 
of our frames $\fg, \fg', \fg''$; let $\{\omega_i\}, \{\omega'_i\}, 
\{\omega''_i\}$ be the dual bases in $\Omega$.   
Define two matrices $\phi, \psi\in GL_n(A)$ 
by $\tau'_i=\phi^{ij}\tau_j,\ \tau''_i=\psi^{ij}\tau'_j$ and set 
$\rho:=\phi^{-1t},\ \sigma:=\psi^{-1t}$. 

The maps $h_{\fg,\fg'}$, etc., are defined by the matrices 
$(h^{ij}_{\fg,\fg'})$, etc.,  
where $h_{\fg,\fg'}(\tau'_i)=h^{ij}_{\fg,\fg'}\omega_j$, etc. 

By $(3.5\gamma)$ and $(5.7.1)_\gamma$ we have 
$$
h_{\fg,\fg'}(a\tau'_i)=ah_{\fg,\fg'}(\tau'_i)-\gamma(a,\tau'_i)=
ah_{\fg,\fg'}(\tau'_i)-\gamma(a,\phi^{ip}\tau_p)=
$$
$$
=ah_{\fg,\fg'}(\tau'_i)+\tau_p(a)\dpar\phi^{ip}+\tau_p(\phi^{ip})\dpar a
$$
Therefore, 
$$
A:=h_{\fg,\fg'}(\tau''_i)=h_{\fg,\fg'}(\psi^{ij}\tau'_j)=
\psi^{ij}h_{\fg,\fg'}(\tau'_j)+\tau_p(\psi^{ij})\dpar\phi^{jp}+
\tau_p(\phi^{jp})\dpar\psi^{ij};
$$ 
$$
B:=h_{\fg',\fg''}(\tau''_i)=h^{ia}_{\fg',\fg''}\omega'_a=
h_{\fg',\fg''}^{ia}\rho^{al}\omega_l=
$$
by (5.7.2)
$$
=
\bigl\{\tau'_p\tau'_a(\psi^{ip})+\frac{1}{2}\tau'_q(\psi^{ip})
\tau'_p(\psi^{rq})\sigma^{ra}\bigr\}\rho^{al}\omega_l
$$
and 
$$
C=-h_{\fg,\fg''}(\tau''_i)=-h_{\fg,\fg''}^{il}\omega_l=
-\bigl\{\tau_p\tau_l((\psi\phi)^{ip})+\frac{1}{2}\tau_q((\psi\phi)^{ip})
\tau_p((\psi\phi)^{rq})(\sigma\rho)^{rl}\bigr\}\omega_l
$$
We have to calculate $\alpha_{\fg,\fg',\fg''}(\tau''_i)=A+B+C$. 
In the computation 
we shall use the same convention for the notation of various terms 
in $A, B, C, \ldots$ as in 5.7. 

So, we have 
$$
A1=\psi^{ij}h^{jl}_{\fg,\fg'}\omega_l=\psi^{ij}\bigl\{\tau_p\tau_l
(\phi^{jp})+\frac{1}{2}\tau_q(\phi^{jp})\tau_p(\phi^{rq})\rho^{rl}\bigr\}
\omega_l
$$
By (5.4.4), 
$$
A2=\tau_p(\psi^{ij})\tau_l(\phi^{jp})\omega_l
$$
and 
$$
A3=\tau_p(\phi^{jp})\tau_l(\psi^{ij})\omega_l
$$
Next, 
$$
B1=\phi^{pu}\tau_u\phi^{av}\tau_v(\psi^{ip})\rho^{al}\omega_l=
\phi^{pu}\tau_u\tau_l(\psi^{ip})\omega_l+
\phi^{pu}\tau_u(\phi^{av})\tau_v(\psi^{ip})\rho^{al}\omega_l
$$
and 
$$
B2=\frac{1}{2}\phi^{qu}\phi^{pv}\tau_u(\psi^{ip})\tau_v(\psi^{rq})
\sigma^{ra}\rho^{al}\omega_l
$$
Finally, 
$$
C1=-\tau_p\tau_l(\psi^{iu}\phi^{up})\omega_l=-\phi^{up}\tau_p\tau_l
(\psi^{iu})\omega_l-\psi^{iu}\tau_p\tau_l(\phi^{up})\omega_l-
$$
$$
-\tau_p(\psi^{iu})\tau_l(\phi^{up})\omega_l-\tau_l(\psi^{iu})
\tau_p(\phi^{up})\omega_l
$$
and
$$  
C2=-\frac{1}{2}\bigl\{\phi^{up}\tau_q(\psi^{iu})+\psi^{iu}\tau_q(\phi^{up})
\bigr\}\bigl\{\phi^{vq}\tau_p(\psi^{rv})+\psi^{rv}\tau_p(\phi^{vq})\bigr\}
\sigma^{rs}\rho^{sl}\omega_l=
$$
$$
=-\frac{1}{2}\bigl\{\phi^{up}\phi^{vq}\tau_q(\psi^{iu})\tau_p(\psi^{rv})
\sigma^{rs}\rho^{sl}+\psi^{iu}\phi^{vq}\tau_p(\psi^{rv})\tau_q(\phi^{up})
\sigma^{rs}\rho^{sl}+
$$
$$
+\phi^{up}\tau_q(\psi^{iu})\tau_p(\phi^{sq})\rho^{sl}+
\psi^{iu}\tau_q(\phi^{up})\tau_p(\phi^{sq})\rho^{sl}\bigr\}\omega_l
$$
We see first of all the terms of the second order cancel out, as they should: 
$A11=-C12$ and $B11=-C11$. Most of the other terms also cancel out, and 
in $A+B+C$ we are left only with 
$$
B12+C23=\frac{1}{2}\phi^{up}\tau_p(\phi^{sq})\tau_q(\psi^{iu})\rho^{sl}\omega_l
\buildrel{(5.4.6)}\over{=}
$$
$$
=\frac{1}{2}\phi^{sp}\tau_p(\phi^{uq})\tau_q(\psi^{iu})\rho^{sl}\omega_l=
\frac{1}{2}\tau_l(\phi^{uq})\tau_q(\psi^{iu})\omega_l
$$
and
$$
C22=-\frac{1}{2}\psi^{iu}\phi^{vq}\tau_q(\phi^{up})\tau_p(\psi^{rv})
\sigma^{rs}\rho^{sl}\omega_l\buildrel{(5.4.6)}\over{=}
-\frac{1}{2}\psi^{iu}\phi^{uq}\tau_q(\phi^{vp})\tau_p(\psi^{rv})
\sigma^{rs}\rho^{sl}\omega_l
$$
Thus, we have
$$
\alpha_{\fg,\fg',\fg''}(\tau''_i)=
\frac{1}{2}\tau_l(\phi^{uq})\tau_q(\psi^{iu})\omega_l
-\frac{1}{2}\psi^{iu}\phi^{uq}\tau_q(\phi^{vp})\tau_p(\psi^{rv})
\sigma^{rs}\rho^{sl}\omega_l
$$
Rewriting the right hand side in the base $\{\omega''_i\}$ we get 
$$
\alpha(\tau''_i):=\alpha_{\fg,\fg',\fg''}(\tau''_i)=\alpha^{ir}\omega''_r
=\alpha_1^{ir}\omega_r''+\alpha_2^{ir}\omega''_r
\eqno{(6.2.1)}
$$
where
$$
\alpha_1^{ir}=\frac{1}{2}\psi^{ra}\phi^{al}\tau_l(\phi^{uq})\tau_q(\psi^{iu})
\eqno{(6.2.2)}
$$
and
$$
\alpha_2^{ir}=-\frac{1}{2}\psi^{iu}\phi^{uq}\tau_q(\phi^{vp})\tau_p(\psi^{rv})
\eqno{(6.2.3)}
$$
So we see that $\alpha_1^{ir}=-\alpha_2^{ri}$, i.e. the matrix $\alpha^{ir}$ 
is skew symmetric, that is, $\alpha\in\Omega^2(\CT)$ as it should be. 

{\bf 6.3.} Let us rewrite the expression (6.2.3) in terms of vector fields 
$\tau''_i$: 
$$
\alpha_2^{ir}=-\frac{1}{2}\tau''_i(\phi^{vp})\phi^{-1pa}\psi^{-1ab}
\tau''_b(\psi^{rv})=
$$
(using the identity $\tau(\phi)=-\phi\tau(\phi^{-1})\phi$)
$$
=-\frac{1}{2}\phi^{vs}\tau''_i(\phi^{-1sa})\psi^{-1ab}\psi^{ru}
\tau''_b(\psi^{-1uc})\psi^{cv}
$$
Using (5.4.13), we have 
$$
\psi^{-1ab}\tau''_b(\psi^{-1uc})=\tau'_a(\psi^{-1uc})=\tau'_u(\psi^{-1ac})
=\psi^{-1ub}\tau''_b(\psi^{-1ac})
$$
whence
$$
\alpha_2(\tau''_i):=\alpha_2^{ir}\omega''_r=-\frac{1}{2}\phi^{vs}
\tau''_i(\phi^{-1sa})\psi^{ru}\psi^{-1ub}\tau''_b(\psi^{-1ac})
\psi^{cv}\omega''_r=
$$
$$
=-\frac{1}{2}\phi^{vs}
\tau''_i(\phi^{-1sa})\tau''_r(\psi^{-1ac})\psi^{cv}\omega''_r=
-\frac{1}{2}tr\bigl\{\phi\tau''_i(\phi^{-1})\tau''_r(\psi^{-1})\psi\bigr\}
\omega''_r
$$ 
Skew symmetrizing, we arrive at the first 
part of 

{\bf 6.4. Theorem.} (a) {\it The cocycle $\alpha_{\fg,\fg',\fg''}$}, (6.1.3) 
{\it is given in coordinates by the expression 
$$
\alpha_{\fg,\fg',\fg''}(\tau''_i)=\alpha(\psi,\phi)(\tau''_i)=
$$
$$
=\frac{1}{2}tr\bigl\{\tau''_i(\psi^{-1})\psi \phi\tau''_r(\phi^{-1})-
\tau''_r(\psi^{-1})\psi \phi\tau''_i(\phi^{-1})\bigr\}\omega''_r
\eqno{(6.4.1)}
$$} 

(b) {\it The $3$-form $\beta_{\fg,\fg'}=\beta(\phi)$} (5.5.3) 
{\it is equal to
$$
\beta_{\fg,\fg'}(\tau'_i,\tau'_j)=\beta(\phi)(\tau'_i,\tau'_j)=
$$
$$
=\frac{1}{2}tr\bigl\{\tau'_i(\phi^{-1})\phi\tau'_j(\phi^{-1})\phi
\tau'_r(\phi^{-1})\phi-
\tau'_j(\phi^{-1})\phi\tau'_i(\phi^{-1})\phi
\tau'_r(\phi^{-1})\phi\bigr\}\omega'_r
\eqno{(6.4.2)}
$$}

Part (b) is proven by the same argument as in 6.3, and we leave it 
to the reader. 

From the expression (6.4.2) we see directly that the form $\beta(\phi)$ 
is closed, 
and from (6.4.1) one checks that the form  $\alpha(\phi,\psi)$ 
satisfies the group (or Cech) cocycle condition
$$
\alpha(\psi,\chi)-\alpha(\phi\psi,\chi)+\alpha(\phi,\psi\chi)-
\alpha(\phi,\psi)=0
\eqno{(6.4.3)}
$$
for all $\phi, \psi, \chi\in GL_n(A)$, which is the same as (6.1.4), 
and 
$$
d\alpha(\psi,\phi)=\beta(\phi)+\beta(\psi)-\beta(\psi\phi)
\eqno{(6.4.4)}
$$ 
Thus, a couple $c(\CT)=(\alpha,\beta)$ 
represents a cohomology class in $H^2(GL_n(A),\Omega^{[2,3\rangle}(\CT))$ 
of the group $GL_n(A)$ with coefficients in the complex
$$
\Omega^{[2,3\rangle}(\CT):=\Omega^2(\CT)\lra\Omega^{3,cl}(\CT)
\eqno{(6.4.5)}
$$
where $\Omega^2(\CT)$ sits in degree $0$  
(the action of $GL_n(A)$ on $\Omega^{[2,3\rangle}$ being trivial). 

The cocycle $\alpha$ is classical, and essentially goes back to 
Atiyah; it is written down explicitly by Harris, [H], p. 280. The class 
$\beta$ resembles "Chern-Simons" form.   
 
The whole cocycle $c(\CT)$ may be thought of as an integration of 
a cocycle $\tc$ from [MSV], (5.16), (5.17). 

\bigskip\bigskip


\centerline{\bf \S 7. Gerbes of Vertex Algebroids}

\bigskip\bigskip 

{\bf 7.1.} Let us reformulate the results of the last three Sections 
in language of Torseurs. This reformulation was inspired by [BD1].   

Let $\CT=(A,T,\ldots)$ be a {\it quasiregular} extended Lie algebroid 
(see 5.1). 
Let us define a groupoid $\Omega^{[2,3\rangle}_\CT$ is follows. 
We set $Ob\ \Omega^{[2,3\rangle}_\CT=\Omega^{3,cl}(\CT)$; for 
$\omega_1,\omega_2\in\Omega^{3,cl}(\CT)$ a morphism $\omega_1\lra\omega_2$ is 
by definition a two-form $\eta\in\Omega^2(\CT)$ such that $d_{DR}(\eta)=
\omega_2-\omega_1$. The composition of morphisms is defined in an obvious 
manner. Note that $\Omega^{[2,3\rangle}_\CT$ is in fact an {\it abelian 
group} in categories. 

Similarly, let $\Omega^{2,cl}_\CT$ denote a groupoid with the unique 
object and the set of morphisms equal to $\Omega^{2,cl}(\CT)$. It is 
also an abelian group in categories. We have a fully faithful monoidal functor 
$$
\Omega^{2,cl}_\CT\lra\Omega^{[2,3\rangle}_\CT
\eqno{(7.1.1)}
$$
sending the unique object in $\Omega^{2,cl}_\CT$ to $0\in\Omega^{3,cl}(\CT)$. 

Consider the groupoid $\CA lg_\CT$. According to 5.3, it is nonempty.   

We can define an {\it Action} of 
$\Omega^{[2,3\rangle}_\CT$ on $\CA lg_\CT$, i.e. a monoidal functor 
$$
\dplus:\ \CA lg_\CT\times\Omega^{[2,3\rangle}_\CT\lra\CA lg_\CT
\eqno{(7.1.2)}
$$
as follows. For $\CA\in\CA lg,\ \omega\in\Omega^{3,cl}(\CT)$, a couple 
$(\CA,\omega)$ goes to the vertex algebroid $\CA\dplus\omega$ defined 
in 4.2. If $\eta\in\Omega^2(\CT)$, $d_{DR}(\eta)=\omega'-\omega$ then 
a morphism 

$\dplus(\eta):\ \CA\dplus\omega\lra \CA\dplus\omega'$ is defined 
according to 4.2.2. 

Let us fix $\CA$ and consider the functor 
$$
\Omega^{[2,3\rangle}_\CT\lra\CA lg_\CT,\ \ \omega\mapsto\CA\dplus\omega
\eqno{(7.1.3)}
$$   
induced by (7.1.2). By 4.2.2 this functor is fully faithful. 

Let $\CA'\in\CA lg_\CT$ be another object. By 5.2, there exists 
a morphism of vertex prealgebroids $g:\ P(\CA)\lra P(\CA')$ lying over 
$Id_\CT$; it is necessarily an isomorphism. Consider the vertex 
algebroid $g_*\CA$ constructed in Theorem 4.5. By definition, 
$g$ is lifted to an isomorphism $\tg:\ \CA\iso g_*\CA$. Since 
$g_*\CA\in\CA lg_{P(\CA')}$, by 4.2.1 $\CA'=g_*\CA\dplus\omega$ for 
some $\omega\in\Omega^{3,cl}(\CT)$. Therefore $\tg$ induces an isomorphism 
$\tg\dplus\omega:\ \CA\dplus\omega\iso g_*\CA\dplus\omega=\CA'$.
 
In other words, we have checked that (7.1.3) is surjective on isomorphism 
classes of objects, hence it is an equivalence of categories. This proves 

{\bf 7.2. Theorem.} {\it Let $\CT$ be a quasiregular extended Lie algebroid. 
Then the Action} (7.1.2) {\it makes the groupoid $\CA lg_\CT$ a nonempty 
$\Omega^{[2,3\rangle}_\CT$-Torseur.} $\btu$

{\bf 7.3.} Let $(X,\CO_X)$ be a topological space ringed by a sheaf of 
commutative $k$-algebras $\CO_X$. Let us call an extended Lie $\CO_X$ 
Lie algebroid $\CT=(\CO_X,T,\ldots)$ (quasi)regular if there exists an open 
covering $X=\bigcup U_i$ of $X$ such that all $\CT(U_i)$ are 
(quasi)regular. 

For example, if $T$ is a Lie $\CO_X$-algebroid which is locally free 
as an $\CO_X$-module then the corresponding extended algebroid $\CT_T$ 
is quasiregular. If $X$ is a smooth $k$-scheme of finite type and 
$T=T_{X/k}$ is the sheaf of vector fields then $\CT_T$ is regular. 

We can sheafify the constructions of the previous Subsections and obtain the 
sheaves ({\it champs}) of groupoids (i.e. {\it gerbes}) $\CA lg_\CT$, etc. 

Let $\CT$ be quasiregular. Consider the gerbe $\CA lg_\CT$. According 
to 7.1 and 7.2, it is locally nonempty but {\it not} locally connected 
in general; its sheaf of 
connected components is an $H^3_{DR}(\CT)$-torseur. 

By the general procedure this gerbe defines a {\it characteristic class} 
$$
c(\CT):=c(\CA lg_\CT)\in H^2(X;\Omega^{[2,3\rangle}(\CT))
\eqno{(7.3.1)}
$$
Here in the right hand side we consider the hypercohomology 
with coefficients in the complex $\Omega^{[2,3\rangle}(\CT)$. 

Let us explain how to define the class (7.3.1). Choose an open 
covering $\CU=\{U_i\}$ of $X$ such that all groupoids $\CA lg_{\CT(U_i)}$ are 
nonempty. Choose an object $\CA_i$ in each $\CA lg_{\CT(U_i)}$. Over double 
intersections $U_{ij}:=U_i\cap U_j$ we get two objects $\CA_i|_{U_{ij}},  
\CA_j|_{U_{ij}}\in \CA lg_{\CT(U_{ij})}$. Choose $3$-forms 
$\omega_{ij}\in\Omega^{3,cl}(\CT(U_{ij})$ such that there exist isomorphisms 
$$
h_{ij}:\ \CA_j|_{U_{ij}}\iso \CA_i|_{U_{ij}}\dplus\omega_{ij}
\eqno{(7.3.2)}
$$
Choose some isomorphisms (7.3.2).  
Then on triple intersections we get isomorphisms 
$$
\CA_i|_{U_{ijk}}\dplus\omega_{ij}\dplus\omega_{ik}\iso
\CA_i|_{U_{ijk}}\dplus\omega_{ik}
\eqno{(7.3.3)}
$$
The isomorphisms (7.5.3) must be given by the $2$-forms $\eta_{ijk}\in\Omega^2
(\CT)$ such that $d_{DR}(\eta_{ijk})=c_{ij}-c_{ik}+c_{jk}$. Then 
$(\omega_{ij},\eta_{ijk})$ is a $2$-cocycle in the Cech complex 
$C^\cdot(\CU;\Omega^{[2,3\rangle}(\CT))$ representing the class 
(7.3.1). 

{\bf 7.4.} Now let us assume that the groupoid $\CT$ is regular.  
Theorem 6.4 calculates the class (7.3.1). 
 
Namely, define 
the "Atiyah-Chern-Simons"  
class $ch_2(\CT)\in H^2(X;\Omega^{[2,3\rangle}(\CT))$ by the following 
procedure. Choose some bases of local sections 
$\tau^{(i)}=\{\tau^{(i)}_{\alpha}\}\subset T(U_i)$ over some open covering. 
Let $\phi_{ij}\in GL_n(\CO(U_{ij}))$ be the transition matrix from 
$\tau^{(i)}$ to $\tau^{(j)}$ over $U_{ij}$. By definition, 
$ch_2(\CT)$ is represented by the Cech $2$-cocycle 
$$
2ch_2(\CT):=(\alpha(\CT),\beta(\CT))
\eqno{(7.4.1)}
$$ 
where  
$$
\alpha(\CT)=(\alpha(\CT)_{ijk})=
tr(\phi_{ij}^{-1}\dpar\phi_{ij}\wedge\dpar\phi_{jk}\cdot
\phi_{jk}^{-1})\}
\in Z^2(\CU;\Omega^{2}(\CT))
\eqno{(7.4.2)}
$$
and
$$
\beta(\CT)=(\beta(\CT)_{ij})=\{\frac{1}{3}tr(\dpar\phi_{ij}\phi_{ij}^{-1}
\dpar\phi_{ij}\phi_{ij}^{-1}\dpar\phi_{ij}\phi_{ij}^{-1})\}\in 
C^1(\CU;\Omega^{3,cl}(\CT))
\eqno{(7.4.3)}
$$ 

Theorem 6.4 implies 

{\bf 7.5. Theorem.} {\it Let $\CT$ be a regular extended Lie $\CO_X$-algebroid. 
Then $c(\CT)=2ch_2(\CT)$. 

Therefore, the gerbe $\CA lg_\CT$ 
admits a global section iff $ch_2(\CT)=0$. If so, then the groupoid of 
global sections $\CA lg_\CT(X)$ is equivalent to the groupoid 
of $\Omega^{[2,3\rangle}(\CT)$-torseurs, whence $\pi_0(\CA lg_\CT(X))
\buildrel\sim\over{=} H^1(X;\Omega^{[2,3\rangle}(\CT))$ and the automorphism 
group of an 
object of this groupoid is isomorphic to $H^0(X,\Omega^{[2,3\rangle}(\CT))$.} 

$\btu$ 

It is instructive to compare the previous discussion with [BB], 2.1. 
In the chiral situation the degree of cohomology goes one step up. 

{\bf 7.6.} Let us identify the class $c(\CT)$ when $X$ is a smooth 
$k$-scheme and $\CT=\CT_X$ is the tangent bundle. If $E$ is an arbitrary 
vector bundle over $X$ given by a Cech $1$-cocycle 
$$
g=(g_{ij})\in Z^1(\CU;GL_r(\CO_X))
\eqno{(7.6.1)}
$$
on some open covering $\CU$ then (7.4.1) - (7.4.3) define a cocycle
$$
c(g)=(\alpha(g),\beta(g))\in Z^2(\CU;\Omega^{[2,3\rangle}_X)
\eqno{(7.6.2)}
$$
If $g_{ij}=\phi_ih_{ij}\phi_j^{-1}$ for some $\phi=(\phi_i)\in 
C^1(\CU; GL(\CO_X)$ then one checks by a direct computation that 

{\bf 7.6.1. Claim.} (a)  
$$
\alpha(g)-\alpha(h)=d_{Cech}\eta
\eqno{(7.6.3)}
$$
{\it where $\eta=\eta(h,\phi)=(\eta_{ij})\in C^2(\CU;\Omega^2_X)$ is given 
by 
$$
\eta_{ij}=tr\bigl\{h_{ij}^{-1}dh_{ij}\phi_j^{-1}d\phi_j-
\phi_i^{-1}d\phi_i dh_{ij} h_{ij}^{-1}+
h_{ij}^{-1}\phi_i^{-1}d\phi_i h_{ij}\phi_j^{-1}d\phi_j\bigr\}
\eqno{(7.6.4)} 
$$}
(b) {\it  
$$
d_{DR}\eta=\beta(g)-\beta(h)-d_{Cech}\gamma
\eqno{(7.6.5)}
$$
where $\gamma=\gamma(h,\phi)=(\gamma_i)\in C^0(\CU;\Omega^{3,cl})$ is 
defined by 
$$
\gamma_i=\frac{1}{3}tr((\phi_i^{-1}d\phi_i)^3)
\eqno{(7.6.6)}
$$}
$\btu$ 

{\bf 7.6.2. Corollary.} {\it 
$$
c(g)-c(h)=de
\eqno{(7.6.7)}
$$
where 
$$
e=(\eta,\gamma)\in C^1(\CU;\Omega^{[2,3\rangle}_X)
\eqno{(7.6.8)}
$$} 
$\btu$

Therefore (7.6.2) gives rise to a well defined characteristic class 
$c(E)\in H^2(X,\Omega^{[2,3\rangle}_X)$. The following Lemma is obvious.  

{\bf 7.7. Lemma.} (a) {\it If  $f:\ Y\lra X$ is an arbitrary morphism 
from another smooth scheme then $c(f^*E)=f^*c(E)$.} 

(b) {\it If $0\lra E'\lra E\lra E''\lra 0$ is a short exact sequence 
of vector bundles over $X$ then $c(E)=c(E')+c(E'')$.} 

(c) {\it If $L$ is a line bundle then $c(L)$ is equal to the image of 
$c_1(L)\otimes c_1(L)\in H^1(X;\CO^*_X)^{\otimes 2}$ under the composition 
$$
H^1(X;\CO^*_X)^{\otimes 2}\buildrel{d log^{\otimes 2}}\over{\lra}
H^1(X;\Omega^{1,cl}_X)^{\otimes 2}\lra H^2(X;\Omega^{2,cl}_X)
\lra H^2(X;\Omega^{[2,3\rangle}_X)
$$}

{\bf 7.8.} Recall (cf. [S]) that for an arbitrary $E$ we have a 
characteristic class 
$$
2ch_2^{(K)}(E):=c_1^{(K)2}(E)-2c_2^{(K)}\in H^2(X;\CK_{2,X})
$$
from which we can get a class $2ch_2(E)\in H^2(X;\Omega^{[2,3\rangle}_X)$ using 
the $dlog$ map 
$$
H^2(X;\CK_2)\lra H^2(X;\Omega^{2,cl}_X)\lra H^2(X,\Omega^{[2,3\rangle}_X)
$$
The class $2ch_2(E)$ also satisfies properties 7.7 (a) - (c). 

It is obvious that the natural map 
$$
H^2(X;\Omega^{[2,3\rangle}_X)\lra H^2(X;\Omega^{[2}_X)
$$
is injective, where $\Omega^{[2}_X:=\Omega^2_X\lra\Omega^3_X\lra\ldots$ 
is the stupidly truncated (and shifted, so that $\Omega^2_X$ 
sits in degree $0$) de Rham complex. 

The proof of the following lemma was provided to us by H.~Esnault.   

{\bf 7.9. Lemma.} {\it The inverse image map
$$
H^2(X;\Omega^{[2}_X)\lra H^2(\BP(E);\Omega^{[2}_{\BP(E)})
$$
is injective.} $\btu$ 

In fact, more is true. One can consider the cohomology theory which 
assigns to a smooth $X$ a collection of cohomology groups 
$\{ H^i(X;\Omega^{[i})\}$. This theory has the standard Grothendieck's 
properties needed 
to define the Chern classes, cf. [Gr]. This fact was communicated to us 
by A.~Beilinson.    

Anyway, the inverse image map $H^2(X;\Omega^{[2,3\rangle}_X)\lra 
H^2(\BP(E);\Omega^{[2,3\rangle}_{\BP(E)})$ is also injective, hence  
by splitting principle we get 

{\bf 7.10. Theorem.} {\it For all vector bundles $E$,  $c(E)=2ch_2(E)$.} $\btu$ 

This theorem was obtained in collaboration with H.~Esnault.  

{\bf 7.11. Corollary.} {\it In the situation of} 7.6  
{\it the class described in} Theorem 7.5 {\it 
is equal to $2ch_2(\CT_X)$.} $\btu$

\bigskip\bigskip


\centerline{\bf \S 8. Vertex Envelope of a Conformal Algebra}

\bigskip\bigskip

{\bf 8.1.} Our aim in this Section will be to construct a left adjoint $U$  
to the forgetful functor (0.8.1) and to prove the 
"Poincar\'e-Birkhoff-Witt" theorem for algebras $UC,\ C\in\CC onf$. 

Let $V$ be a vertex algebra. We have the following two particular 
cases of the OPE formula (0.5.12). The first one corresponds to 
$m=n=-1$: 
$$
[x_{(-1)},y_{(-1)}]=\sum_{j\geq 0}\ (-1)^j\dpar^{(j+1)}(x_{(j)}y)_{(-1)}
\eqno{(8.1.1)}
$$
where we have used (0.5.10). The second one corresponds to 
$m\geq 0,\ n=-1$: 
$$
[x_{(m)},y_{(-1)}]=\sum_{j=0}^m\ \binom{m}{j}(x_{(j)}y)_{(m-j-1)}
\eqno{(8.1.2)}
$$ 

{\bf 8.2.} Let $C=\oplus\ C_i$ be a conformal algebra. Let 
$TC=\sum_{j\in\BZ_{\geq 0}}\ T^jC,\ T^jC:=C^{\otimes j}$ be the 
tensor algebra of $C$ over $k$. The multiplication in $TC$ will be denoted 
$x\cdot y$ or $xy$. The $\BZ_{\geq 0}$-grading of $C$ induces the 
$\BZ_{\geq 0}$-grading of $TC$ such that $TC$ becomes a $\BZ_{\geq 0}$-graded 
associative algebra with the unit $\b1\in T^0C=k$. We have a 
canonical embedding of $k$-modules $C=T^1C\subset TC$. 

There is a unique extension of the operators $\dpar^{(j)}$ on $C$ to the 
whole space $TC$ satisfying 
$$
\dpar^{(j)}(xy)=\sum_{p=0}^j\ \dpar^{(p)}x\cdot\dpar^{(j-p)}y
\eqno{(8.2.0)}
$$
cf. (0.4.6). These operators will satisfy (0.4.1). 

Let $R\subset TC$ be a two-sided ideal generated by all elements 
$$
r(x,y):=xy-yx-\sum_{j\geq 0}\ (-1)^j\dpar^{(j+1)}(x_{(j)}y),\ 
x,y\in C, 
\eqno{(8.2.1)}
$$
cf. (8.1.1). Set $UC=TC/R$. We have a canonical morphism $i:\ C\lra UC$ equal 
to the restriction of the projection $p:\ TC\lra UC$ to $C$. 

Since the relations (8.2.1) are homogeneous, the algebra $UC$ inherits 
a $\BZ_{\geq 0}$-grading from $TC$. Using (0.4.6) and an obvious identity 
$\dpar^{(i)}\dpar^{(j)}=\dpar^{(j)}\dpar^{(i)}$ one sees easily 
that the operators $\dpar^{(i)}$ respect the ideal $R$. Hence they induce 
the operators $\dpar^{(i)}$ on $UC$ of degree $i$ which satisfy (0.4.1).   

{\bf 8.3. Theorem.} {\it There is a unique structure of a vertex algebra 
on the $k$-module $UC$ such that for all $x\in C,\ z\in TC$, 
$$
p(xz)=i(x)_{(-1)}p(z)
\eqno{(8.3.1)}
$$
The corresponcence $C\mapsto UC$ defines a functor $U:\ \CC onf\lra\CV ert$ 
left adjoint to the forgetful functor.} 

The algebra $UC$ will be called a {\it vertex envelope} of a 
conformal algebra $C$. 

This Theorem will be proven in 8.20, after some preparation. 

{\bf 8.4.} Let us define $k$-linear operators $x_{(j)},\ x\in C,\ 
j\in \BZ$ of degree $-j-1$ acting on the module $TC$. If $j<0,\ j=-n-1$, we set 
$$
x_{(-n-1)}z=(\dpar^{(n)}x)z
\eqno{(8.4.1)}
$$
In particular, $x_{(-1)}z=xz$. 

Each element of $TC$ is a linear combination of monomials $z=z_1z_2\ldots z_n,\ 
z_i\in C$. We define $x_{(j)}z,\ j\in\BZ_{\geq 0}$ by induction on the length 
$n$ of the monomial. If $n=1$, i.e. $z\in C$, then we already have   
$x_{(j)}z$ due to the structure of a conformal algebra on $C$. 

If $z=yu,\ y\in C$, we set 
$$
x_{(j)}yu=yx_{(j)}u+\sum_{p=0}^j\ \binom{j}{p}(x_{(p)}y)_{(j-p-1)}u
\eqno{(8.4.2)}
$$
cf. (8.1.2). 

We leave to the reader an easy proof of the lemma below.  

{\bf 8.4.1. Lemma.} {\it For all $i\in\BZ_{\geq 0},\ n\in\BZ$,  
$$
(\dpar^{(i)}x)_{(n)}=(-1)^i\binom{n}{i}x_{(n-i)}
\eqno{(8.4.3)}
$$}

{\bf 8.5. Lemma.} {\it The operators $x_{(j)}$ introduced above respect 
the ideal $R$.} 

The proof will be given in 8.6 --- 8.8 below.  

{\bf 8.6.} Since $R$ is a left ideal, it is respected 
by all operators $x_{(n)}$ with $n<0$.  

It follows from the commutation formula (8.4.2) that 
if all the operators $x_{(n)},\ n\geq 0,$ respect a subset $S\subset TC$ 
then they 
respect the left ideal generated by $S$. Therefore we need to prove that 
$$
u_{(n)}xyz-u_{(n)}yxz-\sum_{j\geq 0}\ (-1)^ju_{(n)}\dpar^{(j+1)}(x_{(j)}y)z\in 
R
\eqno{(8.6.1)}
$$
for all $u, x, y\in C,\ z\in TC,\ n\geq 0$. Let us denote the summands 
in (8.6.1) by $A, B$ and $C$. We have 
$$
A=u_{(n)}xyz=(u_{(n)}x)yz+xu_{(n)}yz+\sum_{p=o}^{n-1}\ \binom{n}{p}
(u_{(p)}x)_{(n-p-1)}yz
\eqno{(8.6.2)}
$$
We shall use the same agreement as in 5.7: in an expression $X$, 
the $n$-th summand will be denoted by $Xn$. We have 
$$
A2=x(u_{(n)}y)z+xyu_{(n)}z+\sum_{q=0}^{n-1}\ \binom{n}{q}\ 
x(u_{(q)}y)_{(n-q-1)}z
\eqno{(8.6.3)}
$$
$$
A3=\sum_{p=0}^{n-1}\ \binom{n}{p}\bigl\{
\bigl((u_{(p)}x)_{(n-p-1)}y\bigr)z+y(u_{(p)}x)_{(n-p-1)}z+
$$
$$
+\sum_{q=0}^{n-p-2}\ \binom{n-p-1}{q}\bigr((u_{(p)}x)_{(q)}y\bigr)_{(n-p-q-2)}z
\bigr\}
\eqno{(8.6.4)}
$$
Similarly, 
$$
B=-u_{(n)}yxz=-(u_{(n)}y)xz-yu_{(n)}xz-
\sum_{p=0}^{n-1}\ \binom{n}{p} (u_{(p)}y)_{(n-p-1)}xz
\eqno{(8.6.5)}
$$
$$
B2=-y(u_{(n)}x)z-yxu_{(n)}z-\sum_{p=0}^{n-1}\ \binom{n}{p}
y(u_{(p)}x)_{(n-p-1)}z
\eqno{(8.6.6)}
$$
and
$$
B3=-\sum_{p=0}^{n-1}\ \binom{n}{p}\bigl\{
\bigl((u_{(p)}y)_{(n-p-1)}\bigr)z+x(u_{(p)}y)_{(n-p-1)}z+
$$
$$
+\sum_{q=0}^{n-p-2}\ 
\binom{n-p-1}{q}\bigl((u_{(p)})_{(q)}x\bigr)_{(n-p-q-2)}z\bigr\}
\eqno{(8.6.7)}
$$
Next, 
$$
C=-\sum_{j\geq 0}\ (-1)^j\bigl\{\bigl(u_{(n)}\dpar^{(j+1)}(x_{(j)}y)\bigr)z+
\dpar^{(j+1)}(x_{(j)}y)u_{(n)}z+
$$
$$
+\sum_{p=0}^{n-1}\ \binom{n}{p}\bigl(u_{(p)}\dpar^{(j+1)}(x_{(j)}y)
\bigr)_{(n-p-1)}z\bigr\}
\eqno{(8.6.8)}
$$
Due to (0.4.7), 
$$
C1=-\sum_{j\geq 0}\ (-1)^j\bigl\{\dpar^{(j+1)}(u_{(n)}x_{(j)}y)+
\sum_{p=1}^{min(n,j+1)}\ \binom{n}{p}\dpar^{(j+1-p)}
(u_{(n-p)}x_{(j)}y)\bigr\}z
\eqno{(8.6.9)}
$$
Next, 
$$
C11=-\sum_{j\geq 0}\ (-1)^j\dpar^{(j+1)}\bigl\{
(u_{(n)}x)_{(j)}y+x_{(j)}u_{(n)}y+\sum_{p=0}^{n-1}\ \binom{n}{p}
(u_{(p)}x)_{(n+j-p)}y\bigr\}z
\eqno{(8.6.10)}
$$ 
and
$$
C12=-\sum_{j\geq 0}\ (-1)^j\sum_{p=1}^{min(n,j+1)}\ 
\binom{n}{p}\dpar^{(j+1-p)}\bigl\{(u_{(n-p)}x)_{(j)}y+
x_{(j)}u_{(n-p)}y+
$$
$$
+\sum_{q=0}^{n-p-1}\ \binom{n-p}{q}(u_{(q)}x)
_{(j+n-p-q)}y\bigr\}z
\eqno{(8.6.11)}
$$
We have $A21+B21+C2=r(x,y)u_{(n)}z\in R$. Next, $A23=-B32$ and $A32=-B23$. 
Next, 
$$
A21+B1=x(u_{(n)}y)z-(u_{(n)}y)xz\buildrel{R}\over{\sim}
\sum_j\ (-1)^j\dpar^{(j+1)}(x_{(j)}u_{(n)}y)z=-C112
$$
Similarly,
$$
A1+B21=(u_{(n)}x)yz-y(u_{(n)}x)z\buildrel{R}\over{\sim}
\sum_j\ (-1)^j\dpar^{(j+1)}\bigl((u_{(n)}x)_{(j)}y\bigr)z=-C111
$$

{\bf 8.7.} Below we shall often use the identity
$$
\binom{n}{p}\binom{p}{r}=\binom{n}{r}\binom{n-r}{n-p}
\eqno{(8.7.1)}
$$ 
We claim that
$$
C113+C121+C123=-A31
\eqno{(8.7.2)}
$$
Indeed, we have
$$
C113=-\sum_{j\geq 0}\sum_{q=0}^{n-1} (-1)^j\binom{n}{q}\dpar^{(j+1)}
\bigl\{(u_{(q)}x)_{(n+j-q)}\bigr\}z
\eqno{(8.7.3)}
$$
and
$$
C121+C123=-\sum_{j\geq 0}\sum_{p=1}^{min(n,j+1)}\ 
$$
$$
\sum_{q=0}^{n-p-1}\ 
(-1)^j\binom{n}{p}\binom{n-p}{q}\dpar^{(j+1-p)}
\bigl\{(u_{(q)}x)_{(j+n-p-q)}y\bigr\}z
\eqno{(8.7.4)}
$$
Using (8.7.1), we have
$$
\binom{n}{p}\binom{n-p}{q}=\binom{n}{n-p}\binom{n-p}{q}=
\binom{n}{q}\binom{n-q}{p}
$$
Consider the part of (8.7.4) 
at a fixed $q$, $0\leq q\leq n-1$: 
$$
(8.7.4)_q=-\binom{n}{q}\sum_{j\geq 0}\sum_{p=1}^{min(n-q,j+1)}\ 
$$
$$ 
(-1)^j\binom{n-q}{p}\dpar^{(j+1-p)}\bigl\{(u_{(q)}x)_{(j+n-p-q)}y\bigr\}z=
$$
$$
=-\binom{n}{q}\sum_{p=1}^{n-q}\sum_{j\geq p-1}
(-1)^j\binom{n-q}{p}\dpar^{(j+1-p)}\bigl\{(u_{(q)}x)_{(j+n-p-q)}y\bigr\}z
\eqno{(8.7.5)}
$$
It is easy to see that the part of (8.7.5) corresponding to 
$j=p-1$ is equal to $-A31_q$ and the part corresponding to $j\geq p$ 
equals $-C113_q$. This proves (8.7.2). 

Using the commutativity formula (0.4.3), one sees that $B31=-C122$. 

{\bf 8.8.} We claim that 
$$
A33+B33+C3=0
\eqno{(8.8.1)}
$$
We have 
$$
C3=-\sum_{j\geq 0}\ (-1)^j\sum_{p=0}^{n-1}\ \binom{n}{p}
\sum_{q=0}^{min(p,j+1)}\ \binom{p}{q}\dpar^{(j+1-q)}(u_{(p-q)}x_{(j)}y)
_{(n-p-1)}z=
$$
$$
=-\sum_{j\geq 0}\ (-1)^j\sum_{p=0}^{n-1}\ \binom{n}{p}
\sum_{q=0}^{min(p,j+1)}\ \binom{p}{q}
(-1)^{j+1-q}\binom{n-p-1}{j-q+1}(u_{(p-q)}x_{(j)}y)_{(n-j-p+q-2)}z=
$$
(we set $r:=p-q$)
$$
=\sum_{j\geq 0}\sum_{p=0}^{n-1}\binom{n}{p}\sum_{r=p-j-1}^p\ (-1)^{p-r}
\binom{p}{r}\binom{n-p-1}{j-p+r+1}(u_{(r)}x_{(j)}y)_{(n-j-r-2)}z=
$$
$$
=\sum_{j\geq 0}\sum_{p=0}^{n-1}\binom{n}{p}\sum_{r=p-j-1}^p\ (-1)^{p-r}
\binom{p}{r}\binom{n-p-1}{j-p+r+1}
$$
$$
\bigl\{(x_{(j)}u_{(r)}y)_{(n-j-r-2)}z+
\sum_{s=0}^r\ \binom{r}{s}\bigl((u_{(s)}x)_{(r+j-s)}y\bigr)_{(n-j-r-2)}z
\bigr\}
\eqno{(8.8.2)}
$$
So, we have written $C3$ as a sum of two terms. 
Now we claim that 
$$
C32=-A33
\eqno{(8.8.3)}
$$
Set $l:=r+j-s$. We have 
$$
C32=\sum_{l\geq 0, s\geq 0, l+s\leq n-2}\ \binom{n}{s} 
\sum_{p=0}^{l+s+1}\sum_{r=s}^{min(p,l+s)}\ 
$$
$$ 
(-1)^{p-r}\binom{n-s}{n-p}\binom{p-s}{r-s}\binom{n-p-1}{l+s-p+1}
\bigl((u_{(s)}x)_{(l)}y\bigr)_{(n-l-s-2)}z
$$
where we have used that
$$
\binom{n}{p}\binom{p}{r}\binom{r}{s}=
\binom{n}{r}\binom{r}{s}\binom{n-r}{n-p}=
\binom{n}{s}\binom{n-s}{n-r}\binom{n-r}{n-p}=
\binom{n}{s}\binom{n-s}{n-p}\binom{p-s}{r-s}
$$
Thus, at fixed $l, s$, 
$$
C32_{l,s}=\binom{n}{s}\sum_{p=0}^{l+s+1}\ (-1)^p\binom{n-s}{n-p}
\binom{n-p-1}{l+s-p+1}
\sum_{r=s}^{min(p,l+s)}\ (-1)^r\binom{p-s}{r-s}
\bigl((u_{s)}x)_{(l)}y)_{(n-l-s-2)}z
$$
The last sum is non-zero in two cases: 

(a) $p=s$, so that $r=s=p$, 
and we get the coefficient $\binom{n}{s}\binom{n-s-1}{l+1}$, 

and

(b) $p=l+s+1$ where we get $-\binom{n}{s}\binom{n-s}{l+1}$.

Since 
$$
\binom{n-s-1}{l+1}-\binom{n-s}{l+1}=-\binom{n-s-1}{l}
\eqno{(8.8.4)} 
$$
we get 
$$
C32_{l,s}=-\binom{n}{s}\binom{n-s-1}{l}=-A33_{l,s}
$$
which proves (8.8.3). 

Next, we claim that 
$$
C31=-B33
\eqno{(8.8.5)}
$$
Indeed, 
$$
C31=\sum_{j\geq 0}\sum_{p=0}^{n-1}\ (-1)^{p-r}\binom{n}{p}
\sum_{r=max(0,p-j-1)}^p\ 
\binom{n-p-1}{j-p+r+1}\binom{p}{r}(x_{(j)}u_{(r)}y)_{(n-j-r-2)}z=
$$
(using (0.4.3))
$$
=\sum_{j\geq 0}\sum_{p=0}^{n-1}\ (-1)^{p-r}\binom{n}{p}
\sum_{r=max(0,p-j-1)}^p\ 
\binom{n-p-1}{j-p+r+1}\binom{p}{r}
$$
$$
(-1)^{j+1}\sum_{s\geq 0}\ (-1)^s\dpar^{(s)}
\bigl((u_{(r)}y)_{j+s}x\bigr)_{(n-j-r-2)}z=
$$
$$
=\sum_{j\geq 0}\sum_{p=0}^{n-1}\ (-1)^{p-r}\binom{n}{p}
\sum_{r=max(0,p-j-1)}^p\ 
\binom{n-p-1}{j-p+r+1}\binom{p}{r}(-1)^{j+1}
$$
$$
\sum_{s\geq 0}\binom{n-j-r-2}{s}\bigl((u_{(r)}y)_{(j+s)}x\bigr)
_{(n-j-s-r-2)}z
$$
Set $l:=j+s$. We get
$$
C31=\sum_{l\geq 0, r\geq 0, l+r\leq n-2}\ 
\sum_{s\geq 0, p\geq r, p+s\leq r+l+1}\ 
$$
$$ 
(-1)^{p-r+l-s+1}\binom{n}{p}\binom{p}{r}\binom{n-p-1}{l+r-s-p+1}
\binom{n-l-r+s-2}{s}\bigl((u_{(r)}y)_{(l)}x\bigr)_{(n-l-r-2)}z
$$
Using the identities $\binom{n}{p}\binom{p}{r}=\binom{n-r}{n-p}$ and 
$$
\binom{n-p-1}{l+r-s-p+1}\binom{n-l-r+s-2}{s}=
\binom{n-p-1}{n-l+s-r-2}\binom{n-l-r+s-2}{s}=
$$
$$
=\binom{n-p-1}{s}\binom{n-p-s-1}{l+r-s-p-1}=
\binom{n-p-1}{n-p-s-1}\binom{n-p-s-1}{n-l-r-2}=
$$
$$
=\binom{n-p-1}{n-l-r-2}\binom{l+r-p+1}{s}
$$
we get at fixed $l,r$
$$
C31_{l,r}=(-1)^{l+r+1}\binom{n}{r}\sum_{p=r}^{r+l+1}\ 
(-1)^p\binom{n-r}{n-p}\binom{n-p-1}{n-l-r-2}
$$
$$
\sum_{s=0}^{min(r+l-p+1,l)}\ (-1)^s\binom{l+r-p+1}{s}
\bigl((u_{(r)}y)_{(l)}x\bigr)_{(n-l-r-2)}z
$$
The last sum is non-zero in two cases: 

(a) $p=l+r+1$, so that $s=0$, and we get the coefficient 
$\binom{n}{r}\binom{n-r}{n-l-r-1}=\binom{n}{r}\binom{n-r}{l+1}$, 

and 

(b) $p=r$, so that we get the coefficient 
$-\binom{n}{r}\binom{n-r-1}{n-l-r-2}=-\binom{n}{r}\binom{n-r-1}{l+1}$.

It follows that 
$$
C31_{r,l}=\binom{n}{r}\binom{n-r-1}{l}\bigl((u_{(r)}y)_{(l)}x\bigr)
_{(n-l-r-2)}z=-B33_{r,l}
$$
which proves (8.8.5). The identities (8.8.3) and (8.8.5) together imply 
(8.8.1). 

After bookkeeping, the computations of 8.6 --- 8.8 show that the left 
hand side of (8.6.1) belongs to the ideal $R$, which finishes 
the proof of Lemma 8.5. $\btu$

{\bf 8.9.} By Lemma 8.5, the operators 
$x_{(j)}$ induce operators 
$$
x_{(j)}:\ UC\lra UC,\ 
x\in C, j\in\BZ
\eqno{(8.9.1)}
$$

{\bf 8.10. Lemma.} {\it The operators} (8.9.1) {\it satisfy the OPE formula} 
(0.5.12). 

This lemma will be proven in 8.11 --- 8.19 below. 

{\bf 8.11.} We have to prove the identity between the operators 
acting on $UC$, 
$$
[x_{(m)},y_{(n)}]=\sum_{j\geq 0}\ \binom{m}{j} (x_{(j)}y)_{(m+n-j)}
\eqno{(8.11.1)}
$$
for all $m,n\in\BZ$. 

Let us discuss separately three cases. 

{\bf 8.12.} {\it Case A. $m<0$ and $n<0$.} First of all, note that 
for $m=n=-1$ the relation (8.11.1) is nothing but (8.2.1) and holds true by 
definition, since $R$ is a left ideal. 

In general, set $m=-1-a,\ n=-1-b$ for $a,b\geq 0$. We have 
by definition (8.4.1) 
$$
[x_{(-1-a)},y_{(-1-b)}]=[\dpar^{(a)}x,\dpar^{(b)}y]=
$$
$$
=\sum_{j\geq 0}\ (-1)^j\dpar^{(j+1)}((\dpar^{(a)}x)_{(j)}
\dpar^{(b)}y)
\eqno{(8.12.1)}
$$
The right hand side of (8.11.1) is equal to 
$$
\sum_{j\geq 0}\ \binom{-1-a}{j}(x_{(j)}y)_{(-2-a-b-j)}=
\sum_{j\geq 0} (-1)^j\binom{a+j}{j}\dpar^{(a+b+j+1)}(x_{(j)}y)
\eqno{(8.12.2)}
$$
where we have used that
$$
\binom{-1-a}{j}=(-1)^j\binom{a+j}{j}
\eqno{(8.12.3)}
$$
for $a, j\geq 0$. 

On the other hand,  
$$
(\dpar^{(a)}x)_{(j)}\dpar^{(b)}y=(-1)^a \binom{j}{a}x_{(j-a)}\dpar^{(b)}y=
(-1)^a\binom{j}{a}\sum_{p=0}^b\ \binom{j-a}{p}\dpar^{(b-p)}
(x_{(j-a-p)}y),
$$
by (0.4.2) and (0.4.7). Therefore, the rhs of (8.12.1) is equal to 
$$
\sum_{j\geq 0}\ (-1)^{j+a}\binom{j}{a}\sum_{p\geq 0}\ 
\binom{j-a}{p}\dpar^{(j+1)}\dpar^{(b-p)}(x_{(j-a-p)}y)=
$$
$$
=\sum_{j\geq 0}\ (-1)^{j+a}\binom{j}{a}\sum_{p=0}^{min(b,j-a)}\ 
\binom{j-a}{p}\binom{j+1+b-p}{j+1}\dpar^{(j+1+b-p)}
(x_{(j-a-p)}y)=
$$
(substituting $k:=j-a-p$)
$$
=\sum_{k\geq 0}\sum_{p=0}^b\ (-1)^{p+k}\binom{k+a+p}{a}
\binom{k+p}{p}\binom{k+a+b+1}{k+a+p+1}\dpar^{(k+a+b+1)}(x_{(k)}y)=
$$
(using that $\binom{k+a+p}{a}\binom{k+p}{p}=\binom{k+a+p}{k+p}
\binom{k+p}{p}=\binom{k+a+p}{p}\binom{k+a}{a}$)
$$
=\sum_{k\geq 0}\ (-1)^k\binom{k+a}{a}\biggl\{\sum_{p=0}^b\ 
(-1)^p\binom{k+a+p}{p}
\binom{k+a+b+1}{k+a+p+1}\biggr\}\dpar^{(k+a+b+1)}(x_{(k)}y)
\eqno{(8.12.4)}
$$

{\bf 8.13. Lemma.} {\it For all $b, q\in\BZ_{\geq 0}$, 
$$
\sum_{p=0}^b\ (-1)^p\binom{p+q}{p}\binom{b+q+1}{p+q+1}=1
$$}

This is easily proved by induction on $b$. $\btu$ 

It follows from this lemma that (8.12.4) is equal to (8.12.2), 
which completes the check of {\it Case A}. 

{\bf 8.14.} {\it Case B. $m\geq 0$ and $n<0$.} If $n=-1$ then (8.11.1) 
is the same as the definition (8.4.2). Now let $n=-1-a,\ a\geq 0$. 
The lhs of (8.11.1) is equal to
$$
[x_{(m)},y_{(-1-a)}]=[x_{(m)},\dpar^{(a)}y]=
\sum_{j=0}^{m-1}\ \binom{m}{j}(x_{(j)}\dpar^{(a)}y)_{(m-j-1)}+
x_{(m)}\dpar^{(a)}y
\eqno{(8.14.1)}
$$
The rhs of (8.11.1) is equal to 
$$
\sum_{k=0}^m\ \binom{m}{k}(x_{(k)}y)_{(m-1-a-k)}=
$$
$$
=\sum_{k=0}^{m-1-a}\ \binom{m}{k}(x_{(k)}y)_{(m-1-a-k)}+
\sum_{k\geq max(m-a,0)}\ \binom{m}{k}(x_{(k)}y)_{(m-1-a-k)}
\eqno{(8.14.2)}
$$
Consider the first sum in the rhs of (8.14.1). It is equal to 
$$
\sum_{j=0}^{m-1}\binom{m}{j}\sum_{p=o}^{min(j,a)}\ \binom{j}{p}
\dpar^{(a-p)}(x_{(j-p)}y)_{(m-j-1)}=
$$
$$
=\sum_{j=0}^{m-1}\binom{m}{j}\sum_{p=max(0,a-m+j+1)}^{min(j,a)}\ \binom{j}{p}
(-1)^{a-p}\binom{m-j-1}{a-p}(x_{(j-p)}y)_{(m-j-a+p-1)}=
$$
(making the substitution $k=j-p$)
$$
=\sum_{k=0}^{m-a-1}\sum_{p=0}^a (-1)^{a-p}\binom{m}{k+p}\binom{k+p}{p}
\binom{m-p-k-1}{a-p}(x_{(k)}y)_{(m-k-a-1)}=
$$
(using that $\binom{m}{k+p}\binom{k+p}{p}=\binom{m}{k}\binom{m-k}{m-p-k}$) 
$$
=\sum_{k=0}^{m-a-1}\binom{m}{k}\biggl\{\sum_{p=0}^a 
(-1)^{a-p}\binom{m-k}{m-p-k}\binom{m-p-k-1}{a-p}\biggr\}
(x_{(k)}y)_{(m-k-a-1)}=
$$
(substituting $r=p-a,\ s=m-a-k-1\geq 0$)
$$
=\sum_{k=0}^{m-a-1}\binom{m}{k}\biggl\{\sum_{r=0}^a (-1)^r
\binom{a+s+1}{r+s+1}\binom{r+s}{r}\biggr\}
(x_{(k)}y)_{(m-k-a-1)}
$$
By Lemma 8.13, this is equal to the first sum in the rhs of (8.14.2). 
 
The second term in the rhs of (8.14.1) is equal to 
$$
x_{(m)}\dpar^{(a)}y=\sum_{p=0}^{min(m,a)}\binom{m}{p}\dpar^{(a-p)}
(x_{(m-p)}y)=
$$
(substituting $k=m-p$)
$$
=\sum_{k\geq max(m-a,0)}\binom{m}{m-k}\dpar^{(a-m+k)}(x_{(k)}y)=
\sum_{k\geq max(m-a,0)}\binom{m}{k}(x_{(k)}y)_{(-1-a+m-k)},
$$
which is the same as the second sum in (8.14.2). 

Therefore, (8.14.1)$=$(8.14.2), which finishes the proof of (8.11.1) in the 
{\it Case B}. 

{\bf 8.15.} {\it Case C. $m<0$ and $n\geq 0$.} First let us treat the case 
$m=-1$. We need to prove that 
$$
[x_{(-1)},y_{(n)}]=\sum_{j\geq 0}\ (-1)^j(x_{(j)}y)_{(n-1-j)}
\eqno{(8.15.1)}
$$
The lhs of (8.15.1) is equal to 
$$
-[y_{(n)},x_{(-1)}]=-\sum_{p\geq 0}\binom{n}{p}(y_{(p)}x)_{(n-p-1)}
\eqno{(8.15.2)}
$$
On the other hand, the rhs of (8.15.1) equals (we use (0.4.3)) 
$$
-\sum_{j\geq 0}\sum_{q\geq 0}(-1)^q\dpar^{(q)}
(y_{(j+q)}x)_{(n-1-j)}\buildrel{(8.4.3)}\over{=}
-\sum_{j\geq 0, q\geq 0}\ \binom{n-1-j}{q}(y_{(j+q)}x)_{(n-1-j-q)}=
$$
(substituting $p=j+q$)
$$
=-\sum_{p\geq 0}(y_{(p)}x)_{(n-1-p)}\cdot\biggl\{\sum_{q=0}^p 
\binom{n-1-p+q}{q}\biggr\}
$$
Therefore, (8.15.1) follows from the identity below, 
which is easily checked by induction on $n$: 
$$
\sum_{q=0}^p\binom{n-1-p+q}{q}=\binom{n}{p}
\eqno{(8.15.3)}
$$
for all $n,p\geq 0$. 

Now assume that $m=-1-a,\ a\geq 0$. We need to prove that 
$$
[x_{(-1-a)},y_{(n)}]=\sum_{p\geq 0}\binom{-1-a}{p}(x_{(p)}y)_{(-1-a+n-p)}
\eqno{(8.15.4)}
$$
The lhs is equal to 
$$
[\dpar^{(a)}x_{(-1)},y_{(n)}]=\sum_{j\geq 0}(-1)^j
(\dpar^{(a)}x_{(j)}y)_{(n-1-j)}=
\sum_{j\geq a}(-1)^{j+a}\binom{j}{a}(x_{(j-a)}y)_{(n-1-j)}=
$$
(substituting $p=j-a$)
$$
=\sum_{p\geq 0}\ (-1)^p\binom{p+a}{a}(x_{(p)}y)_{(n-1-p-a)}
$$
so (8.15.4) follows from (0.0.2). The completes the proof of {\it Case C}. 

{\bf 8.16.} {\it Case D. $m,n\geq 0$.} We have to prove that 
$$
[x_{(m)},y_{(n)}]v=\sum_{j=0}^m\ \binom{m}{j}(x_{(j)}y)_{(m+n-j)}v
\eqno{(8.16.1)}
$$
for all $x,y\in C,\ v\in UC$. In fact, we shall prove a stronger statement, 
which we prefer to formulate as a separate 

{\bf 8.17. Lemma.} {\it If $m$ and $n$ are nonnegative then the OPE identity} 
(8.11.1) {\it holds true on $TC$. In other words,} (8.16.1) {\it 
holds true for all $m,n\geq 0,\ x,y\in C$ and $v\in TC$.}  

We shall prove this by induction on the length 
of $v$. If $v\in C$ then the desired identity holds true by definition 
of a conformal algebra, (0.4.4). Now, let $v=zu$ where $z\in C,\ u\in TC$. 
So, we have to prove that 
$$
[x_{(m)},y_{(n)}]zu=\sum_{j=0}^m\ \binom{m}{j}(x_{(j)}y)_{(m+n-j)}zu
\eqno{(8.17.1)}
$$
Consider the lhs first. We have 
$$
A:=x_{(m)}y_{(n)}zu=x_{(m)}\bigl\{zy_{(n)}u+(y_{(n)}z)u+
\sum_{p=0}^{n-1}\binom{n}{p}(y_{(p)}z)_{(n-p-1)}u\bigr\}=
$$
$$
=zx_{(m)}y_{(n)}u+(x_{(m)}z)y_{(n)}u+\sum_{q=0}^{m-1}\binom{m}{q}
(x_{(q)}z)_{(m-q-1)}y_{(n)}u+
$$
$$
+(y_{(n)}z)x_{(m)}u+(x_{(m)}y_{(n)}z)u+\sum_{q=0}^{m-1}\binom{m}{q}
(x_{(q)}y_{(n)}z)_{(m-q-1)}u+\sum_{p=0}^{n-1}\binom{n}{p}x_{(m)}
(y_{(p)}z)_{(n-p-1)}u
$$
Similarly, 
$$
B:=-y_{(n)}x_{(m)}zu=-zy_{(n)}x_{(m)}u-(y_{(n)}z)x_{(m)}u-
\sum_{p=0}^{n-1}\binom{n}{p}(y_{(p)}z)_{(n-p-1)}x_{(m)}u-
$$
$$
-(x_{(m)}z)y_{(n)}u-(y_{(n)}x_{(m)}z)u-\sum_{p=0}^{n-1}\binom{n}{p}
(y_{(p)}x_{(m)}z)_{(n-p-1)}u-\sum_{q=0}^{m-1}\binom{m}{q}y_{(n)}
(x_{(q)}z)_{(m-q-1)}u
$$
Obviously, the lhs of (8.17.1) is equal to $A+B$, and we should 
compare this with the rhs: 
$$
R:=\sum_{j=0}^m\binom{m}{j}(x_{(j)}y)_{(m+n-j)}zu=\sum_{j=0}^m\binom{m}{j}
\biggl\{z(x_{(j)}y)_{(m+n-j)}u+
$$
$$
+((x_{(j)}y)_{(m+n-j)}z)u+\sum_{r=0}^{m+n-j-1}\binom{m+n-j}{r}
((x_{(j)}y)_{(r)}z)_{(m+n-j-r-1)}u\biggr\}
$$
First of all, $A1+B1=R1$ by induction hypothesis. Next, $A2=-B4$ and 
$A4=-B2$. Further, $A5+B5=R2$ by the axiom of a conformal algebra. 
Next, 
$$
C:=A3+B7=\sum_{q=0}^{m-1}\binom{m}{q}\biggl\{(x_{(q)}z)_{(m-q-1)}y_{(n)}u-
y_{(n)}(x_{(q)}z)_{(m-q-1)}u\biggr\}=
$$
(by induction hypothesis)
$$
=-\sum_{q=0}^{m-1}\binom{m}{q}\sum_{s=0}^{n}\binom{n}{s}
(y_{(s)}x_{(q)}z)_{(m+n-q-s-1)}u
$$
Similarly, 
$$
D:=A7+B3=\sum_{p=0}^{n-1}\binom{n}{p}\sum_{r=0}^m\binom{m}{r}
(x_{(r)}y_{(p)}z)_{(m+n-p-r-1)}u
$$
We have
$$
E:=A6+C_{s=n}=\sum_{q=0}^{m-1}\binom{m}{q}\biggl\{(x_{(q)}y_{(n)}z)_{(m-q-1)}-
(y_{(n)}x_{(q)}z)_{(m-q-1)}\biggr\}u=
$$  
$$
=\sum_{q=0}^{m-1}\sum_{s=0}^q\binom{m}{q}\binom{q}{s}((x_{(s)}y)_{(q+n-s)}z)
_{(m-q-1)}u=
$$
$$
=\sum_{q=0}^{m-1}\sum_{s=0}^q\binom{m}{s}\binom{m-s}{m-q}((x_{(s)}y)_{(q+n-s)}z)
_{(m-q-1)}u
$$
Similarly, 
$$
F:=D_{r=m}+B6=\sum_{p=0}^{n-1}\binom{n}{p}\biggl\{(x_{(m)}y_{(p)}z)_{(n-p-1)}-
(y_{(p)}x_{(m)}z)_{(n-p-1)}\biggr\}u=
$$
$$
=\sum_{p=0}^{n-1}\sum_{r=0}^m\binom{n}{p}\binom{m}{r}((x_{(r)}y)_{(m+p-r)}z)
_{(n-p-1)}u
$$
Finally, 
$$
G:=C_{s<n}+D_{r<m}=\sum_{p=0}^{n-1}\sum_{q=0}^{m-1}\binom{n}{p}
\binom{m}{q}\biggl\{(x_{(q)}y_{(p)}z)_{(m+n-p-q-1)}-
(y_{(p)}x_{(q)}z)_{(m+n-p-q-1)}\bigg\}u=
$$
$$
=\sum_{p=0}^{n-1}\sum_{q=0}^{m-1}\binom{n}{p}\binom{m}{q}
\sum_{s=0}^q\binom{q}{s}\bigl((x_{(s)}y)_{(p+q-s)}z\bigr)_{(m+n-p-q-1)}u=
$$
$$
=\sum_{p=0}^{n-1}\sum_{q=0}^{m-1}\binom{n}{p}
\sum_{s=0}^q\binom{m}{s}\binom{m-s}{m-q}
\bigl((x_{(s)}y)_{(p+q-s)}z\bigr)_{(m+n-p-q-1)}u
$$
A careful examination of the sums involved, together with a formula 
$$
\sum_p\ \binom{a}{r-p}\binom{n}{p}=\binom{a+n}{r}
\eqno{(8.17.2)}
$$
with $a=m-j$, which is obvious from the combinatorial definition 
of the binomial coefficients, shows that $E+F+G=R3$. 

This completes the proof of Lemma. $\btu$

{\bf 8.18. Remark.} Lemmas 8.17 and 8.4.1 mean that $TC$ is canonically 
a {\it module} over conformal algebra $C$. 

{\bf 8.19.} Lemma 8.17 implies {\it Case D}, which completes the proof 
of Lemma 8.10. $\btu$ 

{\bf 8.20.} Now we can finish the proof of Theorem 8.3. According 
to Lemmas 8.5 and 8.10, we have the collection of mutually local 
fields $a(z)=\sum_{n\in\BZ}\ a_{(n)}z^{-n-1},\ a\in C$, acting on the space 
$UC$. These fields satisfy the conditions of Theorem 0.7, which 
gives the desired structure of a vertex algebra on $UC$. 

The other claims of Theorem 8.3 are obvious. $\btu$ 

{\bf 8.21.} Let $C$ be a conformal algebra. 
Let us introduce an operation $[x,y]$ on the space $C$ by  
$$
[x,y]=\sum_{j\geq 0}\ (-1)^j\dpar^{(j+1)}(x_{(j)}y)
\eqno{(8.21.1)}
$$

{\bf 8.22. Theorem.} {\it The operation} (8.21.1) {\it is a Lie bracket 
on $C$.} 

The space $C$ with the Lie algebra structure given by (8.21.1) will be 
denoted $C^{Lie}$.  

The proof is given in 8.23 --- 8.26 below. 

{\bf 8.23.} {\it Skew symmetry.} Using (0.4.3), we have  
$$
[x,y]=\sum_{j\geq 0} (-1)^j\dpar^{(j+1)}\bigl\{(-1)^{j+1}\sum_{p\geq 0}
(-1)^p\dpar^{(p)}(y_{(j+p)}x)\bigr\}=
$$
$$
=-\sum_{j,p\geq 0}(-1)^p\binom{j+1+p}{p}\dpar^{(j+p+1)}(y_{(j+p)}x)=
$$
$$ 
=-\sum_{n\geq 0}\biggl\{\sum_{p=0}^n(-1)^p\binom{n+1}{p}\biggr\}
\dpar^{(n+1)}(y_{(n)}x)=-[y,x],\ QED
$$

{\bf 8.24.} {\it Jacobi identity.} We need to prove that 
$$
[x,[y,z]]=[[x,y],z]+[y,[x,z]]
$$
This is an obvious consequence of two lemmas below. 

{\bf 8.25. Lemma.} {\it For all $x,y,z\in C$} 
$$
[x,[y,z]]=\sum_{a,b\geq 0}\ (-1)^{a+b}\dpar^{(a+b+2)}(x_{(a)}y_{(b)}z)
\eqno{(8.25.1)}
$$

Indeed, we have 
$$
[x,[y,z]]=\sum_{j\geq 0}(-1)^j[x,\dpar^{(j+1)}(y_{(j)}z)]=
\sum_{i,j\geq 0}(-1)^{i+j}\dpar^{(i+1)}\bigl(x_{(i)}\dpar^{(j+1)}(y_{(j)}z)
\bigr)=
$$
(using (0.4.6))
$$
=\sum_{i,j\geq 0}(-1)^{i+j}\dpar^{(i+1)}\biggl\{\sum_{p\geq 0} 
\binom{i}{p}\dpar^{(j+1-p)}(x_{(i-p)}y_{(j)}z)\biggr\}=
$$
$$
=\sum_{i,j\geq 0}(-1)^{i+j}\sum_{p\geq 0}\binom{i+j+2-p}{i+1}
\binom{i}{p}\dpar^{(i+j+2-p)}(x_{(i-p)}y_{(j)}z)=
$$
(substituting $a:=i-p,\ b:=j$) 
$$
=\sum_{a,b\geq 0}(-1)^{a+b}\dpar^{(a+b+2)}(x_{(a)}y_{(b)}z)\cdot
\biggl\{\sum_{p\geq 0}(-1)^p\binom{a+p}{p}\binom{a+b+2}{a+p+1}\biggr\},
$$
and we conclude by Lemma 8.13. 

{\bf 8.26. Lemma.} {\it For all $x,y,z\in C$ 
$$
[[x,y],z]=\sum_{a,b\geq 0}(-1)^{a+b}\dpar^{(a+b+2)}\bigl\{x_{(a)}y_{(b)}z-
y_{(b)}x_{(a)}z\bigr\}
\eqno{(8.26.1)}
$$} 

Indeed, we have 
$$
[[x,y],z]=\sum_{j\geq 0}(-1)^j\sum_{i\geq 0}(-1)^i\dpar^{(i+1)}
\bigl\{\dpar^{(j+1)}(x_{(j)}y)_{(i)}z\bigr\}=
$$
$$
=\sum_{i,j\geq 0}(-1)^{i+1}\binom{i}{j+1}\dpar^{(i+1)}\bigl(
(x_{(j)}y)_{(i-j-1)}z)=
$$
(by (0.4.5))
$$
=\sum_{i,j\geq 0}(-1)^{i+1}\binom{i}{j+1}\sum_{q=0}^j
(-1)^q\binom{j}{q}\dpar^{(i+1)}\bigl\{x_{(j-q)}y_{(i-j+q-1)}z-
(-1)^jy_{(i-q-1)}x_{(q)}z\bigr\},
$$
and we again conclude after an easy rearrangement of the indices and 
using Lemma 8.13. 

This completes the proof of lemma and of the theorem. $\btu$ 

{\bf 8.27. Corollary.} {\it The vertex enevlope $UC$ of a conformal 
algebra $C$ is equal to the associative enveloping algebra $UC^{Lie}$ of the 
Lie algebra $C^{Lie}$.} 

{\it If $i:\ C\lra UC$ denotes the canonical morphism, 
$$
i([x,y])=i(x)_{(-1)}i(y)-i(y)_{(-1)}i(x)
\eqno{(8.27.1)}
$$
for all $x,y\in C$.} 

This follows immediately from the definition in 8.2. 

Note that in an arbitrary vertex algebra, although the operation 
$x_{(-1)}y$ is {\it not} associative, the commutator 
$x_{(-1)}y-y_{(-1)}x$ is a Lie bracket, according to a remark by 
A.Radul cited in [K], 3.1, page 82. 

This corollary supports the point of view advocated by Beilinson 
and Drinfeld, namely, that conformal algebras ({\it $Lie^*$-algebras} in 
the language of BD) are analogs of Lie algebras, while 
vertex algebras ({\it chiral algebras}) are analogs of associative algebras.   

{\bf 8.28.} Let us define a canonical increasing filtration 
$F_\cdot UC$ on $UC$ by induction: $F_{-1}UC=0,\ F_0UC=k\cdot\b1,\ 
F_{i+1}UC=F_iUC+C_{(-1)}F_iUC$. 

Obviously, the operation $x_{(-1)}y$ induces on the associated graded 
space $gr_F(UC):=\oplus_{i\geq 0}\ gr^i_F(UC):=\oplus_i\ F_iUC/F_{i-1}UC$ 
a structure of a commutative associative unitary $k$-algebra. 

We have an evident surjective morphism of commutative algebras 
$$
Sym_k(C)\lra gr_F(UC)
\eqno{(8.28.1)}
$$

{\bf 8.29. Theorem.} {\it If $C$ is a free $k$-module then} (8.28.1) 
{\it is an isomorphism.} 

This is immediate consequence of Corollary 8.27 and the usual 
Poincar\'e-Birkhoff-Witt theorem for Lie algebras, cf. [Se], 
Part I, Chapter III, Theorem 4.3. 

\bigskip\bigskip

\centerline{\bf \S 9. Enveloping Algebra of a Vertex Algebroid}

\bigskip\bigskip

{\it Abelian case}

\bigskip

{\bf 9.1.} Let us define a {\it $\dpar$-module} (over $k$) to be a 
$\BZ_{\geq 0}$-graded $k$-module $C=\oplus_{i\geq 0}\ C_i$ equipped 
with a family 
of endomorphisms $\dpar^{(j)}:\ C\lra C,\ j\geq 0,$ where $\dpar^{(j)}$ 
has degree $j$, $\dpar^{(0)}=Id$ and 
$$
\dpar^{(i)}\dpar^{(j)}=\binom{i+j}{i}\dpar^{(i+j)}
$$
Obviously, a $\dpar$-module is the same as an {\it abelian} conformal 
algebra, i.e. a conformal algebra in which all operations $_{(n)}$ are 
trivial. These objects, with obvious morphisms, form a category 
$\dpar-\CM od$. 

Let us define a {\it $1$-truncated $\dpar$-module} to be a triple 
$\CX=(C_0,C_1,\dpar)$ 
where $C_i$ are $k$-modules and $\dpar:\ C_0\lra C_1$ is a map 
of $k$-modules. These objects form a category $\dpar-\CM od_{\leq 1}$. 
We have an obvious truncation functor 
$$
\dpar-\CM od\lra\dpar-\CM od_{\leq 1}
\eqno{(9.1.1)}
$$
It is easy to construct the left adjoint to (9.1.1). Namely, 
given $\CX=(C_0,C_1,\dpar)\in\dpar-\CM od_{\geq 1}$, consider the direct sum 
$$
M=C_0\oplus C_1\oplus(\oplus_{i\geq 2}\ \dpar^{(i)}C_0)
\oplus(\oplus_{i\geq 1}\ \dpar^{(i)}C_1)
\eqno{(9.1.2)}
$$
where $\dpar^{(i)}X$ denotes a copy of a $k$-module $X$, whose elements are 
denoted $\dpar^{(i)}x,\ x\in X$. Define a $\dpar$-module 
$\CX^{\dpar}=\oplus_{i\geq 0}\ C_i$ as 
the quotient of $M$ modulo the submodule generated by the following 
elements:
$$\
\dpar^{(i+1)}x-(i+1)\dpar^{(i)}\dpar(x),\ i\geq 1, x\in C_0
\eqno{(9.1.3)}
$$
The grading and operators $\dpar^{(j)}:\ C_i\lra C_{i+j}$ 
are defined in the obvious manner. 

Note that if the ground ring $k$ contains the field of rational numbers 
$\BQ$ then $\CX$ is equal simply to 
$$
\CX^\dpar=C_0\oplus C_1\oplus(\oplus_{i\geq 1}\ \dpar^{(i)}C_1)
\eqno{(9.1.4)}
$$ 

{\bf 9.2.} Recall that in 0.9 one introduced the category 
$\dpar-\CA lg$ of $\dpar$-algebras which are identified with 
commutative vertex algebras. A {\it $\dpar$-ideal} $J$ in a $\dpar$-algebra $B$ 
is a $\BZ_{\geq 0}$-graded ideal stable under all endomorphisms 
$\dpar^{(i)}$; this is the same as a vertex ideal if $B$ is 
understood as a vertex algebra.  

We have an obvious forgetful functor 
$$
\dpar-\CA lg\lra\dpar-\CM od
\eqno{(9.2.1)}
$$
This functor admits a left adjoint
$$
Sym:\ \dpar-\CM od\lra\dpar-\CA lg
\eqno{(9.2.2)} 
$$
that assigns to a $\dpar$-module $C$ 
the symmetric algebra (over $k$) $Sym(C)$. It is nothing but the 
restriction of the functor defined in 8.3 to the full subcategory 
of abelian conformal algebras. 

{\bf 9.3. Jet algebra. } Let us call a vertex algebroid 
$\CA=(A,T,\Omega,\ldots)$ 
{\it abelian} if $T=0$. Thus, an abelian vertex algebroid is 
simply a triple $\CA=(A,\Omega,\dpar)$ where $A$ is a commutative algebra, 
$\Omega$ is an $A$-module and $\dpar:\ A\lra\Omega$ is a derivation. 
Thus, abelian vertex algebroids are the same as "1-truncated $\dpar$-algebras". 
The category of abelian vertex algebroids will be denoted $\CA lg\CA b$; 
it is a full subcategory of $\CA lg$.  

We have an obvious forgetful functor 
$$
o:\ \CA lg\CA b\lra\dpar-\CM od
\eqno{(9.3.1)}
$$
and the truncation functor
$$
\dpar-\CA lg\lra \CA lg\CA b
\eqno{(9.3.2)}
$$
Let us construct a left adjoint to (9.3.2). Given $\CA=(A,\Omega,\dpar)\in 
\CA lg\CA b$, consider the $\dpar$-module $(o\CA)^{\dpar}$ and its 
symmetric algebra $S(\CA):=Sym((o\CA)^{\dpar})\in\dpar-\CA lg$. Note that 
the grading on $S(\CA)$ is defined in such a way that $A\subset S(\CA)$ has 
grading $0$ and $\Omega\subset S(\CA)$ has grading $1$. 

Let $R\subset S(\CA)$ be a $\dpar$-ideal generated by all elements 
$$
1_A-1_{S(\CA)};\ ab-a\cdot b;\ a\omega-a\cdot\omega, 
\eqno{(9.3.3)}
$$
$a,b\in A,\ \omega\in\Omega$ (note that these elements are indeed homogeneous). 
We denote the quotient $\dpar$-algebra $S(\CA)/R$ by $J(\CA)$, 
and call it the {\it jet algebra} of $\CA$. The assignment $\CA\mapsto J(\CA)$ 
defines a functor 
$$
J:\ \CA lg\CA b\lra\dpar-\CA lg
\eqno{(9.3.4)}
$$
left adjoint to (9.3.2). 

The composition $A\hra (o\CA)^\dpar=Sym^1((o\CA)^\dpar)\hra S(\CA)\lra J(\CA)$ 
defines a map 
$$
i_A:\ A\lra J(\CA)
\eqno{(9.3.5)}
$$
which is a map of $k$-algebras, due to the relations (9.3.3), which 
identifies $A$ with $J(A)_0$. This map makes $J(\CA)$ an 
$A$-algebra.  

A similar composition defines a map 
$$
i_\Omega:\ \Omega\lra J(\CA)
\eqno{(9.3.6)}
$$
which is a map of left $A$-modules, again due to the relations (9.3.3). 
We have the compatibility: $\dpar(i_A(a))=i_\Omega(\dpar(a))$, and the triple 
$(J(\CA),i_A,i_\Omega)$ has a  

{\bf 9.3.1.} {\bf Universal Property.}  {\it Given a triple $(B,i'_A,
i'_\Omega)$ where $B$ is a $\dpar$-algebra, $i'_A:\ A\lra B$ a morphism 
of $k$-algebras such that $i'_A(A)\subset B_0$, 
$i'_\Omega:\ \Omega\lra B$ a morphism of $A$-modules such that 
$i'_\Omega(\Omega)\subset B_1$ and $i'_\Omega(\dpar(a))=\dpar(i'_A(a))$, 
there exists a unique map of $\dpar$-algebras $f:\ J(A)\lra B$ such that 
$f\circ i_A=i'_A$ and $f\circ i_\Omega=i'_\Omega$.}    

{\bf 9.4.} 
Let $J^+=\oplus_{i\geq 1}\ J(\CA)_i$ be the augmentation ideal. Consider the 
associated graded algebra with respect to the $J^+$-adic filtration
$$
gr J(\CA)=\oplus_{i\geq 0}\ J^{+i}/J^{+i+1}
\eqno{(9.4.1)}
$$
It inherits a $\BZ_{\geq 0}$-grading from $J(\CA)$. Let $\Omega^{(i)}$ denote 
the image of the composition 
$$
\dpar^{(i)}A\oplus \dpar^{(i-1)}\Omega\lra J^+\lra J^+/J^{+2}
$$
Note that $\Omega^{(i)}$ is an $A$-submodule of $J^+/J^{+2}$, and 
is contained inside the homogeneous component $(J^+/J^{+2})_i$. 

Adding up, we get a map of $A$-modules 
$$
\oplus_{i\geq 1}\ \Omega^{(i)}\lra J^+/J^{+2}
\eqno{(9.4.2)}
$$
and hence a morphism of $A$-algebras 
$$
Sym_A(\oplus_{i\geq 1}\ \Omega^{(i)})\lra Sym_A(J^+/J^{+2})\lra gr J(\CA)
\eqno{(9.4.3)}
$$
the second arrow being the evident canonical map. 

{\bf 9.5. Theorem.} {\it Assume that $k\supset\BQ$ and $\Omega$ is a free 
$A$-module. Then} 

(i) {\it the maps $\Omega\lra \Omega^{(i)}$ sending $\omega$ to (the image 
of) $\dpar^{(i-1)}\omega,\ i\geq 1,$ are isomorphisms of $A$-modules}; 

(ii) {\it both maps in} (9.4.3) {\it are isomorphisms of $A$-algebras.} 

Indeed, choose an $A$-base $\{\omega_r\}$ of $\Omega$. Let $\tOmega^{(i)}$ 
denote a free $A$-module with the base $\{\dpar^{(i-1)}\omega_r\}$. 
Let $\tJ$ be the symmetric algebra 
$Sym_A\{A\oplus(\oplus_{i\geq 1}\ \tOmega^{(i)})\}$. It has an obvious 
structure of a $\dpar$-algebra, and satisfies the universal property 
9.3.1. Therefore, we have canonical isomorphism 
$$
J(\CA)\iso\tJ
\eqno{(9.5.1)}
$$
We can apply the construction of 9.4 to the algebra $\tJ$ as well, and for it 
the claims of Theorem 9.5 are clear. On the other hand, the maps 
$\Omega^{(i)}\lra\tOmega^{(i)}$ induced by (9.5.1) are isomorphisms. 
This implies the theorem. $\btu$  

\bigskip

{\it Conformal algebra associated with a vertex algebroid}

\bigskip

{\bf 9.6.} Following the pattern of 3.1, let us define a 
{\it $1$-truncated conformal algebra} to be a quintuple 
$c=(C_0,C_1,\dpar,_{(0)},_{(1)})$ where $C_i$ are $k$-modules, 
$\dpar:\ C_0\lra C_1$ is a $k$-module map and 
$_{(i)}:\ (C_0\oplus C_1)^{\otimes 2}\lra (C_0\oplus C_1)$ are  
bilinear operations of degree $-i-1$. 

Elements of $C_0$ (resp. of $C_1$) will be denoted $a,b,c$ (resp. $x,y,z$). 

These data must satisfy the following axioms. 

{\it (Derivation)} 
$$
(\dpar a)_{(0)}=0;\ (\dpar a)_{(1)}=-a_{(0)};\ 
\dpar(x_{(0)}a)=x_{(0)}\dpar a
\eqno{(Der)}
$$

{\it (Commutativity)} 
$$
x_{(0)}a=-a_{(0)}x;\ x_{(0)}y=-y_{(0)}x+\dpar(y_{(1)}x);\ 
x_{(1)}y=y_{(1)}x
\eqno{(Com)}
$$

{\it (Associativity)} {\it For all $\alpha,\beta,\gamma\in C_0\oplus C_1$} 
and $i=0,1$, 
$$
\alpha_{(0)}\beta_{(i)}\gamma=(\alpha_{(0)}\beta)_{(i)}\gamma+
\beta_{(i)}\alpha_{(0)}\gamma
\eqno{(Ass)}
$$

$1$-truncated conformal algebras, with obvious morphisms, 
form a category $\CC onf_{\leq 1}$. 

We have an evident forgetful functor 
$$
\CV ert_{\leq 1}\lra \CC onf_{\leq 1}
\eqno{(9.6.1)}
$$

{\bf 9.7.} We have an obvious truncation functor 
$$
\CC onf\lra \CC onf_{\leq 1}
\eqno{(9.7.1)}
$$
which assigns to a conformal algebra $C$ its part $C_{\leq 1}$ of 
degree $\leq 1$. 
Let us construct a left adjoint to this functor. 

Given $c=(C_0,C_1,\ldots)\in\CC onf_{\leq 1}$, consider $c$ as a 
$1$-truncated $\dpar$-module (forgetting the operations), 
and consider the corresponding $\dpar$-module $C:=c^\dpar$, cf. 9.1. 

{\bf 9.8. Theorem.} {\it There is a unique structure of a conformal 
algebra on $C$ such that operations $_{(i)}$ and $\dpar$ on the subspace 
$C_{\leq 1}$ coincide with the ones given by the structure of 
a $1$-truncated conformal algebra on $c$. 

The assignement $c\mapsto C$ gives a functor 
$$
\CC onf_{\leq 1}\lra \CC onf
\eqno{(9.8.1)}
$$
left adjoint to} (9.7.1).    

Indeed, uniqueness is clear, due to the axioms (0.4.2) and (0.4.3) 
of a conformal algebra. We leave the proof of existence to the reader. 
The claim about adjointness is evident. 

{\bf 9.9.} Let $\CA=(A,T,\Omega,\ldots)$ be a vertex algebroid. We assign 
to $\CA$ a $1$-truncated conformal algebra $c\CA=(C_0,C_1,\ldots)$ by setting 
$C_0=A,\ C_1=T\oplus\Omega,\ \dpar:\ C_0\lra C_1$ to be the composition 
of $\dpar:\ A\lra\Omega$ with the embedding $\Omega\subset C_1$, and 
defining the operations $_{(i)}$ by
$$
a_{(0)}b=a_{(0)}\omega=\omega_{(0)}\omega'=0;\ 
\tau_{(0)}a=\tau(a);\ \tau_{(0)}\omega=\tau(\omega)
\eqno{(9.9.1)}
$$
$$
\tau_{(0)}\tau'=[\tau,\tau']-c(\tau,\tau')+\frac{1}{2}\dpar\langle\tau,
\tau'\rangle
\eqno{(9.9.2)}
$$
$$
x_{(1)}y=\langle x,y\rangle
\eqno{(9.9.3)}
$$
In other words, $c\CA$ is the result of application of the forgetful 
functor (9.6.1) to the $1$-truncated vertex algebra $u\CA$ defined in 3.3. 

{\bf 9.10.} We may apply the functor 
(9.8.1) to the $1$-truncated conformal algebra $c\CA$, 
and get a conformal algebra, to be denoted $C\CA$. 

Assume for simplicity that $k\supset\BQ$.  
Then, as a $k$-module, 
$$
C\CA=A\oplus\ (\oplus_{i\geq 0}\ \dpar^{(i)}T)\ \oplus\ (\oplus_{i\geq 0}\ 
\dpar^{(i)}\Omega)
\eqno{(9.10.1)}
$$
Here are the explicit formulas for the operations $_{(n)}$ 
(below we agree that $\dpar^{(i)}\alpha=0$ if $i<0$). First of all,   
$\alpha_{(n)}\beta=0$ for all $\alpha, \beta\in A\oplus(\oplus_{i\geq 0} 
\dpar^{(i)}\Omega)$ and $n\geq 0$. Next, 

$$
a_{(n)}\dpar^{(i)}\tau=-\dpar^{(i-n)}\tau(a);\ 
\tau_{(n)}\dpar^{(i)}a=\dpar^{(i-n)}\tau(a)
\eqno{(9.10.2)}
$$
$$
\omega_{(n)}\dpar^{(i)}\tau=-\dpar^{(i-n)}\tau(\omega)+
(i+1)\dpar^{(i-n+1)}\langle\omega,\tau\rangle
\eqno{(9.10.3)}
$$
$$
\tau_{(n)}\dpar^{(i)}\omega=\dpar^{(i-n)}\tau(\omega)+
n\dpar^{(i-n+1)}\langle\tau,\omega\rangle
\eqno{(9.10.4)}
$$
Finally, 
$$
\tau_{(n)}\dpar^{(i)}\tau'=\dpar^{(i-n)}\bigl\{[\tau,\tau']-c(\tau,\tau')
\bigr\}+\frac{1}{2}(i+n+1)\dpar^{(i-n+1)}\langle\tau,\tau'\rangle
\eqno{(9.10.5)}
$$

\bigskip

{\it Construction of envelope}

\bigskip

{\bf 9.11.} Now we can construct the vertex envelope of a vertex 
algebroid $\CA=(A,T,\Omega,\ldots)$. Consider the  
conformal algebra $C\CA$; our vertex envelope $U\CA$ will be 
a quotient of the vertex envelope $UC\CA$ defined in the previous 
Section by certain vertex ideal.  

To define this ideal, let us return to the tensor algebra $T:=TC\CA$, 
as in 8.2. Let us define a {\it left} $\dpar$-ideal $R\subset T$ generated 
by the following elements: 
$$
r_0(u):=(1_A-1_T)u,\ u\in T
\eqno{(9.11.1)}
$$
$$
r(a,x):=a\cdot x-ax,\ a\in A, x\in A\oplus\Omega;\ 
r(a,\tau)=a\cdot\tau-a\tau+\gamma(a,\tau)
\eqno{(9.11.2)}
$$
Here by {\it $\dpar$-ideal} we mean an ideal closed under all operators 
$\dpar^{(i)}$. 

{\bf 9.12. Lemma.} {\it The ideal $R$ is respected by all operations 
$y_{(n)},\ y\in C\CA,\ n\in\BZ$.} 

Indeed, one has to check this for $y\in A\oplus T\oplus\Omega$, $n\geq 0$, 
and only on generators of $R$. This is done by an easy case-by-case computation 
that we will leave to the reader, restricting ourselves by an 
example. 

Let us check that $\tau_{(0)}r(a,\tau')\in S$. We have 
$$
\tau_{(0)}r(a,\tau')=\tau_{(0)}\bigl\{a\cdot\tau'-a\tau'+\gamma(a,\tau')
\bigr\};
$$
$$
\tau_{(0)}(a\cdot\tau')=\tau_{(0)}a\cdot\tau'+a\cdot\tau_{(0)}\tau'=
\tau(a)\cdot\tau'+a\cdot([\tau,\tau']-c(\tau,\tau')+\frac{1}{2}\dpar
\langle\tau,\tau'\rangle);
$$
$$
-\tau_{(0)}(a\tau')=-[\tau,a\tau']+c(\tau,a\tau')-\frac{1}{2}\dpar\langle
\tau,a\tau'\rangle
$$
We have 
$$
-[\tau,a\tau']=-a[\tau,\tau']-\tau(a)\tau'; 
$$
by axiom (A3) of a vertex algebroid, 
$$
c(\tau,a\tau')=-c(a\tau',\tau)=-ac(\tau',\tau)-\gamma(a,[\tau',\tau])+
\gamma(\tau(a),\tau')-\tau(\gamma(a,\tau'))+
$$
$$
+\frac{1}{2}\langle\tau',\tau\rangle\dpar a-\frac{1}{2}\dpar\tau'\tau(a)+
\frac{1}{2}\dpar\langle\tau,\gamma(a,\tau')\rangle
$$
By axiom (A2), 
$$
-\frac{1}{2}\dpar\langle\tau,a\tau'\rangle=-\frac{1}{2}\dpar\bigl\{
a\langle\tau',\tau\rangle+\langle\gamma(a,\tau'),\tau\rangle-\tau'\tau(a)
\bigr\}=
$$
$$
=-\frac{1}{2}\langle\tau',\tau\rangle\dpar a-
\frac{1}{2}a\dpar\langle\tau',\tau\rangle-
\frac{1}{2}\dpar\langle\gamma(a,\tau'),\tau\rangle+
\frac{1}{2}\dpar\tau'\tau(a)
$$
Finally, 
$$
\tau_{(0)}\gamma(a,\tau')=\tau(\gamma(a,\tau'))
$$   
After bookkeeping, we see that
$$
\tau_{(0)}r(a,\tau')=r(a,[\tau,\tau'])+r(\tau(a),\tau')-
r(a,c(\tau,\tau'))+\frac{1}{2}r(a,\dpar\langle\tau,\tau'\rangle)\in S,\ 
QED
$$

{\bf 9.13.} Let $\bR\subset UC\CA$ denote the image of $R$ under 
the canonnical projection $TC\CA\lra UC\CA$. The previous lemma 
means that $\bR$ is a vertex ideal in $UC\CA$. Let us denote by 
$U\CA$ the vertex algebra $UC\CA/\bR$. It is equipped with an obvious 
splitting $T\hra T\oplus\Omega=U\CA_1$. 

Conversely, given a splittable vertex algebra $V$, choose a splitting 
$s:\ T\lra V_1$, and consider the vertex algebroid $\CA(V;s)$, cf. \S 2. 
Let $\CV ert'\subset\CV ert$ denote the full subcategory of splittable 
vertex algebras. We get a functor 
$$
\CA:\ \CV ert'\lra\CA lg
\eqno{(9.13.1)}
$$
which is in fact the composition of the truncation functor (3.1.1) 
(restricted to $\CV ert'$) and of the functor quasiinverse to (3.3.6). 

{\bf 9.14. Theorem.} {\it The assignement $\CA\mapsto U\CA$ provides a functor 
$$
U:\ \CA lg\lra\CV ert'
\eqno{(9.14.1)}
$$
which is left adjoint to} (9.13.1). 

Indeed, this is evident from the construction. 

\bigskip


{\it Poincar\'e-Birkhoff-Witt}

\bigskip

{\bf 9.15.} Let $\CA=(A,T,\Omega,\ldots)$ be a vertex algebroid. 
The enveloping algebra $U\CA$ is generated as a $k$-module by the monomials 
of the form 
$$
x^p\cdot x^{p-1}\cdot\ldots\cdot x^1
\eqno{(9.15.1)}
$$
where $x^i$ has the form $\dpar^{(j)}y,\ y=a,\omega$ or $\tau$. 
Here we have denoted for brevity by $x\cdot x'$ the operation $x_{(-1)}x'$. 

Let us introduce a canonical increasing exhaustive filtration 
$F_0U\CA\subset F_1U\CA\subset\ldots U\CA$ by setting $F_iU\CA$ to be 
equal the $k$-submodule of $U\CA$ generated by all monomials (9.15.1) 
where there are $\leq i$ letters $x_j$ of the form $\dpar^{(a)}\tau$. 

Obviously all submodules $F_iU\CA$ are stable under $\dpar^{(j)}$ and 
$F_iU\CA\cdot F_jU\CA\subset F_{i+j}U\CA$. Consider the 
associated graded module $gr_FU\CA=\oplus_{i\geq 0}\ F_iU\CA/F_{i-1}U\CA$. 

It is easy to see that the operation $x\cdot y$ induces a {\it 
commutative and associative} multiplication on $gr_FU\CA$, i.e. 
$gr_FU\CA$ becomes a $\dpar$-algebra (an abelian vertex algebra). 

Let $\CA^{ab}$ denote an abelian vertex algebroid $(A,\Omega\oplus T,\dpar)$. 
We have 
$$
\CA(gr_FU\CA)=\CA^{ab}
\eqno{(9.15.2)}
$$
hence by adjunction we get a canonical map of $\dpar$-algebras
$$
J(\CA^{ab})\lra gr_FU\CA
\eqno{(9.15.3)}
$$
It is clear the this map is surjective. 

Recall that by 9.4 we have a canonical filtration on $J(\CA^{ab})$ 
(let us denote it by $G$ here) and 
a canonical surjective map of $\dpar$-algebras 
$$
Sym_A\bigl\{(\oplus_{i\geq 1}\ T^{(i)})\oplus 
(\oplus_{i\geq 1}\ \Omega^{(i)})\bigr\}\lra gr_GJ(\CA^{ab})
\eqno{(9.15.4)}
$$
where $X^{(i)}\ (X=T$ or $\Omega$) denotes a copy of $X$ sitting in weight 
$i$. 

By refining the filtration $F$ using the filtration $G$, 
we get a canonical filtration $H$ on $U\CA$ together with a surjective 
map (of $\dpar$-algebras over $A$)
$$
Sym_A\bigl\{(\oplus_{i\geq 1}\ T^{(i)})\oplus 
(\oplus_{i\geq 1}\ \Omega^{(i)})\bigr\}\lra gr_HU\CA
\eqno{(9.15.5)}
$$
The multiplication in the right hand side is induced by the operation 
$_{(-1)}$ on $U\CA$. 

{\bf 9.16. Theorem.} {\it Assume that $k\supset\BQ$ and both $\Omega$ 
and $T$ are free $A$-modules. Then the maps} (9.15.3) {\it and} (9.15.5) 
{\it are isomorphisms. 

Filtration $H$ is compatible with the 
conformal grading, finite on each component with the fixed 
conformal weight and is canonical in following sense: for any morphism 
of vertex algebroids $\CA\lra\CA'$, the corresponding morphism 
of vertex envelopes $U\CA\lra U\CA'$ respects these filtrations.} 

To prove this, we use the strategy of Serre's proof of the usual PBW theorem, 
cf. [Se], pp. 14-16. Let us choose well ordered $A$-bases 
$\{\tau_s\},\ \{\omega_r\}$ of 
$T$ and $\Omega$. They give rise to a well ordered $A$-base 
$$
\{x_u\}_{u\in I}=\{\dpar^{(a)}\tau_s, \dpar^{(b)}\omega_r\}
$$ 
of $(\oplus_{i\geq 1}\ T^{(i)})\oplus (\oplus_{i\geq 1}\ \Omega^{(i)})$. 
For a sequence $M=(u_1,\ldots,u_m),\ u_j\in I,\ u_1\leq u_2\leq\ldots
\leq u_m$ define a monomial $x_M\in U\CA$ by 
$$
x_M=x_{u_1}\cdot\ldots\cdot x_{u_m}
$$
Similarly, $M$ defines monomials $\bx_M\in gr_HU\CA$ and 
$$
\tx_M\in 
S:=Sym_A\bigl\{(\oplus_{i\geq 1}\ T^{(i)})\oplus (\oplus_{i\geq 1}\ 
\Omega^{(i)})\bigr\}
$$
Obviously, the monomials $\{\tx_M\}$ form an $A$-base of $S$. 
On the other hand, it is easy to see that each element of $U\CA$ may 
by written as $\sum a_M\cdot x_M$ with some $a_M\in A$. To prove 
our claim it is enough to show that a relation 
$$
\sum a_M\cdot x_M=0
\eqno{(9.16.1)}
$$ 
implies that all $a_M=0$ (cf. {\it loc. cit.}, 
Lemma 4.5).   
   
Proceeding in a manner similar to 8.4, we define an action of the 
Lie algebra $C\CA^{Lie}$, and hence of its envelope $UC\CA^{Lie}=UC\CA$,  
on $S$. Next one checks the relations (9.11.1) and (9.11.2) and therefore 
gets an action of $U\CA$ on $S$. 

One sees immediately from the definitions that $(a\cdot x_M)\cdot 1_S=a\tx_M$. 
Hence (9.16.1) implies the relation $\sum a_M\tx_M=0$ in $S$ and 
therefore $a_M=0$, which proves the theorem.     
We leave the details to the reader. $\btu$     
     
{\bf 9.17.} Let $X$ be a smooth $k$-scheme, $U\subset X$ a Zariski open 
in $X$. A vertex algebra of the form $U\CA$ where $\CA\in\CA lg_{\CT_X}(U)$ 
is a section of the gerbe $\CA lg_{\CT_X}$ discussed in Section 7, 
is called an {\it algebra of chiral differential operators} on $U$. 

These algebras form a gerbe $\CD iff_X$, by definition isomorphic to the 
gerbe $\CA lg_X$. 

Note that all isomorphisms of algebras of chiral differential operators 
respect the canonical filtrations on them, and there is {\it no} 
obstruction to the gluing of associated graded algebras. 

In particular, we have 

{\bf 9.18. Theorem.} {\it Assume that $k\supset\BQ$. 
Each algebra of chiral do $\CD_X\in\CD iff_X(X)$ 
admits a canonical filtration whose graded algebra is isomorphic to 
$$
gr(\CD_X)=Sym_{\CO_X}\bigl\{(\oplus_{i\geq 1}\ \CT_X^{(i)})
\oplus(\oplus_{i\geq 1}\ \Omega_X^{1(i)})\bigr\}
\eqno{(9.18.1)}
$$
where $\CT_X^{(i)}$ (resp. $\Omega_X^{1(i)}$) denotes a copy of 
the tangent bundle (resp. of the bundle of $1$-forms) sitting in 
conformal degree $i$.} $\btu$ 

\bigskip\bigskip

\bigskip\bigskip

\newpage

\centerline{\bf References}

\bigskip\bigskip

[BB] A.~Beilinson, J.~Bernstein, A proof of Jantzen conjectures, 
I.M.Gelfand Seminar, {\it Adv. in Soviet Math.} {\bf 16}, Part 1, 
AMS, Providence, RI, 1993, pp. 1-50.  

[BD1] A.~Beilinson, V.~Drinfeld, Excerpt from "Chiral Algebras II" on 
chiral tdo and Tate structures, Preprint, 1999. 

[BD2] A.~Beilinson, V.~Drinfeld, Chiral algebras, Preprint, Version 2000.  

[B] R.~Borcherds, Vertex algebras, Kac-Moody Lie algebras, and the 
Monster, {\it Proc. Nat. Acad. Sci. USA}, {\bf 83} (1986), no. 10, 
3068-3071.  

[GMS] V.~Gorbounov, F.~Malikov, V.~Schechtman, Gerbes of chiral differential 
operators, math.AG/9906117; {\it Math. Research Letters}, {\bf 7} (2000), 
55-66. 

[Gr] A.~Grothendieck, Classes de Chern et r\'epresentations lin\'eaires 
des groupes discrets, pp. 215-305 in: Dix expos\'es sur la cohomologie 
des sch\'emas, Masson \& Cie, Paris --- North Holland, Amsterdam, 1968.    

[H] B.~Harris, Group cohomology classes with differential form 
coefficients, 
\newline Algebraic $K$-theory, Evanston 1976, {\it 
Lect. Notes in Math.} {\bf 551}, Springer, 1976; pp. 278-282.   

[K] V.~Kac, Vertex algebras for beginners, 
Second Edition, University Lecture Series, {\bf 10}, 
American Mathematical Society, Providence, Rhode Island, 1998. 

[MS] F.~Malikov, V.~Schechtman, Chiral de Rham complex. II, 
D.B.Fuchs' 60-th Anniversary volume, Differential topology, 
infinite-dimensional Lie algebras, and applications, 149-188, 
{\it AMS Transl. Ser. 2}, {\bf 194}, AMS, Providence, RI, 1999. 
   
[MSV] F.~Malikov, V.~Schechtman, A.~Vaintrob, Chiral de Rham complex,\ 
{\it Comm. Math. Phys.}, {\bf 204} (1999), 439-473. 


[S] V.~Schechtman, Riemann-Roch theorem and the Atiyah-Hirzebruch spectral 
sequence, {\it Usp. Mat. Nauk}, {\bf 35}, no. 6 (1980), 179-180 (Russian). 

[Se] J.-P.~Serre, Lie algebras and Lie groups, 1964 lectures given 
at Harvard University, Second Edition, {\it Lect. Notes in Math.} 
{\bf 1500}, Springer-Verlag, Berlin, 1992.

\bigskip

\bigskip

V.G.: Department of Mathematics, University of Kentucky, 
Lexington, KY 40506, USA;\ vgorb\@ms.uky.edu

F.M.: Department of Mathematics, University of Southern California, 
Los Angeles, CA 90089, USA;\ fmalikov\@mathj.usc.edu   

V.S.: IHES, 35 Route de Chartres, 91440 Bures-sur-Yvette, France;\ 
vadik\@ihes.fr

\enddocument